\documentclass[dvips]{mrv9x6}

\usepackage{makeidx,wsminitoc}
\usepackage{amsxtra}% defines \sphat and \spcheck

\makeindex

\makeatletter  %
\nomtcrule
\def\mtc@bottom@rule{\vskip 5pt \hrule height \z@ width \columnwidth \kern2.6\p@}
\renewcommand*\l@section{\@dottedtocline{1}{1em}{1.3em}}
\renewcommand*\l@subsection{\@dottedtocline{2}{2.3em}{2em}}
\renewcommand*\l@subsubsection{\@dottedtocline{3}{4.3em}{2.6em}} 

\makeatother

\newlength{\displayboxwidth}
\newlength{\hangboxwidth}

\newtheorem{scholium}{Scholium}[subsubsection]
\newlength{\qedskip}
\numberwithin{equation}{subsection}

\begin{document}

\dominitoc       %GENERATE OUT CONTENT TEXT AT PER CHAPTER ONLY

%As requested from old series:-
\makeatletter
\renewcommand\thedefinition{\thesection.\arabic{definition}}
\newtheorem{defn}{Definition}[subsubsection]
\newtheorem{theo}{Theorem}[subsubsection]
\newtheorem{lem}{Lemma}[subsubsection]
\newtheorem{exam}{Example}[subsubsection]
\newtheorem{cor}{Corollary}[subsubsection]
\newtheorem{prop}{Proposition}[subsubsection]
\newtheorem{rem}{Remark}[subsubsection]
\def\x{{\bf x}}
\def\y{{\bf y}}
\newcommand{\X}{{\cal X}}
\newcommand{\F}{\mbox{\bf F}}
\newcommand{\R}{\mbox{\bf R}}
\def\L{{\cal L}}
\def\FX{{\F_q(\X)}}
\def\wt{{\rm wt}}
\def\nP{\nu_P}
\makeatother

%%%%%%%%%%%%%%%%%%%%%% pls. do not remove these lines %%%%%%%%%%%%%%%%%%%%%
%\titlepages %pls. leave this as it is. Pages 1-4 are reserved for publisher
%%%%%%%%%%%%%%%%%%%%%%%%%%%%%%%%%%%%%%%%%%%%%%%%%%%%%%%%%%%%%%%%%%%%%%%%%%%

%\setcounter{page}{5}

%%%%%%%%%%%%%%%%%%%%%% Need not the change anything %%%%%%%%%%%%%%%%%%%%%
%\begin{tblofcontents}
%%%%% Need not to input any data here.
%%%%% System will input automatically.
%\end{tblofcontents}
%%%%%%%%%%%%%%%%%%%%%% Need not the change anything %%%%%%%%%%%%%%%%%%%%%

%\begin{foreword}
%\input foreword.tex
%\end{foreword}

%\begin{preface}
%\input preface.tex
%\end{preface}

\setcounter{page}{1}

\chapter[Unitary Matrix Functions, Algorithms, Wavelets]{UNITARY MATRIX FUNCTIONS,
WAVELET ALGORITHMS, AND STRUCTURAL PROPERTIES OF WAVELETS}

\markboth{P.E.T. Jorgensen}{Unitary Matrix Functions, Algorithms, Wavelets}

\author{Palle E. T. Jorgensen}

\address{Department of Mathematics\\
The University of Iowa\\
Iowa City, Iowa 52242, U.S.A.\\
E-mail: {jorgen@math.uiowa.edu}\\
URL: {http://www.math.uiowa.edu/\!\raisebox{-0.75ex}{\symbol{126}}\,jorgen/}\bigskip\\
\upshape
Program: ``Functional and harmonic analyses of wavelets and frames,'' 4--7~August 2004.\\
Organizers: Judith Packer, Qiyu Sun, Wai Shing Tang\\
Contribution by Palle E. T. Jorgensen to the Tutorial Sessions:}

\begin{abstract}
\label{Abs}\textsc{Abstract:} Some connections between operator theory and wavelet analysis:
Since the mid eighties, it has become clear that key tools in wavelet analysis rely crucially on operator
theory. While isolated variations of wavelets, and wavelet constructions had previously been known,
since Haar in 1910, it was the advent of multiresolutions, and sub-band filtering techniques which
provided the tools for our ability to now easily create efficient algorithms, ready for a rich variety of
applications to practical tasks. Part of the underpinning for this development in wavelet analysis is
operator theory. This will be presented in the lectures, and we will also point to a number of
developments in operator theory which in turn derive from wavelet problems, but which are of
independent interest in mathematics. Some of the material will build on chapters in a new wavelet book,
co-authored by the speaker and Ola Bratteli, see http://www.math.uiowa.edu/\raisebox{-0.75ex}{\symbol{126}}jorgen/\kern0.1em.
\end{abstract}

\vspace*{24pt}   

%\minitoc         %GENERATE OUT CONTENT TEXT AT PER CHAPTER ONLY

\section*{Contents}

Abstract\hfill\pageref{Abs}\newline
\ref{S1}.\enspace Introduction\hfill\pageref{S1}\newline
\ref{S1.2}.\enspace Index of terminology in math and in engineering\hfill\pageref{S1.2}\newline
\ref{S1.2.1}.\enspace Some background on Hilbert space\hfill\pageref{S1.2.1}\newline
\ref{S1.2.2n}.\enspace Connections to group theory\hfill\pageref{S1.2.2n}\newline
\begin{minipage}[b]{0.8\textwidth}\hangindent1in
\rule{0pt}{11pt}\ref{S1.2.2}.\enspace Some background on matrix functions in mathematics and in engineering\end{minipage}\hfill\pageref{S1.2.2}\newline
\ref{S1.3}.\enspace Motivation\hfill\pageref{S1.3}\newline
\ref{S1.3.1}.\enspace Some points of history\hfill\pageref{S1.3.1}\newline
\ref{S1.3.2}.\enspace Some early applications\hfill\pageref{S1.3.2}\newline
\ref{S2}.\enspace Signal processing\hfill\pageref{S2}\newline
\ref{S2.1}.\enspace Filters in communications engineering\hfill\pageref{S2.1}\newline
\ref{S2.2}.\enspace Algorithms for signals and for wavelets\hfill\pageref{S2.2}\newline
\ref{S2.2.1}.\enspace Pyramid algorithms\hfill\pageref{S2.2.1}\newline
\ref{S2.2.2}.\enspace Subdivision algorithms\hfill\pageref{S2.2.2}\newline
\ref{S2.2.3}.\enspace Wavelet packet algorithms\hfill\pageref{S2.2.3}\newline
\ref{S2.2.4}.\enspace Lifting algorithms: Sweldens and more\hfill\pageref{S2.2.4}\newline
\ref{S2.2bis}.\enspace Factorization theorems for matrix functions\hfill\pageref{S2.2bis}\newline
\begin{minipage}[b]{0.8\textwidth}\hangindent1in
\rule{0pt}{11pt}\ref{S2.2bis.1}.\enspace The case of polynomial functions [the polyphase matrix, joint work with 
O Bratteli]\end{minipage}\hfill\pageref{S2.2bis.1}\newline
\ref{S2.2bis.2}.\enspace General results in mathematics on matrix functions\hfill\pageref{S2.2bis.2}\newline
\ref{S2.2bis.3}.\enspace Connection between matrix functions and wavelets\hfill\pageref{S2.2bis.3}\newline
\ref{S2.2bis.3}.1.\enspace Multiresolution wavelets\hfill\pageref{S2.2bis.3.1}\newline
\begin{minipage}[b]{0.8\textwidth}\hangindent1in
\rule{0pt}{11pt}\ref{S2.2bis.3}.2.\enspace Generalized multiresolutions [joint work with L.~Baggett, K.~Merrill,
and J.~Packer]\end{minipage}\hfill\pageref{S2.2bis.3.2}\newline
\ref{SNew5Aug2.3.4}.\enspace Matrix completion\hfill\pageref{SNew5Aug2.3.4}\newline
\ref{S2.2bis.4}.\enspace  Connections between matrix functions and signal processing\hfill\pageref{S2.2bis.4}\newline
Appendix A: Topics for further research\hfill\pageref{App}\newline
\hspace*{1.5em}\ref{S3}.\enspace Connection between the discrete signals and the wavelets
%\hfill\pageref{S3}
\newline
\hspace*{1.5em}\ref{S3.1}.\enspace Wavelet geometry in $L^2(R^n)$
%\hfill\pageref{S3.1}
\newline
\begin{minipage}[b]{0.8\textwidth}\hangindent1in
\rule{0pt}{11pt}\hspace*{1.5em}\ref{S3.2}.\enspace Intertwining operators between sequence spaces $l^2$ and $L^2(R^n)$\end{minipage}
%\hfill\pageref{S3.2}
\newline
\hspace*{1.5em}\ref{S3.3}.\enspace Infinite products of matrix functions
%\hfill\pageref{S3.3}
\newline
\hspace*{1.5em}\ref{S3.3.1}.\enspace Implications for $L^2(R^n)$
%\hfill\pageref{S3.3.1}
\newline
\hspace*{1.5em}\ref{S3.3.2}.\enspace Wavelets in other Hilbert spaces of fractal measures
%\hfill\pageref{S3.3.2}
\newline
\begin{minipage}[b]{0.8\textwidth}\hangindent1in
\rule{0pt}{11pt}\hspace*{1.5em}\ref{S3.4}.\enspace Dependence of the wavelet functions on the matrix functions
which define the wavelet filters\end{minipage}
%\hfill\pageref{S3.4}
\newline
\hspace*{1.5em}\ref{S3.4.1}.\enspace Cycles
%\hfill\pageref{S3.4.1}
\newline
\hspace*{1.5em}\ref{S3.4.2}.\enspace The Ruelle-Lawton wavelet transfer operator
%\hfill\pageref{S3.4.2}
\newline
\hspace*{1.5em}\ref{S4}.\enspace Other topics in wavelets theory
%\hfill\pageref{S4}
\newline
\hspace*{1.5em}\ref{S4.1}.\enspace Invariants
%\hfill\pageref{S4.1}
\newline
\hspace*{1.5em}\ref{S4.1.1}.\enspace Invariants for wavelets: Global theory
%\hfill\pageref{S4.1.1}
\newline
\hspace*{1.5em}\ref{S4.1.2}.\enspace Invariants for wavelet filters: Local theory
%\hfill\pageref{S4.1.2}
\newline
\hspace*{1.5em}\ref{S4.2}.\enspace Function classes
%\hfill\pageref{S4.2}
\newline
\hspace*{1.5em}\ref{S4.2.1}.\enspace Function classes for wavelets
%\hfill\pageref{S4.2.1}
\newline
\hspace*{1.5em}\ref{S4.2.2}.\enspace Function classes for filters
%\hfill\pageref{S4.2.2}
\newline
\hspace*{1.5em}\ref{S4.3}.\enspace Wavelet sets
%\hfill\pageref{S4.3}
\newline
\hspace*{1.5em}\ref{S4.4}.\enspace Spectral pairs
%\hfill\pageref{S4.4}
\newline
Appendix B: Duality principles in analysis\hfill\pageref{Dual}\newline
Acknowledgements\hfill\pageref{Ack}\newline
References\hfill\pageref{Ref}\newline

\vspace{12pt}
\begin{center}
\textit{One cannot expect any serious understanding of what wavelet analysis}%
\linebreak\textit
{means without a deep knowledge of the corresponding operator theory.\bigskip
}\linebreak\settowidth{\qedskip}{\textit
{One cannot expect any serious understanding of what wavelet analysis}%
}\makebox[\qedskip]{\hfill\textsc{---Yves Meyer}$^*$}
\end{center}%

\section{\label{S1}Introduction}

\footnotetext{2000 \textit{Mathematics Subject Classification:}
Primary 42C40, 46L60, 47L30, 42A16, 43A65; Secondary 46L45, 42A65, 41A15        .}
%\footnotetext{2000 \textit{Mathematics Subject Classification:}
%41A15, 42A16, 42A65, 43A65, 46L60, 47D25.}
%\footnotetext{2000 \textit{Mathematics Subject Classification:}
%42C40, 42A16, 43A65, 42A65.}
\footnotetext{\textit{Key words and phrases:} \hyphenpenalty=10000
signal processing, matrix functions, infinite products, pyramid algorithm, 
subdivision algorithm, multiresolution, generalized multiresolution, wavelet 
packets, library of bases, wavelet filters, high-pass, low-pass filters, filter 
bank, Gabor frames, fractal measures, wavelet sets, transfer operator, Ruelle 
operator, Perron--Frobenius, dimension function, homotopy, winding number, index 
theorem, spectral representation, translation invariance,
Hilbert space, biorthogonal wavelet, Cuntz algebra,
completely positive map, Fock space, creation operators.}
\footnotetext{Work supported in part by the U.S.
National Science Foundation under grants DMS-9987777 and DMS-0139473(FRG);
financial support from the National University of Singapore.}
\footnotetext{$^*$\cite{Mey00};
see also the web page http://www.math.uiowa.edu/\raisebox{-0.75ex}{\symbol{126}}jorgen/quotes.html\kern0.1em.}
While this series of four lectures will be on the subject of wavelets, the
emphasis will be on some interconnections between topics in the mathematics of
wavelets and other areas, both within mathematics and outside. Connections to
operator theory, to quantum theory, and especially to signal processing will
be studied. Concepts such as high-pass and low-pass filters have become
synonymous with wavelet tools, but they have also had a significance from the
very start of signal processing, for example early telephone signals over
transatlantic cables. This was long before the much more recent advances in
wavelets which started in the mid-1980's (as a resumption, in fact, of ideas
going back to Alfred Haar \cite{Haa10} much earlier).

\subsection{\label{S1.2}Index of terminology in math and in engineering}

Since the mid-1980's wavelet mathematics has served to some extent as a
clearing house for ideas from diverse areas from mathematics, from
engineering, as well as from other areas of science, such as quantum theory
and optics. This makes the interdisciplinary communication difficult, as the
lingo differs from field to field; even to the degree that the same term might
have a different name to some wavelet practitioners from what is has to
others. In recognition of this fact, Chapter 1 in the recent wavelet book
\cite{BrJo02b} samples a little dictionary of relevant terms. Parts of it are
reproduced here:

\section*{Terminology}

\begin{itemize}

\item \textbf{multiresolution%
\index{multiresolution}:}
\textit{---real world:}
a set of band-pass-filter%
\index{filter}ed component
images, assembled into a mosaic of resolution%
\index{resolution} bands,
each resolution tied to a finer%
\index{resolution!finer} one and a coarser%
\index{resolution!coarser} one.

\textit{---mathematics:}
used in wavelet analysis%
\index{analysis!wavelet} and fractal analysis%
\index{analysis!fractal}, multiresolution%
\index{multiresolution}s
are systems of closed subspaces in a Hilbert%
\index{Hilbert!space decomposition@--- space decomposition} space,
such as $L^{2}\left( \mathbb{R}\right) $, with the subspaces nested,
each subspace representing a resolution%
\index{resolution}, and the
relative complement subspaces representing the detail%
\index{detail}
which is added in getting to the next finer%
\index{resolution!finer}
resolution%
\index{resolution!subspace@--- subspace} subspace.

\item \textbf{matrix%
\index{matrix!function@--- function} function:}
a function from the circle, or the one-torus,
taking values in a group of $N$-by-$N$ complex
matrices%
\index{matrix!invertible}.

\item \textbf{wavelet:}
a function $\psi$, or a finite system of functions
$\left\{ \psi_{i}\right\} $,
such that for some scale%
\index{scaling!number@--- number} number $N$ and a lattice of translation
points on $\mathbb{R}$, say $\mathbb{Z}$, a basis for
$L^{2}\left( \mathbb{R}\right) $
can be built consisting of the functions
$N^{\frac{j}{2}}\psi_{i}\left( N^{j}x-k\right) $, $j,k\in\mathbb{Z}$.

\begin{quote}
\noindent
\textit{Then dulcet music swelled\\
Concordant with the life-strings of the soul;\\
It throbbed in sweet and languid beatings there,\\
Catching new life from transitory death;\\
Like the vague sighings of a wind at even\\
That wakes the wavelets of the slumbering sea}%
$\dots$\\
\hspace*{144pt}---Shelley, \textit{Queen Mab}
\end{quote}

\item \textbf{subband filter%
\index{filter!subband}:}
\textit{---engineering:}
signals are viewed as functions
of time and frequency, the frequency function resulting
from a transform of the time function; the frequency
variable is broken up into bands, and up-sampling
and down-sampling are combined with a filtering of the
frequencies in making the connection from one band
to the next.

\textit{---wavelets:}
scaling is used in passing from one
resolution%
\index{resolution} $V$ to the next; if a scale%
\index{scaling!number@--- number} $N$ is used
from $V$ to the next finer%
\index{resolution!finer} resolution, then scaling by
$\frac{1}{N}$ takes $V$ to a coarser%
\index{resolution!coarser} resolution $V_{1}$ represented by
a subspace of $V$, but there is a set of functions
which serve as multipliers when relating $V$ to $V_{1}$,
and they are called subband filter%
\index{filter!subband}s.

\item \textbf{cascade%
\index{cascade}s:}
\textit{---real world:}
a system of successive refinement%
\index{operator!refinement}s which
pass from a scale to a finer%
\index{scaling!finer} one, and so on;
used for example in graphics algorithm%
\index{algorithm!graphics}s: starting
with control points, a refinement \mbox{matrix} and
masking%
\index{masking coefficients} coefficients are used in a cascade
algorithm%
\index{algorithm!cascade} yielding a cascade%
\index{cascade} of masking points
and a cascade approximation%
\index{approximation!cascade} to a picture.

\textit{---wavelets:}
in one dimension the scaling%
\index{scaling!number@--- number} is by a number and
a fixed simple function, for example of the form
$
\setlength{\unitlength}{14bp}\begin{picture}%
(2,1.25)(-0.5,0)
%\put(-0.5,0){\vector(1,0){3}}
%\thicklines 
\multiput(-0.25,0)(1.25,0){2}{\line(1,0){0.25}}
\put(0,1){\line(1,0){1}}
\multiput(0,0)(1,0){2}{\line(0,1){1}}
%\put(1,-0.0625){\makebox(0,0)[t]{$B$}}
%\put(0.5,1.0625){\makebox(0,0)[b]{$B'$}}
%\put(1.5,1.0625){\makebox(0,0)[b]{$B''$}}
%\put(0.5,0.5){\makebox(0,0){$\varphi$}}
\put(0,-0.0625){\makebox(0,0)[t]{$\scriptstyle 0$}}
\put(1,-0.0625){\makebox(0,0)[t]{$\scriptstyle 1$}}
%\put(2.0625,0.0625){\makebox(0,0)[bl]{$\scriptstyle 2$}}
%\put(2.4375,-0.0625){\makebox(0,0)[t]{$x$}} 
\end{picture}%
$
is chosen as the initial step for the cascade%
\index{cascade}s;
when the masking%
\index{masking coefficients} coefficients are chosen the cascade approximation%
\index{approximation!cascade} leads to a scaling%
\index{scaling!function@--- function} function.

\item \textbf{scaling%
\index{scaling!function@--- function} function:}
a function, or a distribution%
\index{distribution}, $\varphi$, defined on
the real line $\mathbb{R}$ which has the property that, for
some integer $N>1$, the coarser%
\index{resolution!coarser} version $\varphi\left( \frac{x}{N}\right) $
is in the closure (relative to some metric)
of the linear span
of the
set of translated functions
$\dots,
\varphi\left( x+1\right) $,
$\varphi\left( x\right) $,
$\varphi\left( x-1\right) $,
$\varphi\left( x-2\right) ,
\dots$.

\item \textbf{logic gate%
\index{gate!logic}s:}
\textit{---in computation}
the classical logic gate%
\index{gate!logic}s are realized as computers, for example as
electronic switching circuits with two-level voltages, say high and low.
Several gates
have two input voltages and one output,
each one allowing switching between high and low: The
output of the AND gate is high if and only if both inputs
are high. The XOR%
\index{gate!XOR} gate has high output if and only
if one of the inputs, but not more than one, is high.

\item \textbf{qubit%
\index{qubit}s:}
\textit{---in physics and in computation:}
qubit%
\index{qubit}s are the quantum analogue of
the classical bits $0$ and $1$ which are the letters of classical
computers, the qubit%
\index{qubit}s are formed of two-level quantum
systems, electrons
in a magnetic field or polarized
photons, and they are
represented in Dirac's formalism
$\left| 0\right\rangle $ and $\left| 1\right\rangle $; quantum%
\index{quantum!parallelism@--- parallelism}
theory allows superpositions, so state%
\index{state!mixed}s $\left| \psi\right\rangle 
=a\left| 0\right\rangle +b\left| 0\right\rangle $,
$a,b\in\mathbb{C}$, $\left| a\right| ^{2}+\left| b\right| ^{2}=1$,
are also admitted, and computation%
\index{quantum!computing@--- computing}
in the quantum%
\index{quantum!state@--- state} realm allows a continuum of states, as
opposed to just the two classical bits.

\textit{---mathematics:}
a chosen and distinguished basis for the
two-dimensional Hilbert%
\index{Hilbert!space@--- space} space $\mathbb{C}^{2}$ consisting of 
orthogonal%
\index{orthogonal!vectors@--- vectors}
unit vectors, denoted
$\left| 0\right\rangle $, $\left| 1\right\rangle $.

\item \textbf{universality%
\index{universality|textbf}:}
\textit{---classical computing:}
the property of
a set of logic gate%
\index{gate!logic}s that they suffice for the implementation of every program;
or of a single gate that, taken together with the NOT gate,
it suffices for the implementation of every program.

\textit{---quantum%
\index{quantum!computing@--- computing} computing:}
the property of
a set $S$ of basic quantum%
\index{gate!quantum} gates that every (invertible) gate can
be written as a sequence of steps using only gates
from $S$. Usually $S$ may be chosen to consist
of one-qubit%
\index{gate!qubit} gates and a distinguished tensor
gate $t$. An example of a choice for $t$ is
CNOT%
\index{gate!CNOT}. An alternative universal one is the
Toffoli%
\index{gate!Toffoli} gate.

\textit{---mathematics:}
the property of
a set $S$ of basic unitary matrices%
\index{matrix!unitary}
that for every $n$ and every
$u\in\mathrm{U}_{2^{n}}\left( \mathbb{C}\right) $,
there is a
factorization $u=s_{1}s_{2}\cdots s_{k}$, $s_{i}\in S$, with the understanding
that the factors $s_{i}$ are inserted in a chosen
tensor configuration of the quantum%
\index{quantum!register@--- register} register
$
\underset{n\text{ times}}
{\underbrace{\mathbb{C}^{2}\otimes\dots\otimes\mathbb{C}^{2}}}
$. Note that the factors $s_{i}$,
the number $k$, and the configuration of the $s_{i}$'s
all depend on $n$ and the gate
$u\in\mathrm{U}_{2^{n}}\left( \mathbb{C}\right) $
to be studied. The quantum wavelet
algorithm%
\index{algorithm!quantum wavelet} (\ref{eqWaveletMatrix}) is an example of such a
matrix $u$.

\item \textbf{chaos%
\index{chaos}:}
a small variation or disturbance in the
initial states or input of some system giving rise to a
disproportionate, or exponentially growing,
deviation in the resulting output trajectory,
or output data. The term is used more
generally, denoting rather drastic forms of instability; and
it is measured by the use of statistical devices, or
averaging methods.

\item $\mathrm{GL}_{N}\left( \mathbb{C}\right) $\textbf{:}
the \emph{general linear group} of
all complex $N\times N$ invertible
matrices%
\index{matrix!invertible}.

\item $\mathrm{U}_{N}\left( \mathbb{C}\right) $\textbf{:}
$=\left\{\, A\in\mathrm{GL}_{N}\left( \mathbb{C}\right) 
\mid AA^{\ast}=1_{\mathbb{C}^{N}}\,\right\} $
where $A^{\ast}$ denotes the adjoint
matrix, i.e.,
$\left( A^{\ast}\right) _{i,j}=\bar{A}_{j,i}$.

\item \textbf{transfer operator%
\index{operator!transfer} (transition operator%
\index{operator!transition}):}
\textit{---in probability:}
An operator which
transforms signals $s$ from input $s_{\operatorname{in}}$
to output $s_{\operatorname{out}}$. The
signals are represented as functions on some set $E$. In the
simplest case, the operator is linear and given in
terms of conditional probabilities $p\left( x,y\right) $.
The number $p\left( x,y\right) $
may represent the probability of a transition from $y$ to $x$ where
$x$ and $y$ are points in the set $E$. Then
\[
s_{\operatorname{out}}\left( x\right) 
=\sum_{y\in E}p\left( x,y\right) s_{\operatorname{in}}\left( y\right) .
\]

\textit{---in computation:}
Let $X$ and $Y$ be functions on a set
$E$, both taking values in $\left\{ 0,1\right\} $. Let $Y$ be the
initial state of the bit, and $X$ the final state of the
bit. If the process is governed by a probability%
\index{probability!distribution@--- distribution}
distribution $P$, then the transition probabilities $p\left( x,y\right) 
:=P\left( \left\{\, X=x\mid Y=y\,\right\} \right) $
are conditional probabilities: i.e., $p\left( x,y\right) $ is
the probability of a final bit value $x$ given an initial
value $y$, and we have
\[
P\left( \left\{ X=x\right\} \right) 
=\sum_{y\in E}p\left( x,y\right) P\left( \left\{ Y=y\right\} \right) .
\]

\textit{---in wavelet theory:}
Let $N\in\mathbb{Z}_{+}$, and let $W$ be
a positive function on
$\mathbb{T}=\left\{\, z\in\mathbb{C}\mid\left| z\right| =1\,\right\} $,
for example
$W=\left| m_{0}\right| ^{2}$ where $m_{0}$ is some low-pass wavelet
filter%
\index{filter!wavelet}%
\index{filter!low-pass} with $N$ bands.
(Positivity is only in the sense $W\geq 0$, nonnegative,
and the function $W$ may vanish on a subset of $\mathbb{T}$.)
Then define a function $p$ on
$\mathbb{T}\times\mathbb{T}$ as follows:
\[
p\left( z,w\right) 
=\left\{
\begin{array}{ll}
\left( \frac{1}{N}\right) W\left( w\right) 
 &\text{\quad if }w^{N}=z, \\
\,0
 &\text{\quad for all other values of }w.
\end{array}
\right. 
\]
We arrive at the transfer operator%
\index{operator!transfer} $R_{W}$, i.e.,
the operator transforming functions on $\mathbb{T}$ as follows:
\[
s_{\operatorname{out}}\left( z\right) 
=\left( R_{W}s_{\operatorname{in}}\right) \left( z\right) 
=\frac{1}{N}\sum_{w^{N}=z}W\left( w\right) 
s_{\operatorname{in}}\left( w\right) .
\]

\item \textbf{coherence:}
\textit{---in mathematics and physics:}
The vectors $\psi_{i}$
that make up a tight frame%
\index{frame!tight},
one which is not an orthonormal basis, are
said to be subjected to \emph{coherence}.
So
coherent%
\index{coherent!vector@--- vector} vector systems in Hilbert%
\index{Hilbert!space@--- space} space
are viewed as bases which generalize the
more standard concept of orthonormal%
\index{orthonormal!basis@--- basis}
bases from harmonic analysis%
\index{analysis!harmonic}.
A striking
feature of the wavelets with compact%
\index{compact support}
support, which are based on scaling, is that
the variet%
\index{variety!of wavelet@--- of wavelets or wavelet filters}ies of the two kinds of bases
can be well understood geometrically. For
example,
the collapse of the wavelet orthogonal%
\index{orthogonal!wavelet@--- wavelet}ity
relations, degenerating into coherent%
\index{coherent!vector@--- vector}
vectors, happens on a subvariet%
\index{variety!of wavelet@--- of wavelets or wavelet filters}y of a
lower dimension.

%\enlargethispage{2\baselineskip}
More generally, coherent%
\index{coherent!vector@--- vector} vectors in mathematical physics
often arise with a continuous index, even
if the Hilbert%
\index{Hilbert!space@--- space} space is separable, i.e., has a countable
orthonormal%
\index{orthonormal!basis@--- basis} basis.
This is illustrated by a vector
system $\left\{ \psi_{r,s}\right\} $, which
should be thought of as a
continuous analogue, i.e., a
version where a sum gets
replaced with an integral
\[
C_{\psi}^{-1}\iint\limits_{\mathbb{R}^{2}}\frac{dr\,ds}{r^{2}}
\left| \left\langle\,\psi_{r,s}\mid f\,\right\rangle \right| ^{2}
=\left\Vert f\right\Vert ^{2}.
\]
For more details, see also Section 3.3 of \cite{Dau92}%
\index{Daubechies0@Ingrid Daubechies}
and Chapter 3 of \cite{Kai94}.

In quantum%
\index{quantum!mechanics@--- mechanics} mechanics,
one talks, for example,
about coherent state%
\index{state!coherent}s in connection
with wavefunctions of the harmonic oscillator. Combinations of stationary
wavefunctions from different energy eigenvalues vary periodically
in time, and the question is
which of the continuously varying wavefunctions one may
use to expand an unknown function in without encountering overcompleteness
of the basis. The methods of ``coherent state%
\index{state!coherent}s'' are methods for using
these kinds of
%\linebreak
functions (which fit some problems elegantly) while avoiding
the difficulties of overcompleteness.  The term ``coherent'' applies
when you succeed in avoiding those difficulties by some means or
other.
Of course, for students who have just learned about the classic
complete orthonormal%
\index{orthonormal!basis@--- basis}
basis of stationary eigenfunction%
\index{eigenfunction}s, ``coherent state%
\index{state!coherent}'' methods
at first may seem like
a daring relaxation of the rules of orthogonality, so that the term seems to
stand for total freedom!

\end{itemize}

\subsubsection{\label{S1.2.1}Some background on Hilbert space}

Wavelet theory is the art of finding a special kind of basis in Hilbert space.
Let $\mathcal{H}$ be a Hilbert space over $\mathbb{C}$ and denote the inner
product $\left\langle \,\cdot\mid\cdot\,\right\rangle $. For us, it
is assumed linear in the second variable. If $\mathcal{H}=L^{2}\left(
\mathbb{R}\right)  $, then%
\begin{equation}
\left\langle \,f\mid g\,\right\rangle :=\int_{\mathbb{R}}\overline{f\left(
x\right)  }\,g\left(  x\right)  \,dx.\label{eq1.2.1}%
\end{equation}
If $\mathcal{H}=\ell^{2}\left(  \mathbb{Z}\right)  $, then%
\begin{equation}
\left\langle \,\xi\mid\eta\,\right\rangle :=\sum_{n\in\mathbb{Z}}\bar{\xi}_{n}%
\eta_{n}.\label{eq1.2.2}%
\end{equation}
Let $\mathbb{T}=\mathbb{R}/2\pi\mathbb{Z}$. If $\mathcal{H}=L^{2}\left(
\mathbb{T}\right)  $, then%
\begin{equation}
\left\langle \,f\mid g\,\right\rangle :=\frac{1}{2\pi}\int_{-\pi}^{\pi
}\overline{f\left(  \theta\right)  }\,g\left(  \theta\right)  \,d\theta
.\label{eq1.2.3}%
\end{equation}
Functions $f\in L^{2}\left(  \mathbb{T}\right)  $ have Fourier series: Setting
$e_{n}\left(  \theta\right)  =e^{in\theta}$,
\begin{equation}
\hat{f}\left(  n\right)  :=\left\langle \,e_{n}\mid f\,\right\rangle =\frac
{1}{2\pi}\int_{-\pi}^{\pi}e^{-in\theta}f\left(  \theta\right)  \,d\theta
,\label{eq1.2.4}%
\end{equation}
and%
\begin{equation}
\left\Vert f\right\Vert _{L^{2}\left(  \mathbb{T}\right)  }^{2}=\sum
_{n\in\mathbb{Z}}\left\vert \hat{f}\left(  n\right)  \right\vert
^{2}.\label{eq1.2.5}%
\end{equation}
Similarly if $f\in L^{2}\left(  \mathbb{R}\right)  $, then%
\begin{equation}
\hat{f}\left(  t\right)  :=\int_{\mathbb{R}}e^{-ixt}f\left(  x\right)
\,dx,\label{eq1.2.6}%
\end{equation}
and%
\begin{equation}
\left\Vert f\right\Vert _{L^{2}\left(  \mathbb{R}\right)  }^{2}=\frac{1}{2\pi
}\int_{\mathbb{R}}\left\vert \hat{f}\left(  t\right)  \right\vert
^{2}\,dt.\label{eq1.2.7}%
\end{equation}

Let $J$ be an index set. We shall only need to consider the case when $J$ is
countable. Let $\left\{  \psi_{\alpha}\right\}  _{\alpha\in J}$ be a family of
nonzero vectors in a Hilbert space $\mathcal{H}$. We say it is an
\emph{orthonormal basis} (ONB) if%
\begin{equation}
\left\langle \,\psi_{\alpha}\mid\psi_{\beta}\,\right\rangle =\delta
_{\alpha,\beta}\text{\qquad(Kronecker delta)}\label{eq1.2.8}%
\end{equation}
and if
\begin{equation}
\sum_{\alpha\in J}\left\vert \left\langle \,\psi_{\alpha}\mid f\,\right\rangle
\right\vert ^{2}=\left\Vert f\right\Vert ^{2}\text{\qquad holds for all }%
f\in\mathcal{H}.\label{eq1.2.9}%
\end{equation}
If only (\ref{eq1.2.9}) is assumed, but not (\ref{eq1.2.8}), we say that
$\left\{  \psi_{\alpha}\right\}  _{\alpha\in J}$ is a (normalized) \emph{tight
frame}. We say that it is a \emph{frame} with \emph{frame constants} $0<A\leq
B<\infty$ if%
\[
A\left\Vert f\right\Vert ^{2}\leq\sum_{\alpha\in J}\left\vert \left\langle
\,\psi_{\alpha}\mid f\,\right\rangle \right\vert ^{2}\leq B\left\Vert
f\right\Vert ^{2}\text{\qquad holds for all }f\in\mathcal{H}.
\]
Introducing the rank-one operators $Q_{\alpha}:=\left\vert \psi_{\alpha
}\right\rangle \left\langle \psi_{\alpha}\right\vert $ of Dirac's terminology,
see \cite{BrJo02b}, we see that $\left\{  \psi_{\alpha}\right\}  _{\alpha\in
J}$ is an ONB if and only if the $Q_{\alpha}$'s are projections and%
\begin{equation}
\sum_{\alpha\in J}Q_{\alpha}=I\qquad(=\text{the identity operator in
}\mathcal{H}).\label{eq1.2.10}%
\end{equation}
It is a (normalized) tight frame if and only if (\ref{eq1.2.10}) holds but
with no further restriction on the rank-one operators $Q_{\alpha}$. It is a
frame with frame constants $A$ and $B$ if the operator
\begin{equation}
S:=\sum_{\alpha\in J}Q_{\alpha}\label{eq1.2.11}%
\end{equation}
satisfies%
\[
AI\leq S\leq BI
\]
in the order of hermitian operators. (We say that operators $H_{i}=H_{i}%
^{\ast}$, $i=1,2$, satisfy $H_{1}\leq H_{2}$ if $\left\langle
\,f\mid H_{1}f\,\right\rangle \leq\left\langle \,f\mid H_{2}f\,\right\rangle $
holds for all $f\in\mathcal{H}$).

Wavelets in $L^{2}\left(  \mathbb{R}\right)  $ are generated by simple
operations on one or more functions $\psi$ in $L^{2}\left(  \mathbb{R}\right)
$, the operations come in pairs, say scaling and translation, or
phase-modulation and translations. If $N\in\left\{  2,3,\dots\right\}  $ we
set%
\begin{equation}
\psi_{j,k}\left(  x\right)  :=N^{j/2}\psi\left(  N^{j}x-k\right)  \text{\qquad
for }j,k\in\mathbb{Z}.\label{eq1.2.12}%
\end{equation}

\subsubsection{\label{S1.2.2n}Connections to group theory}

We stress the
discrete%
\index{wavelet!discrete} wavelet%
\index{wavelet!transform@--- transform} transform. But
the first line in the two tables
below
is the continuous one. 
It is the only treatment
we give to the continuous%
\index{wavelet!continuous}
wavelet%
\index{wavelet!transform@--- transform} transform, and the
corresponding \emph{coherent%
\index{coherent!vector@--- vector} vector
decomposition%
\index{decomposition}s}. 
But, as is
stressed in \cite{Dau92}%
\index{Daubechies0@Ingrid Daubechies}, \cite{Kai94}, and \cite{KaLe95},
the continuous version came first.

\addvspace{\bigskipamount}
\noindent
$\renewcommand{\arraystretch}{1.25}%\addtolength{\tabcolsep}{-1.1111pt}%
\begin{tabular}
[c]
{p{0.15\textwidth}
||p{0.375\textwidth}
|p{0.375\textwidth}}
\multicolumn{3}{c}{\parbox{0.9\textwidth}
{Summary of and variations on the resolution
of the identity operator $1$ in $L^{2}$ or in $\ell^{2}$, for $\psi$ and
$\tilde{\psi}$ where $\psi_{r,s}\left( x\right) 
=r^{-\frac{1}{2}}\psi\left( \frac{x-s}{r}\right) $,
$C_{\psi}=\int_{\mathbb{R}}\frac{d\omega}{\left| \omega\right| }
\left| \smash{\hat{\psi}}\left( \omega\right) \right| ^{2}<\infty$,
similarly for $\tilde{\psi}$ and 
$C_{\psi,\tilde{\psi}}=\int_{\mathbb{R}}\frac{d\omega}{\left| \omega\right| }
\overline{\hat{\psi}\left( \omega\right) }\,
\Hat{\Tilde{\psi}}\left( \omega\right) $:\rule[-12pt]{0pt}{12pt}}}\\
$N=2$ & Overcomplete Basis & Dual Bases \\\hline\hline
continu\-ous resolution%
\index{resolution} 
& $\displaystyle C_{\psi}^{-1}\!\!\iint\limits_{\mathbb{R}^{2}}^{{}}
\frac{dr\,ds}{r^{2}}
\left| \psi_{r,s}\right\rangle \!\left\langle \psi_{r,s}\right| 
$\par
\hfill
$\displaystyle {}
=1_{L^{2}}\vphantom{\iint}$ 
& $\displaystyle 
C_{\!\psi,\tilde{\psi}}^{-1}\!\iint\limits_{\mathbb{R}^{2}}^{{}}
\frac{dr\,ds}{r^{2}}
\left| \psi_{r,s}\right\rangle 
\!\left\langle \smash{\tilde{\psi}_{r,s}}\right| \vphantom{\tilde{\psi}_{r,s}}
$\par
\hfill
$\displaystyle {}
=1_{L^{2}}\vphantom{\iint}$ \\\hline
discrete resolution%
\index{resolution} 
& {\raggedright $\displaystyle 
\sum_{\vphantom{k}j\in\mathbb{Z}}\sum_{\vphantom{j}k\in\mathbb{Z}}
\left| \psi_{j,k}\right\rangle \left\langle \psi_{j,k}\right| 
=1_{L^{2}}$\,, \\
$\psi_{j,k}$ corresponding to \\
$r=2^{-j}$, 
$s=k2^{-j}$}
& $\displaystyle 
\sum_{\vphantom{k}j\in\mathbb{Z}}\sum_{\vphantom{j}k\in\mathbb{Z}}
\left| \psi_{j,k}\right\rangle 
\left\langle \smash{\tilde{\psi}_{j,k}}\right| \vphantom{\tilde{\psi}_{j,k}}
=1_{L^{2}}$ 
\\\hline\hline
$N\geq 2$ 
& Isometries%
\index{operator!isometric} in $\ell^{2}$ 
& Dual Operator System in $\ell^{2}$ \\\hline\hline
sequence spaces 
& {\raggedright $\displaystyle \sum_{i=0_{\mathstrut}}^{N-1^{\mathstrut}}
S_{i}S_{i}^{\ast}=1_{\ell^{2}}$\,, \\
where
$\displaystyle S_{0},\dots ,S_{N-1}$ \\
are adjoint%
\index{operator!adjoint}s 
to the \\
quadrature mirror filter%
\index{filter!quadrature mirror} \\
operators $F_{i}$, i.e., 
$\displaystyle S_{i}=F_{i}^{\ast}$} 
& {\raggedright $\displaystyle \sum_{i=0_{\mathstrut}}^{N-1^{\mathstrut}}
S_{i}\tilde{S}_{i}^{\ast}=1_{\ell^{2}}$\,, \\
for a dual \\
operator system \\
$\displaystyle S_{0},\dots ,S_{N-1}$, \\
$\displaystyle \smash{\tilde{S}_{0},\dots ,\tilde{S}_{N-1}}$} \\\hline
\end{tabular}
$

\addvspace{\bigskipamount}
\noindent
$\renewcommand{\arraystretch}{1.375}%\addtolength{\tabcolsep}{-1.1111pt}%
\begin{tabular}
[c]
{p{0.40\textwidth}
|p{0.52\textwidth}}
\multicolumn{2}{c}{\parbox{0.9\textwidth}
{Consult Chapter 3 of \cite{Kai94}
for the continuous resolution%
\index{resolution}, and
Section 2.2 of \cite{BrJo02b} for the discrete resolution%
\index{resolution}.
If $h,k$ are vectors in a Hilbert%
\index{Hilbert!space@--- space} space
$\mathcal{H}$, then the operator
$A=\left| h\right\rangle \left\langle k\right| $ is defined by the identity
$\left\langle \,u\mid Av\,\right\rangle =\left\langle \,u\mid h\,\right\rangle \left\langle \,k\mid v\,\right\rangle $
for all $u,v\in\mathcal{H}$. Then the assertions in the first table amount to:\rule[-12pt]{0pt}{12pt}}}\\
$\displaystyle C_{\psi}^{-1}\iint\limits_{\mathbb{R}^{2}}%^{{}}
\frac{dr\,ds}{r^{2}}
\left| \left\langle \,\psi_{r,s}\mid f\,\right\rangle \right| ^{2}
$\par
\hfill
$\displaystyle {}
=\left\Vert f\right\Vert _{L^{2}}^{2}\quad
\forall\,f\in L^{2}\left( \mathbb{R}\right) \vphantom{\iint}$
& $\displaystyle
C_{\!\psi,\tilde{\psi}}^{-1}\!\iint\limits_{\mathbb{R}^{2}}%^{{}}
\frac{dr\,ds}{r^{2}}
\left\langle \,f\mid \psi_{r,s}\,\right\rangle 
\left\langle \,\smash{\tilde{\psi}_{r,s}}\mid g\,\right\rangle \vphantom{\tilde{\psi}_{r,s}}
$\par
\hfill
$\displaystyle {}
=\left\langle \,f\mid g\,\right\rangle \quad
\forall\,f,g\in L^{2}\left( \mathbb{R}\right) \vphantom{\iint}$ \\\hline
$\displaystyle 
\sum_{\vphantom{k}j\in\mathbb{Z}}\sum_{\vphantom{j}k\in\mathbb{Z}_{\mathstrut}}
\left| \left\langle \,\psi_{j,k}\mid f\,\right\rangle \right| ^{2}
$\par
\hfill
$\displaystyle {}
=\left\Vert f\right\Vert_{L^{2}}^{2}\quad
\forall\,f\in L^{2}\left( \mathbb{R}\right) \vphantom{\iint}$
& $\displaystyle 
\sum_{\vphantom{k}j\in\mathbb{Z}}\sum_{\vphantom{j}k\in\mathbb{Z}_{\mathstrut}}
\left\langle \,f\mid \psi_{j,k}\,\right\rangle 
\left\langle \,\smash{\tilde{\psi}_{j,k}}\mid g\,\right\rangle \vphantom{\tilde{\psi}_{j,k}}
$\par
\hfill
$\displaystyle {}
=\left\langle \,f\mid g\,\right\rangle \quad
\forall\,f,g\in L^{2}\left( \mathbb{R}\right) \vphantom{\iint}$ \\\hline
$\displaystyle \sum_{i=0_{\mathstrut}}^{N-1^{\mathstrut}}
\left\Vert S_{i}^{\ast}c\right\Vert ^{2}=\left\Vert c\right\Vert ^{2}
\quad\forall\,c\in\ell^{2}$ 
& $\displaystyle \sum_{i=0_{\mathstrut}}^{N-1^{\mathstrut}}
\left\langle \,S_{i}^{\ast}c\mid \smash{\tilde{S}_{i}^{\ast}}d\,\right\rangle \vphantom{\tilde{S}_{i}^{\ast}}
=\left\langle \,c\mid d\,\right\rangle \quad\forall\,c,d\in\ell^{2}$ \\\hline
\end{tabular}
$

\addvspace{\bigskipamount}
A
function $\psi$ satisfying the resolution identity is called
a \emph{coherent%
\index{coherent!vector@--- vector|textbf} vector} in mathematical physics. The
representation theory for the $\left( ax+b\right) $-group,
i.e., the matrix group
$G=\left\{\, \left( 
\begin{smallmatrix}
a & b \\
0 & 1
\end{smallmatrix}
\right) \mid a\in\mathbb{R}_{+},\; b\in\mathbb{R}\,\right\} $,
serves as its underpinning. Then the tables above
illustrate how the
$\left\{ \psi_{j,k}\right\} $ wavelet
system arises from a discretization of
the following unitary representation%
\index{representation!unitary} of $G$:
\begin{equation}
\left( U_{\left( 
\begin{smallmatrix}
a & b \\
0 & 1
\end{smallmatrix}
\right) }f\right) \left( x\right) 
=a^{-\frac{1}{2}}f\left( \frac{x-b}{a}\right) 
\label{eqUnitaryRep}
\end{equation}
acting on $L^{2}\left( \mathbb{R}\right) $.
This unitary representation%
\index{representation!unitary} also explains
the discretization step in passing
from the first line to the second in
the tables above. The
functions $\left\{\, \psi_{j,k}\mid j,k\in\mathbb{Z}\,\right\} $ which
make up a wavelet system
result from the choice of a
suitable coherent%
\index{coherent!vector@--- vector} vector $\psi\in L^{2}\left( \mathbb{R}\right) $,
and then setting
\begin{equation}
\psi_{j,k}\left( x\right) 
=\left( U_{\left( 
\begin{smallmatrix}
2^{-j} & k\cdot 2^{-j} \\
0 & 1
\end{smallmatrix}
\right) }\psi\right) \left( x\right) 
=2^{\frac{j}{2}}\psi\left( 2^{j}x-k\right) .
\label{eqCoherentVec}
\end{equation}
Even though
this representation lies at the historical origin of the
subject of wavelets
(see \cite{DGM86}%
\index{Meyer0@Yves Meyer}%
\index{Daubechies0@Ingrid Daubechies}),
the $\left( ax+b\right) $-group seems
to be now largely forgotten
in the next generation of the wavelet community.
But Chapters 1--3 of \cite{Dau92}%
\index{Daubechies0@Ingrid Daubechies} still
serve as a beautiful presentation
of this (now much ignored) side
of the subject. It also serves
as a link to mathematical physics
and to classical analysis%
\index{analysis!classical}.

Since the representation%
\index{representation!unitary} $U$ in (\ref{eqUnitaryRep})
on $L^{2}\left( \mathbb{R}\right) $,
when a unitary $U$ is defined from (\ref{eqUnitaryRep}) setting $a=2$, $b=0$,
$\left( Uf\right) \left( x\right) :=2^{-\frac{1}{2}}f\left( \frac{x}{2}\right) $,
leaves invariant the Hardy%
\index{Hardy space} space
\begin{equation}
\mathcal{H}_{+}
=\left\{\, f\in L^{2}\left( \mathbb{R}\right) 
\mid \operatorname{supp}\left( \smash{\hat{f}}\right) 
\subset\left\lbrack 0,\infty\right\rangle \,\right\} ,
\label{eqHardySpace}
\end{equation}
formula (\ref{eqCoherentVec})%
\index{coherent!vector@--- vector}
suggests that it would be
simpler to look for wavelets
in $\mathcal{H}_{+}$. After all, it is a smaller
space, and it is natural to try
to use the causality features of
$\mathcal{H}_{+}$ implied by the support
condition in (\ref{eqHardySpace}). 
Moreover, in the world of the Fourier%
\index{Fourier1@Fourier!transform@--- transform}
transform, the two operations of the
formulas (\ref{eqUnitaryRep}) and (\ref{eqCoherentVec})%
\index{coherent!vector@--- vector} take
the simpler forms
\begin{equation}
\hat{f}\longmapsto a^{\frac{1}{2}}e^{-ibt}\hat{f}\left( at\right) 
\text{\quad and\quad}
\hat{\psi}
\longmapsto 2^{\frac{j}{2}}e^{-i2^{j}kt}\hat{\psi}\left( 2^{j}t\right) .
\label{eqHardyCoherent}
\end{equation}
So in the early
nineties, this was an open problem
in the theory,
i.e., whether or not there are wavelets in the Hardy%
\index{Hardy space} space;
but it received
a beautiful answer in \cite{Aus95}. Auscher
showed that there are no
wavelet%
\index{wavelet!function@--- function} functions $\psi$ in $\mathcal{H}_{+}$ which satisfy
the following mild regularity properties:
\begin{align*}
\left( R_{0}\right) \quad
 & \hat{\psi}\text{ is continuous;}
\\
\left( R_{\varepsilon}\right) \quad
 & \text{for some }\varepsilon\in\mathbb{R}_{+},\;
\hat{\psi}\left( t\right)
=\mathcal{O}\left( \left| t\right| ^{\varepsilon}\right) \\
 & \qquad\qquad\qquad\qquad\qquad
\text{ and }\smash{\hat{\psi}\left( t\right)
=\mathcal{O}\left( \left( 1+\left| t\right| \right) 
^{-\varepsilon-\frac{1}{2}}\right) ,\; t\in\mathbb{R}.}
\end{align*}

Comparison
of formulas (\ref{eqUnitaryRep}) and (\ref{eqCoherentVec})%
\index{coherent!vector@--- vector}
shows that
The traditional discrete%
\index{wavelet!discrete}
wavelet%
\index{wavelet!transform@--- transform} transform may be viewed as
the restriction to a subgroup $H$ of a classical unitary
representation%
\index{representation!unitary} of $G$.
The
unitary representation%
\index{representation!unitary}s of $G$ are
completely understood:
the set of irreducible unitary
representation%
\index{representation!irreducible}%
\index{representation!unitary}s consists of
two infinite-dimensional inequivalent subrepresentation%
\index{representation!sub@sub---}s of
the representation (\ref{eqUnitaryRep}) on $L^{2}\left( \mathbb{R}\right) $,
together with the one-dimensional representations
$\left( 
\begin{smallmatrix}
a & b \\
0 & 1
\end{smallmatrix}
\right) \rightarrow a^{ik}$ parameterized by $k\in\mathbb{R}$.
(The two subrepresentation%
\index{representation!sub@sub---}s of (\ref{eqUnitaryRep}) are obtained by
restricting to $f\in L^{2}\left( \mathbb{R}\right) $ with
$\operatorname{supp}\hat{f}\subseteq\left\langle -\infty,0\right\rbrack $ and
$\operatorname{supp}\hat{f}\subseteq\left\lbrack 0,\infty\right\rangle $,
respectively.)
However,
the subgroup $H$
of $G$ has a rich variety
of inequivalent
infinite-dimensional
representations that
do not arise as restrictions of
(\ref{eqUnitaryRep}), or of any representation of $G$.
The group $H$
considered
in (\ref{eqCoherentVec})%
\index{coherent!vector@--- vector} is a
semidirect product%
\index{product!semidirect}
(as is $G$):
it is of the form
\begin{multline}
H_{N}=\left\{ 
\begin{pmatrix}
a & b \\
0 & 1
\end{pmatrix}
\biggm| a=N^{j},\;
b=\sum_{i\in\mathbb{Z}}n_{i}N^{i},\;
j\in\mathbb{Z},\;
n_{i}\in\mathbb{Z},
\vphantom{\text{ where the }\sum_{i}\text{ summation is finite}}\right. \\
\smash[t]{\left. \vphantom{
\begin{pmatrix}
a & b \\
0 & 1
\end{pmatrix}
\biggm| a=N^{j},\;
b=\sum_{i\in\mathbb{Z}}n_{i}N^{i},\;
j\in\mathbb{Z},\;
n_{i}\in\mathbb{Z},}
\text{ where the }\sum_{i}\text{ summation is finite}\,\right\} .}
\label{eqSemidirectProd}
\end{multline}
(In the jargon of pure algebra, the nonabelian group $H_{N}$ is the semidirect
product%
\index{product!semidirect} of the two abelian groups $\mathbb{Z}$ and
$\mathbb{Z}\left\lbrack \frac{1}{N}\right\rbrack $,
with a naturally defined action of $\mathbb{Z}$ on
$\mathbb{Z}\left\lbrack \frac{1}{N}\right\rbrack $.)

The papers \cite{DaLa98}, \cite{Jor01a}, \cite{BaMe99}, \cite{HLPS99},
\cite{LPT01}%
\index{Packer0@Judy Packer}, and \cite{BreJo91} show that it is
possible to use these nonclassical
representations of $H$ for the construction
of unexpected classes of wavelets,
the wavelet%
\index{wavelet!set@--- set} sets being the
most notable ones. Recall
that a subset $E\subset\mathbb{R}$ of finite measure
is a \emph{wavelet%
\index{wavelet!set@--- set} set} if
$\hat{\psi}=\chi_{E}^{{}}$ is such that, for some
$N\in\mathbb{Z}_{+}$, $N\geq 2$, the functions
$\left\{\, N^{\frac{j}{2}}\psi\left( N^{j}x-k\right) 
\bigm| j,k\in\mathbb{Z}\,\right\} $ form
an orthonormal%
\index{orthonormal!basis@--- basis} basis for $L^{2}\left( \mathbb{R}\right) $.
Until the work of Larson and others,
see \cite{DaLa98} and \cite{HLPS99}, it was
not even clear that
wavelet%
\index{wavelet!set@--- set} sets
$E$ could exist in the case $N>2$.
The paper \cite{LPT01}%
\index{Packer0@Judy Packer} develops and extends the
representation theory for the subgroups $H_{N}$
independently of the ambient group
$G$ and shows that each $H_{N}$
has continuous series of representations
which account for the wavelet%
\index{wavelet!set@--- set}
sets. The role of the representations
of the groups $H_{N}$ and their generalizations
for the study of wavelets was first
stressed in \cite{BreJo91}.

There is a different transform which is analogous
to the wavelet%
\index{wavelet!transform@--- transform} transform of (\ref{eqUnitaryRep})--(\ref{eqCoherentVec})%
\index{coherent!vector@--- vector}, but
yet different in a number of respects. It
is the Gabor transform%
\index{transform!Gabor}, and
it has a history of its own.
Both are special cases of the
following construction: Let $G$ be a nonabelian matrix group
with center $C$, and let $U$ be
a unitary irreducible representation%
\index{representation!irreducible}
of $G$ on the Hilbert%
\index{Hilbert!space@--- space} space $L^{2}\left( \mathbb{R}\right) $.
When $\psi\in L^{2}\left( \mathbb{R}\right) $ is given, we
may define a transform
\begin{equation}
\left( T_{\psi}f\right) \left( \xi\right) 
:=\left\langle \,U\left( \xi\right) \psi\mid f\,\right\rangle ,
\text{\quad for }f\in L^{2}\left( \mathbb{R}\right) 
\text{ and }\xi\in G\diagup C.
\label{eqGabtrans}
\end{equation}
It
turns out that there are classes of
matrix groups, such as the $ax+b$ group,
or the $3$-dimensional group of upper triangular
matrices, which have transforms $T_{\psi}$ admitting effective
discretizations. This means that
it is possible to find a
vector $\psi\in L^{2}\left( \mathbb{R}\right) $, and a discrete
subgroup $\Lambda\subset G\diagup C$, such that
the restriction to $\Lambda$ of the
transform $T_{\psi}$ in (\ref{eqGabtrans}) is
injective from $L^{2}\left( \mathbb{R}\right) $ into
functions on $\Lambda$.

There are many books on transform
theory, and here we are only making the
connection to wavelet theory.
The book \cite{Per86} contains much
more detail on the group-theoretic
approach to these continuous and discrete
coherent%
\index{coherent!vector@--- vector} vector transforms.
%\addtolength{\textheight}{0.125\baselineskip}

\subsubsection{\label{S1.2.2}Some background on matrix functions in mathematics and in engineering}

\label{Hom}
One of our coordinates for the landscape of multiresolution%
\index{multiresolution!wavelet@--- wavelet}
wavelets takes the form of a geometric index%
\index{index!geometric}. In fact, it
involves a traditional operator-theoretic index%
\index{index} with values in $\mathbb{Z}$.
When it is
identified with a winding%
\index{winding number} number or a counting of homotopy%
\index{wavelet!homotopy of@homotopy of ---s} classes,
it serves also as
a Fredholm index%
\index{index!Fredholm} of an associated Toeplitz operator%
\index{operator!Toeplitz}. An
orthogonal%
\index{orthogonal!wavelet basis@--- wavelet basis} dyadic wavelet%
\index{wavelet!dyadic}%
\index{wavelet!basis@--- basis} basis has its wavelet%
\index{wavelet!function@--- function} function $\psi$ satisfying
the normalization
$\left\Vert \psi\right\Vert _{L^{2}\left( \mathbb{R}\right) }=1$,
i.e., $\psi$ is a vector of norm one in the Hilbert%
\index{Hilbert!space@--- space} space
$L^{2}\left( \mathbb{R}\right) $. In the lingo of
quantum%
\index{quantum!state@--- state} theory, $\psi$ is therefore a pure state%
\index{state!pure quantum@pure (quantum)}, and the
$x$-coordinate is an observable called the position. The integral
$E_{\psi}\left( x\right) 
=\int_{\mathbb{R}}x\left| \psi\left( x\right) \right| ^{2}\,dx$
is the expected value of the position.
If $\psi_{H}$ denotes the standard Haar%
\index{Haar1@Haar!wavelet@--- wavelet} function in (\ref{eqHaarWavelet}),
then
clearly $E_{\psi_{H}}\left( x\right) =\frac{1}{2}$.
Also note the translation formula
$E_{\psi\left( \,\cdot \,-k\right) }\left( x\right) 
=E_{\psi}\left( x\right) +k$.
We showed in Corollary 2.4.11 of \cite{BrJo02b}, completely generally, that
the other orthonormal%
\index{orthonormal!wavelet@--- wavelet} wavelets $\psi$ have expected values in
the set $\frac{1}{2}+\mathbb{Z}$. Hence, after $\psi$ is translated by
an integer, you cannot distinguish it from the Haar%
\index{Haar1@Haar!wavelet@--- wavelet} wavelet $\psi_{H}$
in (\ref{eqHaarWavelet})
by looking only at the expected value of its
position coordinate. The translation integer $k$ turns
out to be a winding%
\index{winding number} number. Our result
holds more generally when the definition of
$E_{\psi}\left( x\right) $ is adapted to a wider wavelet context, as we
showed in Chapter 6 of \cite{BrJo02b}; but in all cases, there
is a winding%
\index{winding number} number which produces the above-mentioned integer
translate $k$.

The issue of connectedness for various classes of wavelets
is a general question which has been addressed previously in
the wavelet literature; see, e.g., \cite{HLPS99}, \cite{HeWe96},
\cite{StZh01}%
\index{Strang0@Gilbert Strang}, and \cite{ReWe98}. Here we bring homotopy%
\index{wavelet!homotopy of@homotopy of ---s} to bear on the
question, and we identify the connected%
\index{connected components} components when the compact support%
\index{compact support} is fixed and given. We show among other things that
for a fixed $K%
\index{K 1@$K\sb 1$!class@--- -class}_{1}$-class a homotopy%
\index{wavelet!homotopy of@homotopy of ---s} may take place within a variet%
\index{variety!of wavelet@--- of wavelets or wavelet filters}y of
wavelets which is specified by a slightly bigger support than the
initially given one.

An important point of our present discussion, beyond the mere fact of compact
support%
\index{compact support}, is the size of the support of the wavelets in question. Consider two
wavelets $A$ and $B$ of a certain support size. Then our
first results in 
%Section \ref{Hom} below
this section
also specify the paths $C\left(  t\right)  $, if any, which connect $A$ and
$B$, and in particular the size of the support of the wavelets corresponding
to $C\left(  t\right)  $.
In \cite{BrJo02b}, we 
%will then, in Section \ref{Glo},
treat connectivity in
the wider context of noncompactly supported wavelets,
following at the outset
\cite{Gar99}%
\index{Garrigos0@G. Garrig\'{o}s}, which considers scale%
\index{scaling!number@--- number} number $N=2$, and wavelets $\psi$
satisfying%
\begin{equation}
\left\{  2^{\frac{j}{2}}\psi\left(  2^{j}x-k\right)  \right\}  _{j,k\in
\mathbb{Z}}\text{ is an orthonormal%
\index{orthonormal!basis@--- basis} basis (ONB) for }L^{2}\left(
\mathbb{R}\right)  . \label{eqChHomMay15.1}%
\end{equation}
Garrig\'{o}s%
\index{Garrigos0@G. Garrig\'{o}s} considers, for $\frac{1}{2}<\alpha\leq\infty$, the class
$\mathcal{W}_{\alpha}$ of wavelets $\psi$ such that%
\begin{equation}
\int_{\mathbb{R}}\left\vert \psi\left(  x\right)  \right\vert ^{2}\left(
1+\left\vert x\right\vert ^{2}\right)  ^{\alpha}\,dx<\infty,
\label{eqChHomMay15.2}%
\end{equation}
and there is an $\varepsilon=\varepsilon\left(  \psi\right)  $ such that%
\begin{equation}
\int_{\mathbb{R}}\left\vert \hat{\psi}\left(  t\right)  \right\vert
^{2}\left(  1+\left\vert t\right\vert ^{2}\right)  ^{\varepsilon}\,dt<\infty,
\label{eqChHomMay15.3}%
\end{equation}
i.e., the wavelet is supposed to have some degree of smoothness in the sense
of Sobolev.

We now turn to the group of functions $U\colon\mathbb{T}\rightarrow
\mathrm{U}\left(  N\right)  $, where $\mathrm{U}\left(  N\right)  $ denotes
the group of all complex $N$-by-$N$ matrices. The functions will not be
assumed continuous in general. The continuous functions will be designated
$C\left(  \mathbb{T},\mathrm{U}\left(  N\right)  \right)  $. Each function in
$C\left(  \mathbb{T},\mathrm{U}\left(  N\right)  \right)  $ has a $K_{1}%
$-class, also called a winding number; see \cite{BrJo02b}. The functions in
$C\left(  \mathbb{T},\mathrm{U}\left(  N\right)  \right)  $ with finite
Fourier expansion will be called \emph{Fourier polynomials}, also if they are
functions which take values in $\mathrm{U}\left(  N\right)  $.

\begin{prop}
\label{ProHom.1}Let $U\in C\left(  \mathbb{T},\mathrm{U}\left(  N\right)
\right)  $ be a Fourier polynomial%
\index{Fourier1@Fourier!polynomial@--- polynomial}, and assume that $K%
\index{K 1@$K\sb 1$!class@--- -class}_{1}\left(  U\right)
=d\in\mathbb{Z}$. Then $U$ is homotopic in $C\left(  \mathbb{T},\mathrm{U}%
\left(  N\right)  \right)  $ to%
\begin{equation}
V\left(  z\right)  =z^{d}p\oplus\left(  1_{N}-p\right)  \label{eqHom.1}%
\end{equation}
where $p$ is the one-dimensional projection%
\index{operator!projection} onto the first coordinate slot in $\mathbb{C}^{N}$, and if
$U$ has the form%
\begin{equation}
U\left(  z\right)  =\sum_{k=-D}^{D}z^{k}a_{k} ,\label{eqHom.2}%
\end{equation}
then $U$ may be homotopically deform%
\index{homotopy deformation}ed to $V$ in $C\left(  \mathbb{T}%
,\mathrm{U}\left(  N\right)  \right)  $ through Fourier polynomials%
\index{Fourier1@Fourier!polynomial@--- polynomial} of degree%
\index{degree!of a Fourier polynomial@--- of a Fourier polynomial} at most $\left|  d\right|  +ND$.

This proposition remains true if the word ``Fourier polynomial%
\index{Fourier1@Fourier!polynomial@--- polynomial}'' is replaced by ``polynomial'' and $a_{k}=0$ for
$k=-D,-D+1,\dots,-1$. In that case $d\in\mathbb{Z}_{+}$ and $U$ may be
homotopically deform%
\index{homotopy deformation}ed to $V$ in the loop semigroup of polynomial unitaries in
$C\left(  \mathbb{T},\mathrm{U}\left(  N\right)  \right)  $ through
polynomials of degree%
\index{degree!of a Fourier polynomial@--- of a Fourier polynomial} at most $d$.
\end{prop}

\begin{proof}
Multiplying $U$ by $z^{D}$, we obtain a polynomial $z^{D}U\left(  z\right)  $
of degree%
\index{degree!of a Fourier polynomial@--- of a Fourier polynomial} $2D$ mapping $\mathbb{T}$ into $\mathrm{U}\left(  N\right)  $. Then
$K%
\index{K 1@$K\sb 1$!class@--- -class}_{1}\left(  z^{D}U\right)  =d+ND$. By Proposition 3.3 of \cite{BrJo02a}, there
exist $d+ND$ one-dimensional projections%
\index{operator!projection} $p_{1},p_{2},\dots,p_{d+ND}$ in $M_{N}\left(
\mathbb{C}\right)  $ and a unitary%
\index{operator!unitary} $V_{0}\in M_{N}\left(  \mathbb{C}\right)  $
such that%
\begin{equation}
z^{D}U\left(  z\right)  =V_{0}\prod_{k=1}^{d+ND}\left(  1-p_{i}+zp_{i}\right)  .
\label{eqHom.3}%
\end{equation}
(See 
%Exercise 
\S\ \ref{Exe6DauSwe-11} for a related, but different, decomposition%
\index{decomposition}.)
Now, deform%
\index{deformation!continuous}ing each of the $p_{i}$'s continuously through one-dimensional
projection%
\index{operator!projection}s to the projection%
\index{operator!projection} $p_{0}$ onto the first coordinate direction, and deform%
\index{deformation!continuous}ing
$V_{0}$ in $\mathrm{U}\left(  N\right)  $ into $1_{N}$, we see that $z^{D}U\left(
z\right)  $ can be deform%
\index{deformation!continuous}ed into
\begin{equation}
\prod_{k=1}^{d+ND}\left(  1-p_{0}+zp_{0}\right)  =1-p_{0}+z^{d+ND}p_{0}.
\label{eqHom.4}%
\end{equation}
Thus $U\left(  z\right)  $ itself is deform%
\index{deformation!continuous}ed into%
\begin{equation}
z^{-D}\left(  1-p_{0}\right)  +z^{d+\left(  N-1\right)  D}p_{0}.
\label{eqHom.5}%
\end{equation}
But writing $\left(  1-p_{0}\right)  $ as a sum of $N-1$ one-dimensional
projection%
\index{operator!projection}s $q_{1},\dots,q_{N-1}$, we have that the unitary%
\index{operator!unitary} that
$U\left(  z\right)  $ is deform%
\index{deformation!continuous}ed into is%
\begin{equation}
\prod_{k=1}^{N-1}\left(  \left(  1-q_{k}\right)  +z^{-D}q_{k}\right)
\cdot\left(  1+z^{d+\left(  N-1\right)  D}p_{0}\right)  , \label{eqHom.6}%
\end{equation}
and next deform%
\index{deformation!continuous}ing each of the $q_{k}$ in this decomposition%
\index{decomposition} into $p_{0}$, we
see that $U\left(  z\right)  $ is deform%
\index{deformation!continuous}ed into%
\begin{equation}
\prod_{k=1}^{N-1}\left(  \left(  1-p_{0}\right)  +z^{-D}p_{0}\right)
\cdot\left(  1+z^{d+\left(  N-1\right)  D}p_{0}\right)  =\left(
1-p_{0}\right)  +z^{d}p_{0}. \label{eqHom.7}%
\end{equation}
The crude estimate $\left|  d\right|  +ND$ on the degree%
\index{degree!of a Fourier polynomial@--- of a Fourier polynomial} of the Fourier
polynomials%
\index{Fourier1@Fourier!polynomial@--- polynomial} occurring during the deformation%
\index{homotopy deformation} is straightforward.

To prove the last statement in the proposition one does not need to multiply
$U$ by $z^{D}$, and the proof simplifies. Note in particular that $D\leq d$
(assuming $a_{D}\neq0$).
\end{proof}

\begin{rem}
\label{RemHom.2}We do not know if Proposition \textup{\ref{ProHom.1}} is true
if $C\left(  \mathbb{T},\mathrm{U}\left(  N\right)  \right)  $ is replaced by
$C\left(  \mathbb{T},\mathrm{GL}\left(  N\right)  \right)  $. It is known from
Lemma 11.2.12 of \cite{RLL00} that if $A\in C\left(  \mathbb{T},\mathrm{GL}%
\left(  N\right)  \right)  $ is a polynomial of degree%
\index{degree!of a Fourier polynomial@--- of a Fourier polynomial} $1$ in $z$,
then $A$ can
be homotopically deform%
\index{homotopy deformation}ed through first-order polynomials in $C\left(
\mathbb{T},\mathrm{GL}\left(  N\right)  \right)  $ to a unitary%
\index{operator!unitary} of the form
$z\rightarrow zp+\left(  1_{N}-p\right)  $ for some projection%
\index{operator!projection} $p$, and hence Proposition \textup{\ref{ProHom.1}} for
$C\left(  \mathbb{T},\mathrm{GL}\left(  N\right)  \right)  $ would follow if
any polynomial $A\in C\left(  \mathbb{T},\mathrm{GL}\left(  N\right)  \right)
$ could be factored into first-order polynomials. It is also clear, since any
element $A\in C\left(  \mathbb{T},\mathrm{GL}\left(  N\right)  \right)  $ can
be homotopically deform%
\index{homotopy deformation}ed into $z^{d}p\oplus\left(  1_{N}-p\right)  $ in
$C\left(  \mathbb{T},\mathrm{GL}\left(  N\right)  \right)  $, that if $A$ is a
Fourier polynomial%
\index{Fourier1@Fourier!polynomial@--- polynomial}, then $A$ can be homotopically deform%
\index{homotopy deformation}ed into
$z^{d}p\oplus\left(  1_{N}-p\right)  $ through Fourier polynomials%
\index{Fourier1@Fourier!polynomial@--- polynomial}. This follows by compactness and the
Stone--Weierstra\ss\ theorem%
\index{theorem!Stone Weierstrass@Stone--Weierstra\ss } \textup{(Lemma 11.2.3 of \cite{RLL00})}.
For our purposes in
wavelet theory, though, we would need a computable upper bound for the
degree%
\index{degree!of a Fourier polynomial@--- of a Fourier polynomial} of the Fourier polynomials%
\index{Fourier1@Fourier!polynomial@--- polynomial}.
\end{rem}

%Because of the difficulty mentioned in Remark \ref{RemHom.2} we will only
%formulate the main theorem%
%\index{theorem!homotopy} on homotopy%
%\index{wavelet!homotopy of@homotopy of ---s} for the orthogonal case, although bits
%and pieces of it obviously are also true in the biorthogonal%
%\index{biorthogonal} case with suitable modifications. 
For ease of reference
we will now list the correspondences between the various objects that interest
us in this case. 
%(The formalism below for the more general biorthogonal%
%\index{biorthogonal} case can be read off (\ref{eqPFeFeb8.20}%
%)--(\ref{eqRTO.43}) in Section \ref{Cyc} and Theorem \ref{CorCas.1}.) 
These
objects are:

\begin{enumerate}
\renewcommand{\theenumi}{\roman{enumi}}

\item \label{HomCorr(1)}matrix%
\index{matrix!function@--- function} functions, $A\colon\mathbb{T}\rightarrow\mathrm{U}_{N}\left(
\mathbb{C}\right)  $, satisfying the normalization%
\begin{equation}
A\left(  1\right)  =H,\qquad H_{k,l}=\frac{1}{\sqrt{N}}e^{i2\pi kl/N}%
,\;k,l=0,\dots,N-1, \label{eqHom.12}%
\end{equation}

\item \label{HomCorr(2)}high- and low-pass wavelet filter%
\index{filter!wavelet}s%
\index{filter!high-pass}%
\index{filter!low-pass} $m_{i}$, $i=0,1,\dots,N-1$, satisfying
\begin{equation}
\sum_{w^{N}=z}\overline{m_{i}\left(  w\right)  }\,m_{j}\left(  w\right)
=N\delta_{ij},\qquad i,j=0,\dots,N-1, \label{eqHom.13}%
\end{equation}
and%
\begin{equation}
m_{0}\left(  1\right)  =\sqrt{N}, \label{eqHom.14}%
\end{equation}

\item \label{HomCorr(3)}scaling%
\index{scaling!function@--- function} functions $\varphi$ together with wavelet%
\index{wavelet!generator@--- generator} generators $\psi_{i}$.
\end{enumerate}

\noindent 
We did not specify the continuity%
\index{continuity} and regularity requirements of the functions
$A$, $m_i$, $\varphi$, $\psi_i$ above. This will be done
differently in different contexts and the classes clearly depend on these added requirements. We will now restrict to the case
that the functions $\varphi$
and $\psi_{i}$ have compact support%
\index{compact support} in $\left[  0,\infty\right\rangle $, i.e., that $A$
and $m_{i}$ are polynomials in $z$. Thus $z\rightarrow A\left(  z\right)  $ is
a polynomial function with%
\begin{equation}
\left(  A\left(  z\right)  \right)  ^{\ast}A\left(  z\right)  =1,\qquad
z\in\mathbb{T}. \label{eqHom.11}%
\end{equation}

\noindent
\textbf{Scaling%
\index{scaling!function@--- function} functions/wavelet%
\index{wavelet!generator@--- generator} generators to wavelet filter%
\index{filter!wavelet}s} $\left(
\varphi,\psi\right)  \mapsto m$

\addvspace{\medskipamount} \noindent One defines $a_{n}$ by
\begin{equation}
\varphi\left(  x\right)  =\sqrt{N}\sum_{n\in\mathbb{Z}}a_{n}\varphi\left(
Nx-n\right)  , \label{eqInt.1b}%
\end{equation}
(cf.\ (\ref{eqInt.1})) and then $m_{0}$ by 
\begin{equation}
m_{0}\left(  z\right)  =\sum_{n}a_{n}z^{n}, \label{eqIntNew.2}%
\end{equation}
or one uses%
\begin{equation}
\sqrt{N}\hat{\varphi}\left(  Nt\right)  =m_{0}\left(  t\right)  \hat{\varphi
}\left(  t\right)  \label{eqHom.27}%
\end{equation}
directly. Then the high-pass filters%
\index{filter!high-pass} $m_{i}$, $i=1,\dots,N-1$, can be derived from
(\ref{eqHom.21})
below.
If we are in the generic case (\ref{eqMul.12}),
%of Theorem \ref{ThmMul.1},
we may also recover the Fourier coefficient%
\index{Fourier1@Fourier!coefficients@--- coefficient(s)}s
$a_{n}^{\left( i\right) }$ of $m_i$ by
\begin{align*}
a_{n}^{\left( i\right) }
& = \left( 1/\sqrt{N}\right) \left\langle\,\varphi\left( \,\cdot \, -
n\right) \mid \psi_{i}\left( \,\cdot \, /N\right) \right\rangle
\\
& = \left\langle\, \varphi\left( \,\cdot \, - n\right) \mid U\psi_{i}\,\right\rangle 
\text{\qquad (with }\psi_{0} = \varphi\text{),}
\end{align*}
where $U\psi_{i}\left( x\right) :=N^{-1/2}\psi_{i}\left( x/N\right) $.
In particular it follows in this generic case that if the
scaling and wavelet function%
\index{scaling!function@--- function}%
\index{wavelet!function@--- function}s have compact support%
\index{wavelet!compactly supported} and the filter%
\index{filter}s are
Lipschitz%
\index{Lipschitz!function@--- function}, then the filter%
\index{filter}s are Fourier polynomial%
\index{Fourier1@Fourier!polynomial@--- polynomial}s.
Is this true also in the nongeneric tight frame%
\index{frame!tight} case? 

\addvspace{\medskipamount}

Now, if $D\in\mathbb{N}$, define:%
\settowidth{\hangboxwidth}{$\displaystyle \operatorname{WF}\left(  D\right)  $}
\begin{align}
\llap{$\bullet\mkern6mu$}
\makebox[\hangboxwidth][r]{$\displaystyle\operatorname{MF}\left( D\right) $} &=
\begin{minipage}[t]{0.895\displayboxwidth}\raggedright
the set of polynomial functions in $z\in\mathbb{T}$ in $C\left(  \mathbb
{T},\mathrm{U}_{N}\left(\mathbb{C}\right)\right
)  $ of degree%
\index{degree!of a Fourier polynomial@--- of a Fourier polynomial} at most $D$ satisfying (\ref{eqHom.12});\end{minipage}%
\label{eqHom.28}\\[\smallskipamount]
\llap{$\bullet\mkern6mu$}
\makebox[\hangboxwidth][r]{$\displaystyle\operatorname{WF}\left( D\right) $} &=
\begin{minipage}[t]{0.895\displayboxwidth}\raggedright
the set of $N$-tuples of wavelet filter%
\index{filter!wavelet}s $\left(  m_{0},\dots,m_{N-1}
\right)  $ such that all $m_{i}$ are polynomials in $z\in\mathbb
{T}$ of degree%
\index{degree!of a Fourier polynomial@--- of a Fourier polynomial} at most $D$ satisfying (\ref{eqHom.13})
and (\ref{eqHom.14});\end{minipage}%
\label{eqHom.29}\\[\smallskipamount]
\llap{$\bullet\mkern6mu$}
\makebox[\hangboxwidth][r]{$\displaystyle\operatorname{SF}\left( D\right) $} &=
\begin{minipage}[t]{0.895\displayboxwidth}\raggedright
the set of $N$-tuples $\left(  \varphi,\psi_{1},\dots,\psi_{N-1}%
\right)  $ of scaling%
\index{scaling!function@--- function} functions/ wavelet%
\index{wavelet!function@--- function} functions with support in $\left
[0,D\right]$.\end{minipage}
\label{eqHom.30}%
\end{align}
The spaces $\operatorname{MF}\left(  D\right)  $, $\operatorname{WF}\left(
D\right)  $, and $\operatorname{SF}\left(  D\right)  $ may be equipped with
the obvious topologies, coming in the first two cases from,
for example, the
$L^{\infty}$-norm over $z$, and in the last case either from the $L^{2}\left(
\mathbb{R}\right)  $-norm or, as will be more relevant, the
tempered-distribution%
\index{distribution} topology. By virtue of Proposition 3.2 in \cite{BrJo02a},
$\operatorname{MF}\left(  D\right)  $ has the structure of a compact algebraic
variet%
\index{variety!algebraic}y, and so by (\ref{eqHom.17}) below,
%and (\ref{eqHom.18}), 
$\operatorname{WF}%
\left(  D\right)  $ is a compact algebraic variet%
\index{variety!algebraic}y. It is clear from (\ref{eqHom.17}) that 
the map $A\rightarrow m$ maps $\operatorname{MF}\left(
D\right)  $ into $\operatorname{WF}\left(  \left(  D+1\right)  N-1\right)  $,
and 
%from (\ref{eqHom.18}) 
that $m\rightarrow A$ maps $\operatorname{WF}\left(
\left(  D+1\right)  N-1\right)  $ into $\operatorname{MF}\left(  D\right)  $.
Furthermore, it is clear from (\ref{eqInt.1b}) and (\ref{eqHom.21}) that
$m\rightarrow\left(  \varphi,\psi\right)  $ maps $\operatorname{WF}\left(
\left(  N-1\right)  D\right)  $ into $\operatorname{SF}\left(  D\right)  $,
and conversely $\left(  \varphi,\psi\right)  \rightarrow m$ maps
$\operatorname{SF}\left(  D\right)  $ into $\operatorname{WF}\left(  \left(
N-1\right)  D\right)  $.

Now, let a subindex $0$ denote the subsets of these various spaces such that
the condition 
\begin{equation}
\operatorname*{Spec}\left(  R_{0}\right)  \cap\mathbb{T}=\left\{  1\right\}
\text{\quad and\quad}\dim\left\{  g\in\mathcal{K}_{\left\lfloor \frac{D}%
{N-1}\right\rfloor },\;R\left(  g\right)  =g\right\}  =1 \label{eqHom.23}%
\end{equation}
holds. It is known that the set of points such
that (\ref{eqHom.23}) does not hold is a lower-dimensional subvariet%
\index{variety!of wavelet@--- of wavelets or wavelet filters}y of the
various variet%
\index{variety!of wavelet@--- of wavelets or wavelet filters}ies, see
Section 6 of \cite{Jor01b}, and hence $\operatorname{MF}%
_{0}\left(  D\right)  $, $\operatorname{WF}_{0}\left(  D\right)  $, and
$\operatorname{SF}_{0}\left(  D\right)  $ contain the generic points in
$\operatorname{MF}\left(  D\right)  $, $\operatorname{WF}\left(  D\right)  $,
and $\operatorname{SF}\left(  D\right)  $.

We now summarize the local connectivity results by stating the following theorem.
The proof may be found in \cite{BrJo02b}, where this is Theorem 2.1.3.

\begin{theo}
\label{ThmHom.3}Let $k\in\mathbb{N}$. Equip the space $\operatorname{SF}%
\left(  kN+1\right)  $ of scaling%
\index{scaling!function@--- function} functions/wavelet function%
\index{wavelet!function@--- function}s with support in $\left[  0,kN+1\right]  $ with the
tempered-distribution%
\index{distribution} topology. Then $\operatorname{SF}\left(  kN+1\right)  $
is homeomorphic%
\index{homeomorphism} to a compact algebraic variet%
\index{variety!algebraic}y. Furthermore, for two
elements $\left(  \varphi_{0},\psi_{0}\right)  ,\left(  \varphi_{1},\psi
_{1}\right)  \in\operatorname{SF}\left(  kN+1\right)  $, the following
conditions are equivalent:
\renewcommand{\theenumi}{\alph{enumi}}

\begin{enumerate}
\item \label{ThmHom.3(1)}The elements $\left(  \varphi_{0},\psi_{0}\right)  $
and $\left(  \varphi_{1},\psi_{1}\right)  $ can be connected to each other by
a continuous path in $\operatorname{SF}\left(  NkN+1\right)  $;

\item \label{ThmHom.3(2)}$K%
\index{K 1@$K\sb 1$!class@--- -class}_{1}\left(  \varphi_{0},\psi_{0}\right)
=K_{1}\left(  \varphi_{1},\psi_{1}\right)  $;

\item \label{ThmHom.3(3)}The elements $\left(  \varphi_{0},\psi_{0}\right)  $
and $\left(  \varphi_{1},\psi_{1}\right)  $ can be connected to each other by
a continuous path in some $\operatorname{SF}\left(  K\right)  $.
\end{enumerate}

\noindent Thus, $\operatorname{SF}\left(  kN+1\right)  $ is divided into
$Nk\left(  N-1\right)  +1$ components which are connected over
$\operatorname{SF}\left(  NkN+1\right)  $.
\end{theo}

\subsection{\label{S1.3}Motivation}

In addition to the general background
material in the present section, the
reader may find a more detailed
treatment of some of the current research
trends in wavelet analysis in the
following papers: \cite{Jor03a} (a book review), \cite{Jor03b} (a survey), and
the research papers  \cite{DuJo03}, \cite{DuJo04a}, \cite{DuJo04b}, \cite{DuJo04c}, \cite{Jor04a}, and \cite{Jor04b}.

As a mathematical subject, the theory of
wavelets draws on tools from mathematics itself, such as harmonic analysis and
numerical analysis. But in addition there are exciting links to areas outside
mathematics. The connections to electrical and computer engineering, and to
image compression and signal processing in particular, are especially
fascinating. These interconnections of research disciplines may be illustrated
with the two subjects (1)~wavelets and (2)~subband filtering [from signal
processing]. While they are quite different, and have distinct and independent
lives, and even have different aims, and different histories, they have in
recent years found common ground. It is a truly amazing success story.
Advances in one area have helped the other: subband filters are absolutely
essential in wavelet algorithms, and in numerical recipes used in subdivision
schemes, for example, and especially in JPEG 2000---an important and
extraordinarily successful image-compression code. JPEG uses nonlinear
approximations and harmonic analysis in spaces of signals of bounded
variation. Similarly, new wavelet approximation techniques have given rise to
the kind of data-compression which is now used by the FBI [via a patent held
by two mathematicians] in digitizing fingerprints in the {U.S}. It is the
happy marriage of the two disciplines, signal processing and wavelets, that
enriches the union of the subjects, and the applications, to an extraordinary
degree. While the use of high-pass and low-pass filters has a long history in
signal processing, dating back more than fifty years, it is only relatively
recently, say the mid-1980's, that the connections to wavelets have been made.
Multiresolutions from optics are the bread and butter of wavelet algorithms,
and they in turn thrive on methods from signal processing, in the quadrature
mirror filter construction, for example. The effectiveness of multiresolutions
in data compression is related to the fact that multiresolutions are modelled on the
familiar positional number system: the digital, or dyadic, representation of
numbers. Wavelets are created from scales of closed subspaces of the Hilbert
space $L^{2}\left(  \mathbb{R}\right)  $ with a scale of subspaces
corresponding to the progression of bits in a number representation. While
oversimplified here, this is the key to the use of wavelet algorithms in
digital representation of signals and images. The digits in the classical
number representation in fact are quite analogous to the frequency subbands
that are used \emph{both} in signal processing and in wavelets.

The two functions%
\begin{equation}%
\begin{array}
[c]{ccc}%
\varphi\left(  x\right)  =\left\{  \renewcommand{\arraystretch}{1.125}%
\begin{array}
[c]{ll}%
1\phantom{-} & 0\leq x<1\\
0 & \text{elsewhere}%
\end{array}
\right.  & \text{and} & \psi\left(  x\right)  =\left\{
\renewcommand{\arraystretch}{1.125}%
\begin{array}
[c]{ll}%
1 & 0\leq x<\frac{1}{2}\\
-1 & \frac{1}{2}\leq x<1\\
0 & \text{elsewhere}%
\end{array}
\right. \\
&  & \\
\setlength{\unitlength}{36pt}\begin{picture}(2.25,3)(-0.5,-1.5)\put(-0.5,0){\vector(1,0){2.25}}\put(0,-1.5){\vector(0,1){3}}\thicklines\put(-0.125,0){\line(1,0){0.125}}\put(0,0){\line(0,1){1}}\put(0,1){\line(1,0){1}}\put(1,1){\line(0,-1){1}}\put(1,0){\line(1,0){0.125}}\put(0.5,1.0625){\makebox(0,0)[b]{$\varphi$}}\end{picture} &
&
\setlength{\unitlength}{36pt}\begin{picture}(2.25,3)(-0.5,-1.5)\put(-0.5,0){\vector(1,0){2.25}}\put(0,-1.5){\vector(0,1){3}}\thicklines\put(-0.125,0){\line(1,0){0.125}}\put(0,0){\line(0,1){1}}\put(0,1){\line(1,0){0.5}}\put(0.5,1){\line(0,-1){2}}\put(0.5,-1){\line(1,0){0.5}}\put(1,-1){\line(0,1){1}}\put(1,0){\line(1,0){0.125}}\put(0.5,1.0625){\makebox(0,0)[b]{$\psi$}}\end{picture}\\
\text{Father function\quad} &  & \text{Mother function\quad}\\
\text{(a)\quad} &  & \text{(b)\quad}%
\end{array}
\label{eq0}%
\end{equation}
capture in a glance the refinement identities%
\[
\varphi\left(  x\right)  =\varphi\left(  2x\right)  +\varphi\left(
2x-1\right)  \text{\qquad and\qquad}\psi\left(  x\right)  =\varphi\left(
2x\right)  -\varphi\left(  2x-1\right)  .
\]
The two functions are clearly orthogonal in the inner product of $L^{2}\left(
\mathbb{R}\right)  $, and the two closed subspaces $\mathcal{V}_{0}$ and
$\mathcal{W}_{0}$ generated by the respective integral translates%
\begin{equation}
\left\{  \varphi\left(  \,\cdot\,-k\right)  :k\in\mathbb{Z}\right\}
\text{\qquad and\qquad}\left\{  \psi\left(  \,\cdot\,-k\right)  :k\in
\mathbb{Z}\right\}  \label{eq1}%
\end{equation}
satisfy%
\begin{equation}
U\mathcal{V}_{0}\subset\mathcal{V}_{0}\text{\qquad and\qquad}U\mathcal{W}%
_{0}\subset\mathcal{V}_{0} \label{eq2}%
\end{equation}
where $U$ is the dyadic scaling operator $Uf\left(  x\right)  =2^{-1/2}%
f\left(  x/2\right)  $. The factor $2^{-1/2}$ is put in to make $U$ a unitary
operator in the Hilbert space $L^{2}\left(  \mathbb{R}\right)  $. This version
of Haar's system naturally invites the question of what other pairs of
functions $\varphi$ and $\psi$ with corresponding orthogonal subspaces
$\mathcal{V}_{0}$ and $\mathcal{W}_{0}$ there are such that the same
invariance conditions (\ref{eq2}) hold. The invariance conditions hold if
there are coefficients $a_{k}$ and $b_{k}$ such that the scaling identity%
\begin{equation}
\varphi\left(  x\right)  =\sum_{k\in\mathbb{Z}}a_{k}\varphi\left(
2x-k\right)  \label{eq3}%
\end{equation}
is solved by the father function, called $\varphi$, and the mother function
$\psi$ is given by%
\begin{equation}
\psi\left(  x\right)  =\sum_{k\in\mathbb{Z}}b_{k}\varphi\left(  2x-k\right)  .
\label{eq4}%
\end{equation}
A fundamental question is the converse one: Give simple conditions on two
sequences $\left(  a_{k}\right)  $ and $\left(  b_{k}\right)  $ which
guarantee the existence of $L^{2}\left(  \mathbb{R}\right)  $-solutions
$\varphi$ and $\psi$ which satisfy the orthogonality relations for the
translates (\ref{eq1}). How do we then get an orthogonal basis from this? The
identities for Haar's functions $\varphi$ and $\psi$ of (\ref{eq0})(a) and
(\ref{eq0})(b) above make it clear that the answer lies in a similar tiling
and matching game which is implicit in the more general identities (\ref{eq3})
and (\ref{eq4}). Clearly we might ask the same question for other scaling
numbers, for example $x\rightarrow3x$ or $x\rightarrow4x$ in place of
$x\rightarrow2x$. Actually a direct analogue of the visual interpretation from
(\ref{eq0}) makes it clear that there are no nonzero locally integrable
solutions to the simple variants of (\ref{eq3}),%
\begin{align}
\varphi\left(  x\right)   &  =\frac{3}{2}\left(  \varphi\left(  3x\right)
+\varphi\left(  3x-2\right)  \right) \label{eq5}\\%
%TCIMACRO{\TeXButton{or}{\intertext{or}}}%
%BeginExpansion
\intertext{or}%
%EndExpansion
\varphi\left(  x\right)   &  =2\left(  \varphi\left(  4x\right)
+\varphi\left(  4x-2\right)  \right)  . \label{eq6}%
\end{align}
There \emph{are} nontrivial solutions to (\ref{eq5}) and (\ref{eq6}), to be
sure, but they are versions of the Cantor Devil's Staircase functions, which
are prototypes of functions which are not locally integrable.

\begin{figure}[ptb]
\setlength{\unitlength}{0.675bp} \begin{picture}(478,644)
\put(0,526){\includegraphics[bb=0 1 118 117,width=79.65bp]{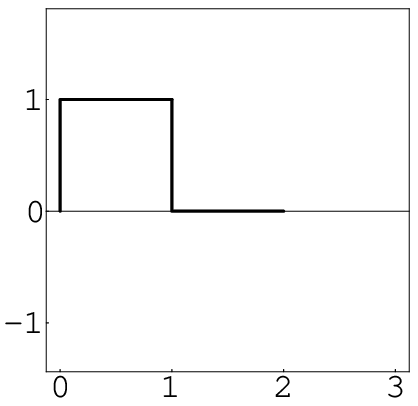}}
\put(120,526){\includegraphics[bb=0 1 118 117,width=79.65bp]{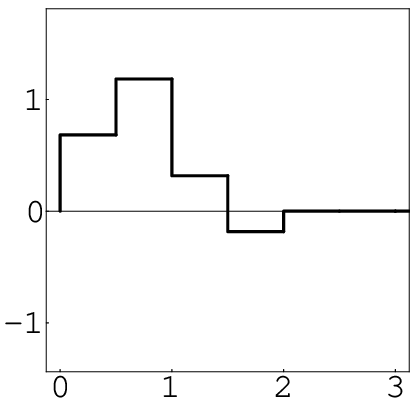}}
\put(240,526){\includegraphics[bb=0 1 118 117,width=79.65bp]{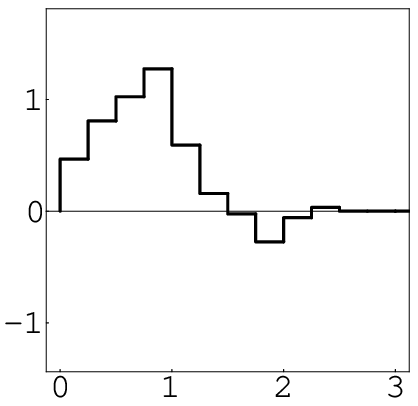}}
\put(360,526){\includegraphics[bb=0 1 118 117,width=79.65bp]{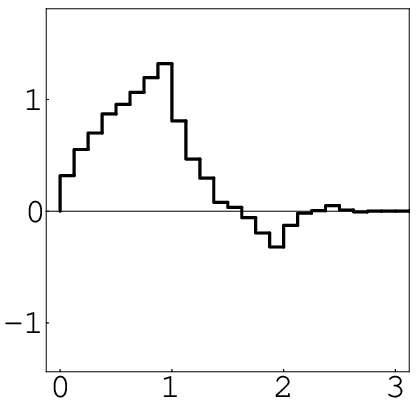}}
\put(0,406){\includegraphics[bb=0 1 118 117,width=79.65bp]{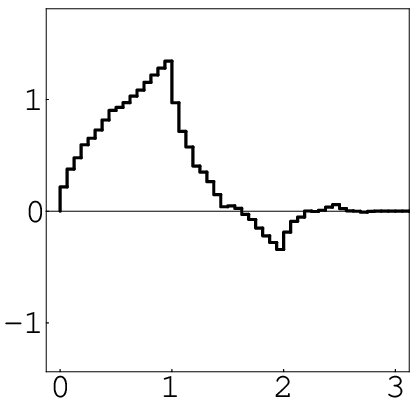}}
\put(120,406){\includegraphics[bb=0 1 118 117,width=79.65bp]{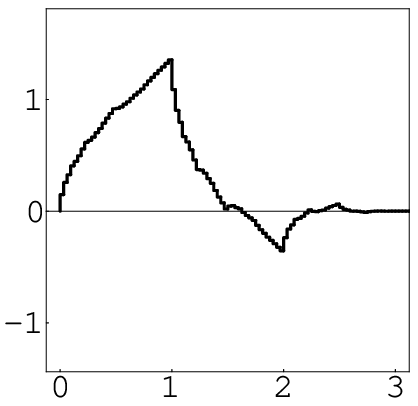}}
\put(298,344){\includegraphics[bb=0 1 180 179,width=121.5bp]{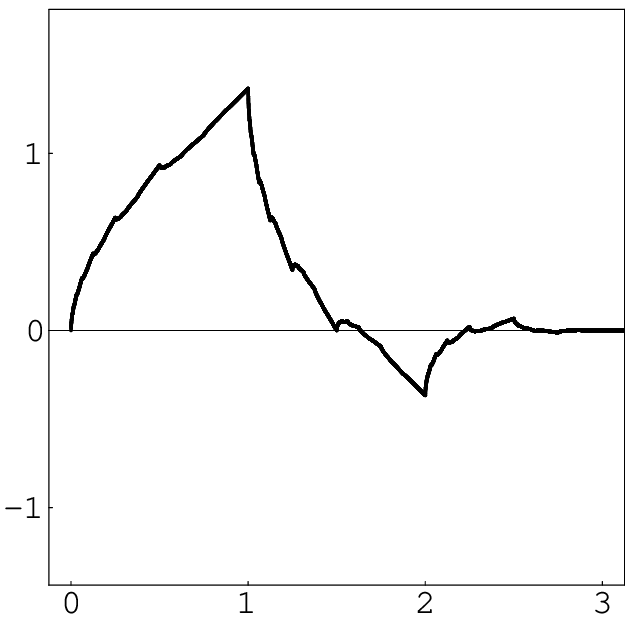}}
\put(269,458){\makebox(0,0)[b]{$\cdots$}}
\put(388,479){\makebox(0,0)[b]{$\varphi$}}
\put(392,327){\makebox(0,0)[b]{\small (a): Father function}}
\put(0,199){\includegraphics[bb=0 1 118 117,width=79.65bp]{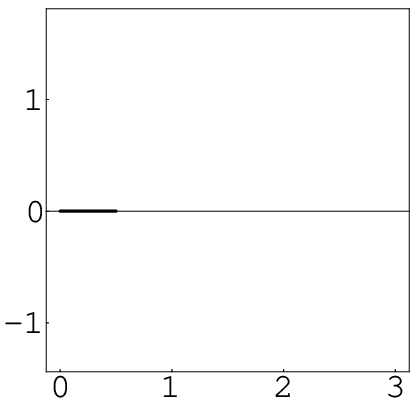}}
\put(120,199){\includegraphics[bb=0 1 118 117,width=79.65bp]{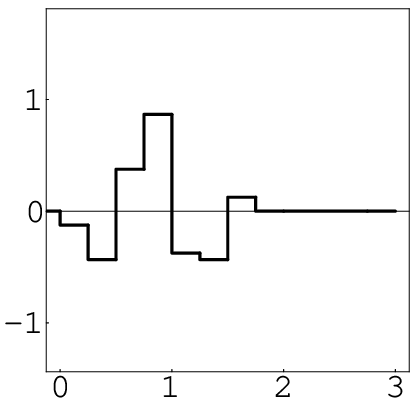}}
\put(240,199){\includegraphics[bb=0 1 118 117,width=79.65bp]{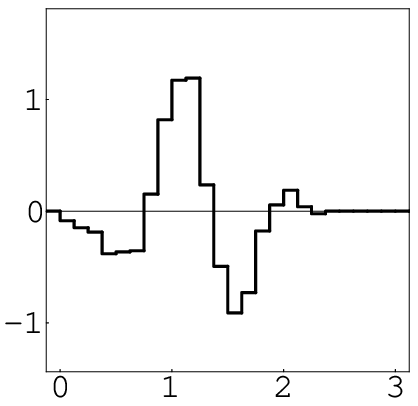}}
\put(360,199){\includegraphics[bb=0 1 118 117,width=79.65bp]{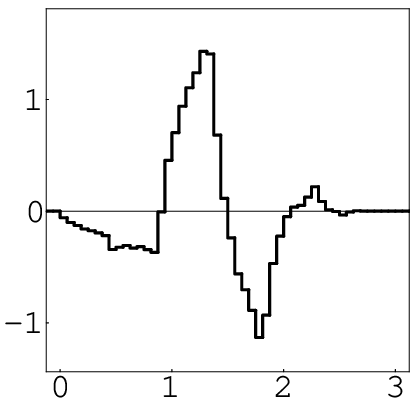}}
\put(0,79){\includegraphics[bb=0 1 118 117,width=79.65bp]{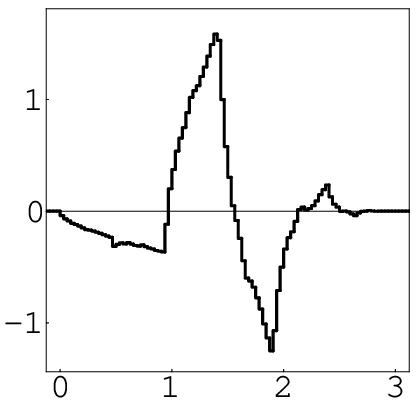}}
\put(120,79){\includegraphics[bb=0 1 118 117,width=79.65bp]{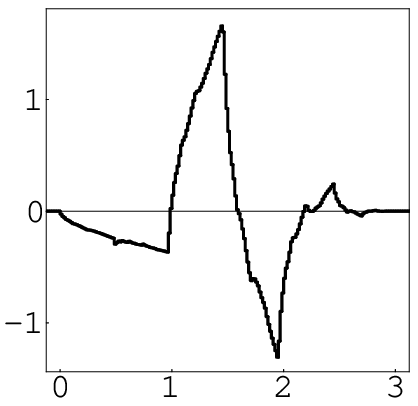}}
\put(298,17){\includegraphics[bb=0 1 180 179,width=121.5bp]{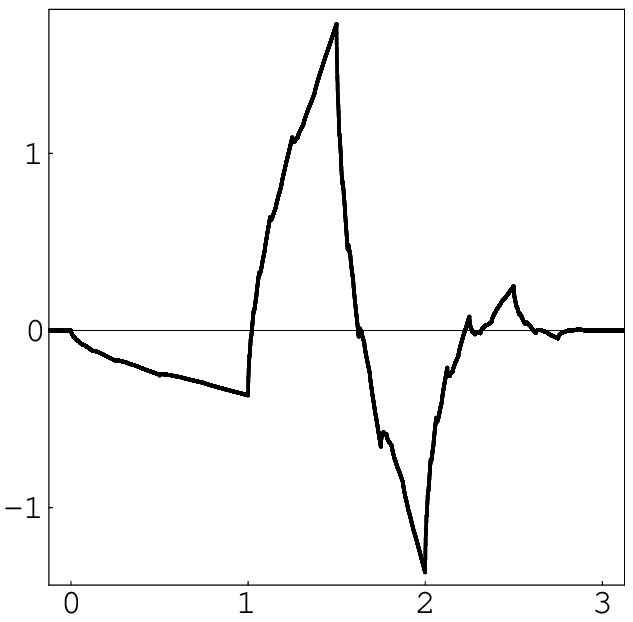}}
\put(269,131){\makebox(0,0)[b]{$\cdots$}}
\put(409,152){\makebox(0,0)[b]{$\psi$}}
\put(392,0){\makebox(0,0)[b]{\small (b): Mother function}}
\end{picture}\caption{Daubechies wavelet functions and series of cascade
approximants}%
\label{FigDaubechies}%
\end{figure}

Since the Haar example is based on the fitting of copies of a fixed
\textquotedblleft box\textquotedblright\ inside an expanded one, it would
almost seem unlikely that the system (\ref{eq3})--(\ref{eq4}) admits finite
sequences $\left(  a_{k}\right)  $ and $\left(  b_{k}\right)  $ such that the
corresponding solutions $\varphi$ and $\psi$ are continuous or differentiable
functions of compact support. The discovery in the mid-1980's of compactly
supported differentiable solutions, see \cite{Dau92}, was paralleled by
applications in seismology, acoustics \cite{EsGa77}, and optics \cite{Mar82},
as discussed in \cite{Mey93a}, and once the solutions were found, other
applications followed at a rapid pace: see, for example, the ten books in
Benedetto's review \cite{Ben00}. It is the solution $\psi$ in (\ref{eq4}) that
the fuss is about, the mother function; the other one, $\varphi$, the father
function, is only there before the birth of the wavelet. The most famous of
them are named after Daubechies, and look like the graphs in Figure
\ref{FigDaubechies}. With the multiresolution idea, we arrive at the closed
subspaces%
\begin{equation}
\mathcal{V}_{j}:=U^{-j}\mathcal{V}_{0},\qquad j\in\mathbb{Z},\label{eq7}%
\end{equation}
as noted in (\ref{eq1})--(\ref{eq2}), where $U$ is some scaling operator.
There are extremely effective iterative algorithms for solving the scaling
identity (\ref{eq3}): see, for example, Example 2.5.3, pp.\ 124--125, of
\cite{BrJo02b}$^*$\footnotetext{$^*$See an implementation of the \textquotedblleft
cascade\textquotedblright\ algorithm using Mathematica, and
a ``cartoon'' of wavelets computed with it, at
\newline
http://www.math.uiowa.edu/\raisebox{-0.75ex}{\symbol{126}}%
jorgen/wavelet\_motions.pdf~.}, \cite{Dau92}, and \cite{StNg96}, and Figure
\ref{FigDaubechies}. A key step in the algorithms involves a clever choice of
the kind of resolution pictured in (\ref{eq12}), but digitally encoded. The
orthogonality relations can be encoded in the numbers $\left(  a_{k}\right)  $
and $\left(  b_{k}\right)  $ of (\ref{eq3})--(\ref{eq4}), and we arrive at the
doubly indexed functions
\begin{equation}
\psi_{j,k}^{{}}\left(  x\right)  :=2^{j/2}\psi\left(  2^{j}x-k\right)  ,\qquad
j,k\in\mathbb{Z}.\label{eq8}%
\end{equation}
It is then not difficult to establish the combined orthogonality relations%
\begin{equation}
\int_{\mathbb{R}}\overline{\psi_{j,k}^{{}}\left(  x\right)  }\,\psi
_{j^{\prime},k^{\prime}}^{{}}\left(  x\right)  ~dx=\left\langle \psi_{j,k}%
^{{}}\mid\psi_{j^{\prime},k^{\prime}}^{{}}\right\rangle =\delta_{j,j^{\prime}%
}^{{}}\delta_{k,k^{\prime}}^{{}}\label{eq9}%
\end{equation}
plus the fact that the functions in (\ref{eq8}) form an orthogonal basis for
$L^{2}\left(  \mathbb{R}\right)  $. This provides a painless representation of
$L^{2}\left(  \mathbb{R}\right)  $-functions
\begin{equation}
f=\sum_{j\in\mathbb{Z}}\sum_{k\in\mathbb{Z}}c_{j,k}^{{}}\psi_{j,k}^{{}%
}\label{eq10}%
\end{equation}
where the coefficients $c_{j,k}$ are%
\begin{equation}
c_{j,k}=\int_{\mathbb{R}}\overline{\psi_{j,k}^{{}}\left(  x\right)
}\,f\left(  x\right)  \,dx=\left\langle \psi_{j,k}^{{}}\mid f\right\rangle
.\label{eq11}%
\end{equation}
What is more significant is that the resolution structure of closed subspaces
of $L^{2}\left(  \mathbb{R}\right)  $%
\begin{equation}
\cdots\subset\mathcal{V}_{-2}\subset\mathcal{V}_{-1}\subset\mathcal{V}%
_{0}\subset\mathcal{V}_{1}\subset\mathcal{V}_{2}\subset\cdots\label{eq12}%
\end{equation}
facilitates powerful algorithms for the representation of the numbers
$c_{j,k}$ in (\ref{eq11}). Amazingly, the two sets of numbers $\left(
a_{k}\right)  $ and $\left(  b_{k}\right)  $ which were used in (\ref{eq3}%
)--(\ref{eq4}), and which produced the magic basis (\ref{eq8}), the wavelets,
are the same magic numbers which encode the quadrature mirror filters of
signal processing of communications engineering. On the face of it, those
signals from communication engineering really seem to be quite unrelated to
the issues from wavelets---the signals are just sequences, time is discrete,
while wavelets concern $L^{2}\left(  \mathbb{R}\right)  $ and problems in
mathematical analysis that are highly non-discrete. Dual filters, or more
generally, subband filters, were invented in engineering well before the
wavelet craze in mathematics of recent decades. These dual filters in
engineering have long been used in technology, even more generally than merely
for the context of quadrature mirror filters (QMF's), and it turns out that
other popular dual wavelet bases for $L^{2}\left(  \mathbb{R}\right)  $ can be
constructed from the more general filter systems; but the best of the wavelet
bases are the ones that yield the strongest form of orthogonality, which is
(\ref{eq9}), and they are the ones that come from the QMF's. The QMF's in turn
are the ones that yield perfect reconstruction of signals that are passed
through filters of the analysis-synthesis algorithms of signal processing.
They are also the algorithms whose iteration corresponds to the resolution
sytems (\ref{eq12}) from wavelet theory.

While Fourier invented his transform for the purpose of solving the heat
equation, i.e., the partial differential equation for heat conduction, the
wavelet transform (\ref{eq10})--(\ref{eq11}) does not diagonalize the
differential operators in the same way. Its effectiveness is more at the level
of computation; it turns integral operators into sparse matrices, i.e.,
matrices which have \textquotedblleft many\textquotedblright\ zeros in the
off-diagonal entry slots. Again, the resolution (\ref{eq12}) is key to how
this matrix encoding is done in practice.

\subsubsection{\label{S1.3.1}Some points of history}

The first wavelet was discovered by
Alfred
Haar%
\index{Haar0@Alfred Haar} long ago, but its use was limited since it was based on step-functions,
and the step-functions jump from one step to the next. The implementation of
Haar%
\index{Haar1@Haar!wavelet@--- wavelet}'s wavelet in the approximation%
\index{approximation!wavelet} problem for continuous functions was
therefore rather bad,
and for differentiable functions it is atrocious,
and so Haar%
\index{Haar0@Alfred Haar}'s method was forgotten for many years. And
yet it had in it the one idea which proved so powerful in the recent rebirth
(since the 1980's) of wavelet analysis%
\index{analysis!wavelet}: the idea of a \emph{multiresolution%
\index{multiresolution|textbf}}.
You see it in its simplest form by noticing that a box function%
\index{Haar1@Haar!scaling function@--- scaling function} $B$ of
(\ref{eqHaarScaling}) may be scaled down by a half such that two copies
$B^{\prime}$ and $B^{\prime\prime}$ of the smaller box then fit precisely inside
$B$. See (\ref{eqHaarScaling}).
\begin{gather}%
\begin{array}
[c]{c}%
\setlength{\unitlength}{30bp}%\settowidth{\qedskip}{$B$}
\begin{picture}%
(3,1.559)(-0.5,-0.309) \put(-0.5,0){\vector(1,0){3}}
\thicklines \put(-0.125,0){\line(1,0){2.25}} \put(0,1){\line(1,0){2}}
\multiput(0,0)(1,0){3}{\line(0,1){1}} \put(1,-0.0425){\makebox(0,0)[t]{$B$}}
\put(0.3657,1.09375){\makebox(0,0)[bl]{$B'$}}
\put(1.3657,1.09375){\makebox(0,0)[bl]{$B''$}}
\put(0.5,0.5){\makebox(0,0){$\varphi$}}
\put(-0.0625,0.0625){\makebox(0,0)[br]{$\scriptstyle 0$}}
\put(1.0625,0.0625){\makebox(0,0)[bl]{$\scriptstyle 1$}}
\put(2.0625,0.0625){\makebox(0,0)[bl]{$\scriptstyle 2$}}
\put(2.4375,-0.125){\makebox(0,0)[t]{$x$}} \end{picture}%
\end{array}
\label{eqHaarScaling}
\\
\begin{array}
[c]{c}%
\setlength{\unitlength}{30bp}\begin{picture}%
(3,1.61)(-0.5,-0.79) \put(-0.5,0){\vector(1,0){3}}
\thicklines \put(-0.125,0){\line(1,0){0.125}} \put(1,0){\line(1,0){1.125}}
\multiput(0,0)(0.5,0){2}{\line(0,1){1}} \multiput(0,1)(0.5,-2){2}%
{\line(1,0){0.5}} \multiput(0.5,-1)(0.5,0){2}{\line(0,1){1}}
\put(0.25,0.5){\makebox(0,0){$\psi$}}
\put(0,-0.0625){\makebox(0,0)[t]{$\scriptstyle 0$}}
\put(0.53125,-0.0625){\makebox(0,0)[tl]{$\scriptstyle \frac{1}{2}$}}
\put(1.0625,-0.0625){\makebox(0,0)[tl]{$\scriptstyle 1$}}
\put(2.4375,-0.125){\makebox(0,0)[t]{$x$}} \end{picture}%
\end{array}
\label{eqHaarWavelet}%
\end{gather}
This process may be continued if you scale by powers of $2$ in both
directions, i.e., by $2^{k}$ for integral $k$, $-\infty<k<\infty$. So for
every $k\in\mathbb{Z}$, there is a finer resolution%
\index{resolution!finer}, and if you take an up-
and a shifted mirror image down-version of the dyadic scaling as in
(\ref{eqHaarWavelet}), and allow all linear combinations, you will notice that
arbitrary functions $f$ on the line $-\infty<x<\infty$, with reasonable
integrability properties, admit a representation
\begin{equation}
f\left(  x\right)  =\sum_{k,n}c_{k,n}\psi\left(  2^{k}x-n\right)  ,
\label{eqIntTut.1}%
\end{equation}
where the summation is over all pairs of integers $k,n\in\mathbb{Z}$, with $k$
representing scaling and $n$ translation. The very simple idea of turning this
construction into a multiresolution (\textquotedblleft
multi\textquotedblright \ for the variety of scales in (\ref{eqIntTut.1}))
leads not only to an algorithm%
\index{algorithm!multiresolution} for the analysis%
\index{analysis!synthesis@--- and synthesis}/synthesis problem,
\begin{equation}
f\left(  x\right)  \longleftrightarrow c_{k,n}, \label{eqIntTut.2}%
\end{equation}
in (\ref{eqIntTut.1}), but also to
a construction of the
single functions $\psi$ which solve the
problem in
(\ref{eqIntTut.1}),
and which can be chosen differentiable, and
yet with support contained in a fixed finite interval. These two features, the
algorithm%
\index{algorithm} and the finite support (called \emph{compact%
\index{compact support|textbf}} support), are crucial
for computations: Computers do algorithm%
\index{algorithm}s, but they do not do infinite
intervals well. Computers do summations and algebra well, but they do not do
integrals and differential equations, unless the calculus problems are
discretized and turned into algorithm%
\index{algorithm}s.

In the discussion to follow,
the multiresolution analysis%
\index{analysis!multiresolution}
viewpoint is dominant, which increases the role of algorithm%
\index{algorithm}s; for example,
the so-called pyramid algorithm%
\index{algorithm!pyramid} for analyzing signals, or shapes, using
wavelets,
is an outgrowth of multiresolution%
\index{multiresolution}s.

Returning to (\ref{eqHaarScaling}) and (\ref{eqHaarWavelet}),
%(see also
%(\ref{eqIntAug12.3})),
we see that the
scaling%
\index{scaling!function@--- function} function
$\varphi$ itself may be expanded in the wavelet%
\index{wavelet!basis@--- basis} basis which is defined from
$\psi$,
and we arrive at the infinite series
\begin{equation}
\varphi \left( x\right) =\sum_{k=1}^{\infty}2^{-k}\psi\left( 2^{-k}x\right) 
\label{eqInfiniteSeries}
\end{equation}
which is pointwise converge%
\index{convergence!pointwise}nt for $x\in\mathbb{R}$.
(It is a special case of the expansion (\ref{eqIntTut.1}) when $f=\varphi$.)
In view of the picture
(\kern2pt\includegraphics[bb=0 13 36 46,width=9bp]{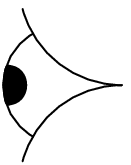}\kern2pt)
below, (\ref{eqInfiniteSeries}) gives
an alternative meaning to the traditional
concept of a \emph{telescoping%
\index{telescoping infinite sum|textbf}} infinite sum.
If, for example,
$0<x<1$, then the representation (\ref{eqInfiniteSeries}) yields
$\varphi \left( x\right) =1=\frac{1}{2}+\left( \frac{1}{2}\right) ^{2}+\cdots$,
while for $1<x<2$, 
$\varphi \left( x\right) =0
=-\frac{1}{2}
+\left( \frac{1}{2}\right) ^{2}
+\left( \frac{1}{2}\right) ^{3}+\cdots $. 
More generally, 
if $n\in\mathbb{N}$, and $2^{n-1}<x<2^{n}$, then
\[
\varphi \left( x\right) =0=-\left( \frac{1}{2}\right) ^{n}
+\sum_{k>n}\left( \frac{1}{2}\right) ^{k}.
\]
So the function $\varphi$ is itself in the space 
$\mathcal{V}_{0}\subset L^{2}\left( \mathbb{R}\right) $,
and $\varphi$
represents the \emph{initial resolution%
\index{resolution!initial|textbf}}. The tail terms
in (\ref{eqInfiniteSeries}) corresponding to
\begin{equation}
\sum_{k>n}2^{-k}\psi\left( 2^{-k}x\right) 
=\frac{1}{2^{n}}\varphi \left( \frac{x}{2^{n}}\right) 
\label{eqTail}
\end{equation}
represent
the \emph{coarser%
\index{resolution!coarser|textbf} resolution}. The
finite sum 
\[
\sum_{k=1}^{n}2^{-k}\psi\left( 2^{-k}x\right) 
\]
represents
the \emph{missing detail%
\index{detail|textbf}} of $\varphi$ as a ``bump signal''.
While the sum on the left-hand side in (\ref{eqTail})
is \emph{infinite}, i.e., the summation index $k$
is in the range $n<k<\infty$, the expression
$2^{-n}\varphi\left( 2^{-n}x\right) $ on the right-hand side is
merely a coarser%
\index{scaling!coarser} scaled version of the
original function $\varphi$ from the subspace
$\mathcal{V}\subset L^{2}\left( \mathbb{R}\right) $
which specifies the initial resolution%
\index{resolution!initial}.
Infinite sums are \emph{analysis%
\index{analysis|textbf} problems} while a
scale operation is a single simple \emph{algorithmic
step}. And so we have encountered a
first (easy) instance of the magic of
a resolution algorithm%
\index{algorithm!resolution}; i.e., an instance
of a transcendental step (the analysis%
\index{analysis} problem) which is
converted into a programmable operation,
here the operation of scaling. (Other
more powerful uses of the scaling operation
may be found in the recent book \cite{Mey98}
\index{Meyer0@Yves Meyer}
by Yves Meyer%
\index{Meyer0@Yves Meyer}, especially Ch.~5,
and \cite{HwMa94}%
\index{Mallat0@St\'{e}phane G. Mallat}.)

The sketch below
allows you to visualize more clearly
this resolution%
\index{resolution} versus detail%
\index{detail}
concept which is so central to the wavelet algorithm%
\index{algorithm!wavelet}s,
also for general wavelets which otherwise may
be computationally more difficult than the Haar%
\index{Haar1@Haar!wavelet@--- wavelet}
wavelet.
\begin{gather*}
\setlength{\unitlength}{0.675bp}
\begin{picture}(384,164)
\put(0,0){\includegraphics[bb=8 0 376 164,width=248.4bp]{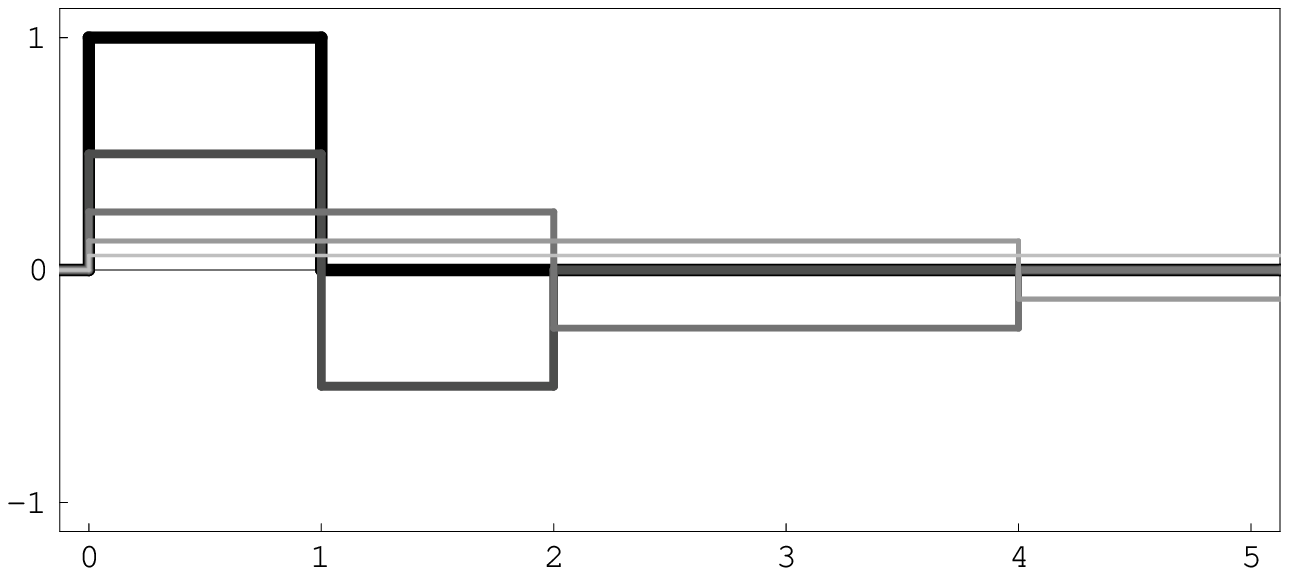}}
\put(98,146){$\varphi\!\left( x\right) $}
\put(104,39){$\frac{1}{2}\psi\!\left( \frac{x}{2}\right) $}
\put(208,56){$\frac{1}{4}\psi\!\left( \frac{x}{4}\right) $}
\put(312,65){$\frac{1}{8}\psi\!\left( \frac{x}{8}\right) $}
\put(308,102){$\frac{1}{16}\psi\!\left( \frac{x}{16}\right) $}
\put(413,73){\includegraphics[bb=0 0 36 46,width=18bp]{eye.eps}}
\end{picture}
\\
\text{\footnotesize The wavelet%
\index{wavelet!decomposition@--- decomposition} decomposition of Haar%
\index{Haar1@Haar!scaling function@--- scaling function}'s bump function $\varphi
$ in (\ref{eqHaarScaling}) and (\ref{eqInfiniteSeries})}
\end{gather*}
Using
the sketch we see for example
that the simple step function
\begin{gather}
f\left( x\right) =a\varphi \left( x\right) +b\varphi \left( x-1\right) 
=a\chi_{\left\lbrack 0,1\right\rangle}^{{}}\left( x\right) 
+b\chi_{\left\lbrack 1,2\right\rangle}^{{}}\left( x\right) 
\label{eqTwoStep}\\
\setlength{\unitlength}{0.675bp}
\begin{picture}(352,160)(0,4)
\put(0,0){\includegraphics[bb=16 0 368 164,width=237.6bp]{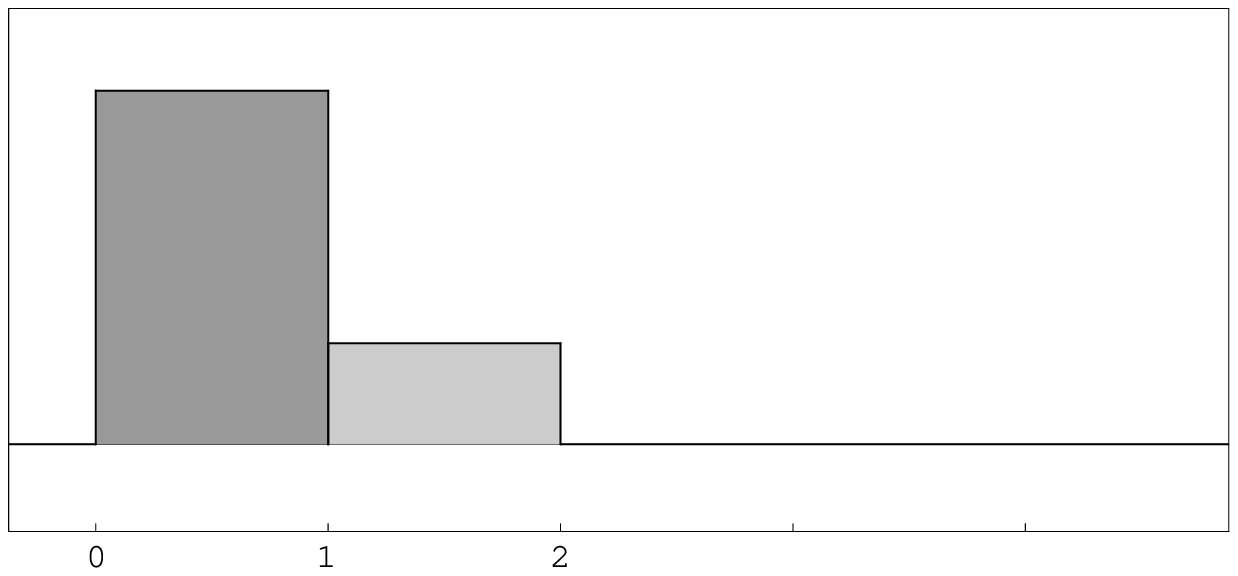}}
\put(14,87){$a$}
\put(96,49){$b$}
\end{picture}
\nonumber
\end{gather}
has the wavelet%
\index{wavelet!decomposition@--- decomposition} decomposition into a sum of
a
\emph{coarser%
\index{resolution!coarser|textbf} resolution} and an
\emph{intermediate detail%
\index{detail|textbf}} as follows:
\begin{equation}
f\left( x\right) 
=\underset{\text{intermediate detail%
\index{detail}}}
{\underbrace{\frac{a-b}{2}\psi\left( \frac{x}{2}\right) }}
+\underset{\text{coarser%
\index{resolution!coarser} version}}
{\underbrace{\frac{a+b}{2}\varphi \left( \frac{x}{2}\right) }},
\qquad x\in\mathbb{R}.
\label{eqDetail}
\end{equation}
Thus
the details are measured as differences. This is a
general feature
that is
valid for other functions and other wavelet resolutions.
See, for instance, \S\ \ref{Exe1Sep7-4} below.

\subsubsection{\label{S1.3.2}Some early applications}

While the Haar wavelet is built from flat pieces, and the orthogonality
properties amount to a visual tiling of the graphs of the two functions
$\varphi$ and $\psi$, this is not so for the Daubechies wavelet nor the other
compactly supported smooth wavelets. By the Balian--Low theorem \cite{Dau92},
a time-frequency wavelet cannot be simultaneously localized in the two dual
variables: if $\psi$ is a time-frequency Gabor wavelet,
then the two quantities $\int_{\mathbb{R}}\left\vert x\psi\left(  x\right)
\right\vert ^{2}\,dx$ and $\int_{\mathbb{R}}\left\vert t\hat{\psi}\left(
t\right)  \right\vert ^{2}\,dt$ cannot both be finite. Since $\left(
\frac{d\psi}{dx}\right)  \sphat\left(  t\right)  =it\hat{\psi}\left(
t\right)  $, this amounts to poor differentiability properties of
well-localized Gabor wavelets, i.e., wavelets built using the two operations
translation and frequency modulation over a lattice.

But with the multiresolution viewpoint, we can understand the first of
Daubechies's scaling functions as a one-sided differentiable solution
$\varphi$ to%
\begin{equation}
\varphi\left(  x\right)  =h_{0}\varphi\left(  2x\right)  +h_{1}\varphi\left(
2x-1\right)  +h_{2}\varphi\left(  2x-2\right)  +h_{3}\varphi\left(
2x-3\right)  ,\label{eq1.3.1}%
\end{equation}
where the four real coefficients satisfy%
\begin{equation}
\left. 
\begin{aligned}
h_{0}+h_{1}+h_{2}+h_{3}  & =2,\\
h_{3}-h_{2}+h_{1}-h_{0}  & =0,\\
h_{3}-2h_{2}+3h_{1}-4h_{0}  & =0,\\
h_{1}h_{3}+h_{0}h_{2}  & =0.
\end{aligned}
\right\}\label{eq1.3.2}%
\end{equation}
The system (\ref{eq1.3.2}) is easily solved:%
\begin{equation}
\left.
\begin{aligned}
4h_{0}  & =1+\sqrt{3}, & 4h_{2} & =3-\sqrt{3},\\
4h_{1}  & =3+\sqrt{3}, & 4h_{3} & =1-\sqrt{3},
\end{aligned}
\right\}\label{eq1.3.3}%
\end{equation}
and Daubechies showed that (\ref{eq1.3.1}) has a solution $\varphi$ which is
supported in the interval $\left[  0,3\right]  ,$ is one-sided differentiable,
and satisfies the conditions%
\begin{equation}
\int_{\mathbb{R}}\varphi\left(  x\right)  \,dx=1,\quad\int_{\mathbb{R}}%
\psi\left(  x\right)  \,dx=0,\text{\quad and\quad}\int_{\mathbb{R}}%
x\psi\left(  x\right)  \,dx=0.\label{eq1.3.4}%
\end{equation}

The first applications served as motivating ideas as well: optics, seismic
measurements, dynamics, turbulence, data compression; see the book
\cite{KaLe95} 
%has somewhat of the same flavor, but it is written for
%mathematicians. It has formulas,
%and in addition a lot of excellent writing. Like
%\cite{Hub98}, it stresses the intuitive ideas behind the formulas. 
Actually,
it is two books: the first one (primarily by Kahane) is classical Fourier
analysis%
\index{analysis!Fourier}, and the second one (primarily by P.-G. Lemari\'{e}-Rieusset) is the
wavelet book. It will help you, among other things, to get a better feel for
the French connection, the Belgian connection, and the diverse and early
impulses from applications in the subject. Enjoy!

For a list of more recent applications we recommend \cite{Mey00}.

\section{\label{S2}Signal processing}

\setcounter{equation}{0}
\renewcommand{\theequation}{\thesection.\arabic{equation}}
If we idealize and view time as discrete, a copy of $\mathbb{Z}$, then a
signal is a sequence $\left(  \xi_{n}\right)  _{n\in\mathbb{Z}}$ of numbers. A
filter is an operator which calculates weighted averages%
\begin{equation}
\left(  \xi_{n}\right)  \longmapsto\sum_{k\in\mathbb{Z}}a_{k}\xi
_{n-k}.\label{eq2.1}%
\end{equation}
But working instead with functions of $z\in\mathbb{T}$, this is
multiplication, $f\left(  z\right)  \mapsto m\left(  z\right)  f\left(
z\right)  $, where $m\left(  z\right)  =\sum_{k\in\mathbb{Z}}a_{k}z^{k}$ and
$f\left(  z\right)  =\sum_{k\in\mathbb{Z}}\xi_{k}z^{k}$ are the usual Fourier
representation of the corresponding generating functions. Similarly,
down-sampling
$\vcenter{\hbox{\Large $\bigcirc \llap{$\vcenter{\hbox{\normalsize $\scriptscriptstyle N$}}\mkern-1mu\vcenter{\hbox{\normalsize $\scriptscriptstyle\downarrow$}}\mkern3.5mu$}$}}$
and up-sampling
$\vcenter{\hbox{\Large $\bigcirc \llap{$\vcenter{\hbox{\normalsize $\scriptscriptstyle N$}}\mkern-1mu\vcenter{\hbox{\normalsize $\scriptscriptstyle\uparrow$}}\mkern3.5mu$}$}}$
as operators on sequences take the form%
\begin{equation}
f\longmapsto\frac{1}{N}\sum_{w\in\mathbb{T},\,w^{N}=z}f\left(  w\right)
\label{eq2.2}%
\end{equation}
and%
\begin{equation}
f\longmapsto f\left(  z^{N}\right)  .\label{eq2.3}%
\end{equation}
Since the operators
$\vcenter{\hbox{\Large $\bigcirc \llap{$\vcenter{\hbox{\normalsize $\scriptscriptstyle N$}}\mkern-1mu\vcenter{\hbox{\normalsize $\scriptscriptstyle\downarrow$}}\mkern3.5mu$}$}}$
and
$\vcenter{\hbox{\Large $\bigcirc \llap{$\vcenter{\hbox{\normalsize $\scriptscriptstyle N$}}\mkern-1mu\vcenter{\hbox{\normalsize $\scriptscriptstyle\uparrow$}}\mkern3.5mu$}$}}$
are clearly dual to one another on the Hilbert space $\ell^{2}\left(
\mathbb{Z}\right)  $ of sequences (i.e., time-signals), we get the
corresponding duality for $L^{2}\left(  \mathbb{T}\right)  $, i.e.,%
\begin{equation}
\int_{\mathbb{T}}f\left(  z^{N}\right)  g\left(  z\right)  \,d\mu\left(
z\right)  =\int_{\mathbb{T}}f\left(  z\right)  \frac{1}{N}\sum_{w^{N}%
=z}g\left(  w\right)  \,d\mu\left(  z\right)  ,\label{eq2.4}%
\end{equation}
where $\mu$ denotes the normalized Haar measure on $\mathbb{T}$, or
equivalently the following identity for $2\pi$-periodic functions:%
\begin{equation}
\int_{0}^{2\pi}f\left(  N\theta\right)  g\left(  \theta\right)  \,d\theta
=\int_{0}^{2\pi}f\left(  \theta\right)  \frac{1}{N}\sum_{k=0}^{N-1}g\left(
\frac{\theta+k\cdot2\pi}{N}\right)  \,d\theta.\label{eq2.5}%
\end{equation}

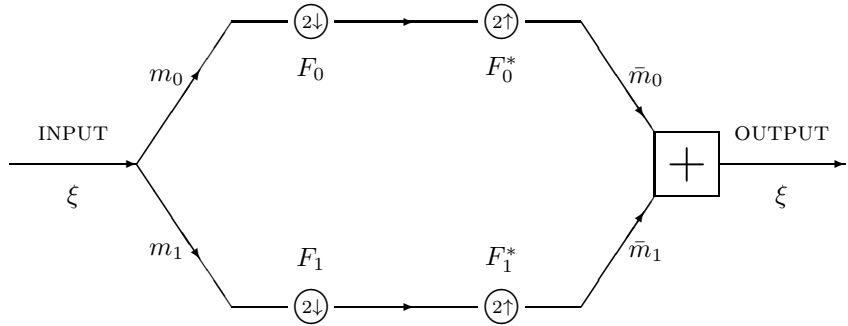
\begin{figure}[ptb]
\begin{center}
\setlength{\unitlength}{2pc} \begin{picture}(13.166667,5)(1.5,-2.5)
\put(1.5,0){\vector(1,0){2}}
\put(1.5,0){\makebox(2,1){\textsc{input}}}
\put(1.5,-1){\makebox(2,1){$\xi$}}
\put(3.5,0){\line(2,3){1.5}}
\put(3.5,0){\vector(2,3){1}}
\put(5,2.25){\line(1,0){0.875}}
\put(3.5,0){\line(2,-3){1.5}}
\put(3.5,0){\vector(2,-3){1}}
\put(5,-2.25){\line(1,0){0.875}}
\put(5.75,1.75){\makebox(1,1){$\vcenter{\hbox{\Large $\bigcirc \llap{$\vcenter{\hbox{\normalsize $\scriptstyle 2$}}\mkern-0.1mu\vcenter{\hbox{\normalsize $\scriptstyle\downarrow$}}\mkern3.5mu$}$}}$}}
\put(3.25,0.875){\makebox(1,1)[r]{$m_{0}$}}
\put(3.25,-1.875){\makebox(1,1)[r]{$m_{1}$}}
\put(5.75,-2.75){\makebox(1,1){$\vcenter{\hbox{\Large $\bigcirc \llap{$\vcenter{\hbox{\normalsize $\scriptstyle 2$}}\mkern-0.1mu\vcenter{\hbox{\normalsize $\scriptstyle\downarrow$}}\mkern3.5mu$}$}}$}}
\put(6.625,2.25){\line(1,0){2.25}}
\put(6.625,2.25){\vector(1,0){1.25}}
\put(6.625,-2.25){\line(1,0){2.25}}
\put(6.625,-2.25){\vector(1,0){1.25}}
\put(5.75,1){\makebox(1,1){$F_{0}$}}
\put(5.75,-2){\makebox(1,1){$F_{1}$}}
\put(12.666667,0){\vector(1,0){2}}
\put(12.666667,0){\makebox(2,1){\textsc{output}}}
\put(12.666667,-1){\makebox(2,1){$\xi$}}
\put(10.5,2.25){\line(2,-3){1.166667}}
\put(10.5,2.25){\vector(2,-3){1}}
\put(9.625,2.25){\line(1,0){0.875}}
\put(10.5,-2.25){\line(2,3){1.166667}}
\put(10.5,-2.25){\vector(2,3){1}}
\put(9.625,-2.25){\line(1,0){0.875}}
\put(8.75,1.75){\makebox(1,1){$\vcenter{\hbox{\Large $\bigcirc \llap{$\vcenter{\hbox{\normalsize $\scriptstyle 2$}}\mkern-0.1mu\vcenter{\hbox{\normalsize $\scriptstyle\uparrow$}}\mkern3.5mu$}$}}$}}
\put(11.25,0.875){\makebox(1,1)[l]{$\bar{m}_{0}$}}
\put(11.25,-1.875){\makebox(1,1)[l]{$\bar{m}_{1}$}}
\put(8.75,-2.75){\makebox(1,1){$\vcenter{\hbox{\Large $\bigcirc \llap{$\vcenter{\hbox{\normalsize $\scriptstyle 2$}}\mkern-0.1mu\vcenter{\hbox{\normalsize $\scriptstyle\uparrow$}}\mkern3.5mu$}$}}$}}
\put(8.75,1){\makebox(1,1){$F_{0}^{*}$}}
\put(8.75,-2){\makebox(1,1){$F_{1}^{*}$}}
\put(11.666667,-0.5){\framebox(1,1){\huge$+$}}
\end{picture}
\end{center}
\caption{Perfect reconstruction of signals}%
\label{FigPerfectReconstruction}%
\end{figure}

Quadrature mirror filters with $N$ frequency subbands $m_{0},m_{1}%
,\dots,m_{N-1}$ give perfect reconstruction when signals are analyzed into
subbands and then reconstructed via the up-sampling and corresponding dual
filters. In engineering formalism this is expressed in the diagram in
Fig.\ \ref{FigPerfectReconstruction}, for $N=2$, and $m_{0}$, resp.\ $m_{1}$,
are called low-pass, resp.\ high-pass, filters. In operator language, this
takes the form%
\[
F_{0}^{\ast}F_{0}+F_{1}^{\ast}F_{1}=I,
\]
where $F_{0}$ and $F_{1}$ are the operators in
Fig.\ \ref{FigPerfectReconstruction}, with dual operators $F_{0}^{\ast}$ and
$F_{1}^{\ast}$. The quadrature conditions may be expressed as%
\begin{equation}
F_{0}F_{0}^{\ast}=F_{1}F_{1}^{\ast}=I\label{eq2.6}%
\end{equation}
and%
\begin{equation}
F_{0}F_{1}^{\ast}=F_{1}F_{0}^{\ast}=0.\label{eq2.7}%
\end{equation}
In operator theory there is tradition for working instead with the operators
$S_{j}:=F_{j}^{\ast}$. When viewed as operators on $L^{2}\left(
\mathbb{T}\right)  $ they are therefore isometries with orthogonal ranges, and
they satisfy%
\begin{equation}
\sum_{j=0}^{N-1}S_{j}S_{j}^{\ast}=I\label{eq2.8}%
\end{equation}
with $I$ now representing the identity operator acting on $L^{2}\left(
\mathbb{T}\right)  $. The relations on the $S_{j}$-operators are known as the
Cuntz relations because of their use in $C^{\ast}$-algebra theory; see
\cite{Cun77}. In the present application they take the form%
\begin{equation}
\left(  S_{j}f\right)  \left(  z\right)  =m_{j}\left(  z\right)  f\left(
z^{N}\right)  ,\qquad f\in L^{2}\left(  \mathbb{T}\right)  ,\label{eq2.9}%
\end{equation}
and%
\begin{equation}
\left(  S_{j}^{\ast}f\right)  \left(  z\right)  =\frac{1}{N}\sum_{w^{N}%
=z}\overline{m_{j}\left(  w\right)  }\,f\left(  w\right)  ,\label{eq2.10}%
\end{equation}
and the Cuntz relations are equivalent to the conditions%
\begin{equation}
\sum_{w^{N}=z}\left\vert m_{j}\left(  w\right)  \right\vert ^{2}%
=N\label{eq2.11}%
\end{equation}
and%
\begin{equation}
\sum_{w^{N}=z}\overline{m_{j}\left(  w\right)  }\,m_{k}\left(  w\right)
=0\text{\qquad for all }z\in\mathbb{T}\text{ and }j\neq k.\label{eq2.12}%
\end{equation}
The last conditions are known in engineering as the quadrature conditions for
the subband filters $m_{0},m_{1},\dots,m_{N-1}$, with $m_{0}$ denoting the
low-pass filter. The low-pass and band-pass conditions on the functions
$m_{j}$ are perhaps more familiar in the additive notation given by the
substitution $z:=e^{-i\theta}$. Then the functions $m_{j}$ are viewed as
$2\pi$-periodic, and
\[
m_{j}\left(  j\cdot\frac{2\pi}{N}\right)  =\sqrt{N},
\]
while%
\[
m_{j}\left(  k\cdot\frac{2\pi}{N}\right)  =0\text{\qquad for }j\neq k,
\]
with both of the indices $j$, $k$ ranging over $0,1,\dots,N-1$.
\renewcommand{\theequation}{\thesubsection.\arabic{equation}}

\subsection{\label{S2.1}Filters in communications engineering}

The coefficients of the functions $m_{j}\left(  \,\cdot\,\right)  $ are called
\emph{impulse response coefficients} in communications engineering, and when
used in wavelets and in subdivision algorithms, they are called \emph{masking
coefficients}. In the finite case, the $m_{j}\left(  \,\cdot\,\right)  $'s are
also called FIR for finite impulse response. The model illustrated in
Fig.\ \ref{FigPerfectReconstruction} is used in filter design in either
hardware or software:

\begin{enumerate}
\renewcommand{\theenumi}{\arabic{enumi}} \renewcommand{\labelenumi}{$\lbrack\!\lbrack\theenumi\rbrack\!\rbrack$}

\item \label{S2.1(1)}Try filters $m_{0}$, $m_{1}$ in
Fig.\ \ref{FigPerfectReconstruction}, and approximate the output to the input;

\item \label{S2.1(2)}Choose a specific structure in which the filter will be
realized and then quantize the coefficients, length and numerical values;

\item \label{S2.1(3)}Verify by simulation that the resulting design meets
given performance specifications.
\end{enumerate}

Once filters are constructed, we saw that they are also providing us with
wavelet algorithms. When the steps of Fig.\ \ref{FigPerfectReconstruction} are
iterated, we arrive at wavelet subdivision algorithms. Relative to a given
resolution (pictured as a closed subspace $\mathcal{V}_{1}$, say, in
$L^{2}\left(  \mathbb{R}\right)  $), signals, i.e., functions in $L^{2}\left(
\mathbb{R}\right)  $, decompose into coarser ones and intermediate details.
Relative to the subspaces $\mathcal{W}_{0}$ and $\mathcal{V}_{1}$, this
amounts to
\begin{equation}
\underset{%
\begin{array}
[c]{c}%
\uparrow\\
\text{given}\\
\text{resolution}%
\end{array}
}{{}_{\phantom{1}}\mathcal{V}_{1}}=\underset{%
\begin{array}
[c]{c}%
\uparrow\\
\text{coarser}\\
\text{resolution}%
\end{array}
}{{}_{\phantom{0}}\mathcal{V}_{0}}+\underset{%
\begin{array}
[c]{c}%
\uparrow\\
\text{intermediate}\\
\text{detail}%
\end{array}
}{{}_{\phantom{0.}}\mathcal{W}_{0}.}\label{eq2.1.1}%
\end{equation}
Ideally, we wish the decomposition in (\ref{eq2.1.1}) to be orthogonal in the
sense that%
\begin{equation}
\left\langle \,f\mid g\,\right\rangle =0\text{\qquad for all }f\in
\mathcal{V}_{0}\text{ and all }g\in\mathcal{W}_{0}.\label{eq2.1.2}%
\end{equation}
Since the subdivisions involve translations by discrete steps, we specialize
the resolution such that both of the spaces $\mathcal{V}_{0}$ and
$\mathcal{W}_{0}$ are invariant under translations by points in $\mathbb{Z}$,
i.e., such that%
\begin{equation}
T\colon f\longmapsto f\left(  \,\cdot\,-1\right)  \label{eq2.1.3}%
\end{equation}
leaves both of the subspaces $\mathcal{V}_{0}$ and $\mathcal{W}_{0}$
invariant. The \emph{multiresolution analysis} case (MRA) corresponds to the
setup when $\mathcal{V}_{0}$ is singly generated, i.e., there is a function
$\varphi\in\mathcal{V}_{0}$ such that the closed linear span of%
\begin{equation}
T_{n}\varphi\left(  \,\cdot\,\right)  =\varphi\left(  \,\cdot\,-n\right)
,\qquad n\in\mathbb{Z},\label{eq2.1.4}%
\end{equation}
is all of $\mathcal{V}_{0}$. If $N=2$, then there is then also a $\psi
\in\mathcal{W}_{0}$ such that the closed linear span of $\left\{
\,\psi\left(  \,\cdot\,-n\right)  :n\in\mathbb{Z}\,\right\}  $ is all of
$\mathcal{W}_{0}$. If $N>2$, we may need functions $\psi_{1},\dots,\psi_{N-1}$
in $\mathcal{W}_{0}$ such that $\left\{  \,\psi_{i}\left(  \,\cdot\,-n\right)
:i=1,\dots,N-1,\;n\in\mathbb{Z}\,\right\}  $ has a closed span equal to
$\mathcal{W}_{0}$.

\subsection{\label{S2.2}Algorithms for signals and for wavelets}

\label{Exe1Sep7-4}
\textbf{The pyramid algorithm%
\index{algorithm!pyramid} and the Cuntz%
\index{Cuntz relations} relations.}
Since the two Hilbert%
\index{Hilbert!space@--- space} spaces $L^{2}\left(  \mathbb{T}\right)  $ and $\ell
^{2}\left(  \mathbb{Z}\right)  $ are isomorphic via the Fourier%
\index{Fourier1@Fourier!series@--- series} series
representation, 
%see \textup{(\ref{eqIntJul17.21})\textendash
%(\ref{eqIntJul17.23})}, 
it follows that the system $\left\{  S_{i}\right\}
_{i=0}^{1}$ 
%of \textup{(\ref{eqExe1-3(2)})} in Exercise \textup{\ref{Exe1-3}}
%above 
is equivalent to a system $\left\{  
\smash{\hat{S}_{i}}
\right\}  _{i=0}^{1}$
acting on $\ell^{2}\left(  \mathbb{Z}\right)  $. Specifically, $\left(
S_{i}f\right)  \sphat=\hat{S}_{i}\hat{f}$, $i=0,1$, where $\hat{f}\left(
n\right)  :=\int_{\mathbb{T}}z^{-n}f\left(  z\right) \,d\mu\left(  z\right)
$.
%, satisfy the identities \textup{(\ref{Exe1-3-a})\textendash (\ref{Exe1-3-b}%
%)} in Exercise \textup{\ref{Exe1-3}}. 
For $c:=\left(  c_{n}\right)
_{n\in\mathbb{Z}}$ in $\ell^{2}\left(  \mathbb{Z}\right)  $, and functions $f$
on $\mathbb{R}$, set
\begin{align*}
f_{-1}\left(  x\right)  &:=\left( Uf\right) \left( x\right) 
=2^{-\frac{1}{2}}f\left(
\frac{x}{2}\right)  \text{,
and}\\
\left(  c\ast f\right)  \left(  x\right)
&:=\sum_{n\in\mathbb{Z}}c_{n}f\left(  x-n\right)  .
\end{align*}
For the present,
let $\left\{
m_{i}\right\}  _{i=0}^{1}$ be
the
low-pass and high-pass wavelet filter%
\index{filter!wavelet}%
\index{filter!high-pass}%
\index{filter!low-pass}s,
%given in \textup{(\ref{eqLowPassFilter})} and
%\textup{(\ref{eqHighPassFilter})}, 
and let
$\varphi$, $\psi$ be the corresponding
scaling%
\index{scaling!function@--- function} function, resp., wavelet%
\index{wavelet!function@--- function} function,
also called  father function, resp., mother function.
%In Chapter \textup{\ref{ChHom}} we will repeat this exercise for
%the more general filter%
%\index{filter}s and scaling%
%\index{scaling!function@--- function}/wavelet%
%\index{wavelet!function@--- function}
%functions introduced in Theorem \ref{ThmMul.1}. See
%Exercise \textup{\ref{Exe2Oct5-21}}.
Now introduce the corresponding
operators $S_{i}$ 
%of \textup{(\ref{eqExe1-3(2)})}, 
and their cousins $\hat
{S}_{i}$. The adjoint%
\index{operator!adjoint}s $\hat{S}_{i}^{\ast}$ are
also
called \emph{filters}.

%\begin{enumerate}
%\item \label{Exe1Sep7-4(1)}
%Show that%
Then
\begin{equation}
c\ast\varphi=\underset{\text{coarser%
\index{resolution!coarser} resolution}}{\underbrace{\left(  \left(
\hat{S}_{0}^{\ast}c\right)  \ast\varphi\right)  _{-1}}}+\underset
{\text{detail%
\index{detail}}}{\underbrace{\left(  \left(  \hat{S}_{1}^{\ast}c\right)
\ast\psi\right)  _{-1}}}\text{\qquad for all }c\in\ell^{2}\left(
\mathbb{Z}\right)  .\label{eqExe1Sep7-4(2)}%
\end{equation}
%
%\emph{Hint:}
%If you compute \textup{(\ref{Exe1Sep7-4(1)})} directly
%you will probably get into a complicated mess.
%If so, try doing
%\textup{(\ref{Exe1Feb8-11(2)})},
%\textup{(\ref{Exe1Feb8-11(3)})}, and
%\textup{(\ref{Exe1Feb8-11(4)})} first.
%To get a better idea of what \textup{(\ref{Exe1Sep7-4(1)})} says,
%look at
%\textup{(\ref{Exe1Feb8-11(5)})} and
%\textup{(\ref{Exe1Sep7-4(2)})}.
%
%\item 
Define $W\colon\ell^{2}\rightarrow\ell^{2}$ by
\begin{equation}
W\left( c\right) \left( x\right) 
=\left(  c\ast \varphi\right)  \left(  x\right)
=\sum_{n\in\mathbb{Z}}c_{n}\varphi\left(  x-n\right)  .
\label{Exe1Feb8-11(2)}
\end{equation}
%Show that 
Then
$W$ maps $\ell^{2}$ isometrically onto $\mathcal{V}_{0}$
%\textup{(}defined before \textup{(\ref{eqIntAug12.1}))} 
in the orthogonal case
and
\[
W\hat{S}_{0}=UW.
\]
%
%\item \label{Exe1Feb8-11(3)}
%Use \textup{(\ref{Exe1Feb8-11(2)})} to show that
Further
\[
W\hat{S}_{0}
\hat{S}_{0}^{\ast}c
=\left(
\hat{S}_{0}^{\ast}c\ast \varphi\right) _{-1}.
\]
%
%\item \label{Exe1Feb8-11(4)}
%%\addtolength{\textheight}{-0.25\baselineskip}\label{textheightshavee6}%
Embedding $\ell^{2}$ into $\ell^{2}\oplus\ell^{2}$ as $\ell^{2}\oplus 0$,
extend $W$ to $\ell^{2}\oplus\ell^{2}$ by putting
\[
W\left( c\oplus d\right) 
=c\ast \varphi
+d\ast \psi .
\]
%Show that 
Then
the extended $W$ maps $\ell^{2}\oplus\ell^{2}$
isometrically onto $U^{-1}\mathcal{V}_{0}$ and 
%that
\[
W\left( \hat{S}_{0}c+\hat{S}_{1}d\right) 
=UW\left( c\oplus d\right) 
\]
for all $c,d\in\ell^{2}$, where the left $W$ is the one from
\textup{(\ref{Exe1Feb8-11(2)})} and the right is the extension of
$W$ to $\ell^{2}\oplus\ell^{2}$.

At this point you can use 
$1_{\ell^{2}}=\hat{S}_{0}\hat{S}_{0}^{\ast}+\hat{S}_{1}\hat{S}_{1}^{\ast}$
to show \textup{(\ref{eqExe1Sep7-4(2)})}. Note that if $c_{0}=a$ and $c_{1}=b$ and
$c_{i}=0$ for other $i$, the formula \textup{(\ref{eqExe1Sep7-4(2)})} reduces to
\textup{(\ref{eqDetail})}.

%\item \label{Exe1Feb8-11(5)}
%Show that the identity 
The subdivision relations
\textup{(\ref{eqExe1Sep7-4(2)})} 
%is 
are
equivalent
to the system
\begin{align}
\sqrt{2}\varphi\left( 2x\right) 
&=\sum_{k\in\mathbb{Z}}\bar{a}_{2k}\varphi\left( x+k\right) 
+\sum_{k\in\mathbb{Z}}\bar{b}_{2k}\psi\left( x+k\right) ,\label{eqExe1Feb8-11(5).1}
\\
\sqrt{2}\varphi\left( 2x-1\right) 
&=\sum_{k\in\mathbb{Z}}\bar{a}_{2k+1}\varphi\left( x+k\right) 
+\sum_{k\in\mathbb{Z}}\bar{b}_{2k+1}\psi\left( x+k\right) ,\label{eqExe1Feb8-11(5).2}
\end{align}
where the coefficients $a_{n}$, $b_{n}$ are those of the quantum
wavelet algorithm%
\index{algorithm!quantum wavelet}, i.e.,
the coefficients in
the ``large'' unitary matrix%
\index{matrix!unitary} 
%listed above 
\textup{(\ref{eqBigWaveletMatrix})}.
Thus the quantum%
\index{algorithm!quantum}
algorithm does the wavelet%
\index{wavelet!decomposition@--- decomposition} decomposition within
a fixed resolution%
\index{resolution!subspace@--- subspace} subspace.

%\item \label{Exe1Sep7-4(2)}
The scaling%
\index{scaling!function@--- function} function $\varphi$ defines a resolution%
\index{resolution!subspace@--- subspace} subspace
$\mathcal{V}_{0}\subset L^{2}\left(  \mathbb{R}\right)  $.
%by the method described above \textup{(\ref{eqIntAug12.1})}. Show that
Then
\textup{(\ref{eqExe1Sep7-4(2)})}, 
or equivalently \textup{(\ref{eqExe1Feb8-11(5).1})--(\ref{eqExe1Feb8-11(5).2})}, 
represents the orthogonal%
\index{orthogonal!decomposition@--- decomposition} decomposition of
functions in $\mathcal{V}_{0}$
into an orthogonal sum of a function with
coarser%
\index{resolution!coarser} resolution and a function in the intermediate detail%
\index{detail!subspace@--- subspace} subspace.

Let $m_{0}$, $m_{1}$ be a dyadic wavelet%
\index{wavelet!dyadic} filter%
\index{filter!wavelet}, and let
$\mathbb{T}\ni z\mapsto A\left( z\right) 
\in\mathrm{U}_{2}\left( \mathbb{C}\right) $
be the corresponding matrix%
\index{matrix!function@--- function} function,
$A_{i,j}\left( z\right) 
=\frac{1}{2}\sum_{w^{2}=z}w^{-j}m_{i}\left( w\right) $.
If the low-pass filter%
\index{filter!low-pass} $m_{0}\left( z\right) =a_{0}+a_{1}z+\dots +a_{2n+1}z^{2n+1}$,
then a choice for $m_{1}\left( z\right) =\sum_{k=0}^{2n+1}b_{k}z^{k}$
is $b_{k}=\left( -1\right) ^{k}\bar{a}_{2n+1-k}$. We
then have $A\left( z\right) =\sum_{k=0}^{n}A_{k}z^{k}$ where
$A_{k}=\left( 
\begin{matrix}
a_{2k} & a_{2k+1} \\
b_{2k} & b_{2k+1}
\end{matrix}
\right) $, and the
following $2^{n+2}\times 2^{n+2}$ scalar matrix can be checked to be unitary%
\index{matrix!unitary}:
\begin{equation}
\left( 
\begin{array}{c|c|c|c|c|c|c|c|c|c|c|c}
\begin{matrix}
a_{1}\\
b_{1}
\end{matrix}
 & A_{1} & A_{2} & \cdots & A_{n-1} & A_{n} & 0 &  & \cdots 
 &  & 0 & 
\begin{matrix}
a_{0}\\
b_{0}
\end{matrix}
 \\\hline
\begin{matrix}
0\\
0
\end{matrix}
 & A_{0} & A_{1} & \cdots & A_{n-2} & A_{n-1} & A_{n} & 0 &  & \cdots 
 & 0 & \begin{matrix}
0\\
0
\end{matrix}
 \\\hline
\begin{matrix}
0\\
0
\end{matrix}
 & 0 & A_{0} & \cdots & A_{n-3} & A_{n-2} & A_{n-1} & A_{n} & 0 & \cdots 
 & 0 & \begin{matrix}
0\\
0
\end{matrix}
 \\\hline
\begin{matrix}
0\\
0
\end{matrix}
 &  &  &  &  &  &  &  &  &  &  & \begin{matrix}
0\\
0
\end{matrix}
 \\\hline
\vdots
 &  &  &  & \ddots &  &  &  &  & \ddots &  & \vdots \\\hline
\begin{matrix}
0\\
0
\end{matrix}
 &  &  &  &  &  &  &  &  &  &  & \begin{matrix}
0\\
0
\end{matrix}
 \\\hline
\begin{matrix}
a_{2n+1}\\
b_{2n+1}
\end{matrix}
 & 0 &  & \cdots 
 &  & 0 & A_{0} & A_{1} & A_{2} & \cdots  & A_{n-1} & \begin{matrix}
a_{2n}\\
b_{2n}
\end{matrix}
 \\\hline
\begin{matrix}
a_{2n-1}\\
b_{2n-1}
\end{matrix}
 & A_{n} & 0 &  & \cdots 
 &  & 0 & A_{0} & A_{1} & \cdots & A_{n-2} & \begin{matrix}
a_{2n-2}\\
b_{2n-2}
\end{matrix}
 \\\hline
\begin{matrix}
a_{2n-3}\\
b_{2n-3}
\end{matrix}
 & A_{n-1} & A_{n} & 0 &  & \cdots 
 &  & 0 & A_{0} & \cdots & A_{n-3} & \begin{matrix}
a_{2n-4}\\
b_{2n-4}
\end{matrix}
 \\\hline
\vdots &  &  & \ddots &  &  &  &  &  & \ddots &  & \vdots \\\hline
\begin{matrix}
a_{3}\\
b_{3}
\end{matrix}
 & A_{2} & A_{3} & \cdots & A_{n} & 0 &  & \cdots 
 &  & 0 & A_{0} & \begin{matrix}
a_{2}\\
b_{2}
\end{matrix}
\end{array}
\right) 
\label{eqBigWaveletMatrix}
\end{equation}
Except for the scalar entries in the two extreme left and right columns,
all the other entries of the big combined matrix $U_{A}$ are taken
from the cyclic arrangements of the $2\times 2$ matrices of coefficients
$A_{0},A_{1},\dots,A_{n}$ in the expansion of $A\left( z\right) $.
For the case of $n=1$ this amounts to the simple $8\times 8$ wavelet%
\index{wavelet!matrix@--- matrix} matrix
\begin{equation}
\begin{array}[c]{c}
\setlength{\unitlength}{1bp}
%\begin{picture}(216,160)(-12,0)
\begin{picture}(168,130)(12,14)
\put(96,80){\makebox(0,0){$\displaystyle
\left( 
\begin{array}{c|c|c|c|c}
\begin{matrix}
a_{1}\\
b_{1}
\end{matrix}
 & A_{1} & 0 & 0 & 
\begin{matrix}
a_{0}\\
b_{0}
\end{matrix}
 \\\hline
\begin{matrix}
0\\
0
\end{matrix}
 & A_{0} & A_{1} & 0 & \begin{matrix}
0\\
0
\end{matrix}
 \\\hline
\begin{matrix}
0\\
0
\end{matrix}
 & 0 & A_{0} & A_{1} & \begin{matrix}
0\\
0
\end{matrix}
 \\\hline
\begin{matrix}
a_{3}\\
b_{3}
\end{matrix}
 & 0 & 0 & A_{0} & \begin{matrix}
a_{2}\\
b_{2}
\end{matrix}
\end{array}
\right) 
$}}
%\put(0,136){\includegraphics[bb=0 0 192 25]{cycloop1.eps}}
\put(24,126){\includegraphics[bb=0 0 192 25,width=144bp]{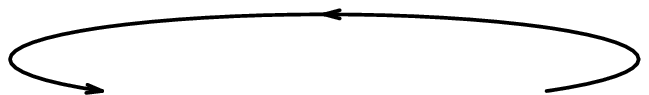}}
%\put(-2,148){\makebox(0,0)[r]{$A_{0}$}}
\put(22,134){\makebox(0,0)[r]{$A_{0}$}}
%\put(0,0){\includegraphics[bb=0 0 192 25]{cycloop2.eps}}
\put(24,14){\includegraphics[bb=0 0 192 25,width=144bp]{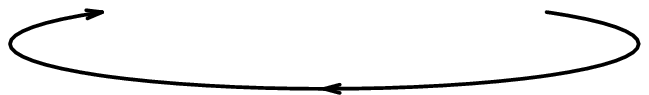}}
%\put(194,12){\makebox(0,0)[l]{$A_{1}$,}}
\put(170,26){\makebox(0,0)[l]{$A_{1}$,}}
\end{picture}
\end{array}
\label{eqWaveletMatrix}
\end{equation}
which is the one that produces the sequence of quantum%
\index{gate!quantum} gates.
%in
%Exercise \ref{Exe1Oct24-9}.
The quantum algorithm%
\index{algorithm!quantum wavelet} of a wavelet filter%
\index{filter!wavelet}
is thus represented by a $2^{n+2}\times 2^{n+2}$ unitary matrix%
\index{matrix!unitary} $U_{A}$
acting on the quantum%
\index{quantum!register@--- register} qubit%
\index{qubit!register@--- register} register
$\underset{n+2\text{ times}}
{\underbrace{\mathbb{C}\otimes\dots\otimes\mathbb{C}}}
=\mathbb{C}^{2\left( n+2\right) }$,
i.e., it acts on a configuration of $n+2$ qubit%
\index{qubit}s. The realization
of a wavelet algorithm%
\index{algorithm!wavelet}%
\index{algorithm!quantum wavelet} in the quantum%
\index{gate!quantum} realm
thus amounts to spelling out the steps
in factoring $U_{A}$ into a product of
qubit%
\index{gate!qubit} gates. By Shor's theorem%
\index{theorem!Shor}, we
know that this can be done, and $U_{A}$ may be built out of
one-qubit%
\index{gate!qubit} gates and CNOT%
\index{gate!CNOT} gates
following the ideas sketched above.
The reader may find more discussion
of the matrix $U_{A}$ in Section 3 of \cite{Fre00}.

The generalization of classical and quantum wavelet resolution algorithm%
\index{algorithm!wavelet resolution}%
\index{algorithm!quantum wavelet}s
from $N=2$ to $N>2$ is immediate: Then
$m_{i}\left(  z\right)  =\sum_{k\in\mathbb{Z}}a_{k}^{\left(  i\right)  }z^{k}%
$,
\begin{equation}
\left(  S_{i}f\right)  \left(  z\right)  =m_{i}\left(  z\right)  f\left(
z^{N}\right)  ,\qquad i=0,\dots,N-1, \label{eqIntMay7.7}%
\end{equation}
and the transformation rules%
\begin{equation}
\xi_{Nk+i}=\sum_{l\in\mathbb{Z}}a_{l-Nk}^{\left(  i\right)  }\varepsilon
_{l},\qquad i=0,1,\dots,N-1, \label{eqIntApr23.star3}%
\end{equation}
permute the set of ONB's in $\ell^{2}\left(  \mathbb{Z}\right)  $ and define a
unitary%
\index{operator!unitary} commuting with the $N$-shift%
\index{operator!shift}. Hence, the standard formulas from
\cite{Wic93}%
\index{Wickerhauser0@M.V. Wickerhauser}, \cite{Kla99}%
\index{Klappenecker0@Andreas Klappenecker}, and \cite{FiWi99} for the quantum computing algorithm%
\index{algorithm!quantum computing}
%,
%which are based on (\ref{eqIntApr23.star1}), 
naturally generalize to the case
$N>2$ via (\ref{eqIntApr23.star3}). Instead of $k$-registers $\underset
{k\text{ times}}{\underbrace{\mathbb{C}^{2}\otimes\dots\otimes\mathbb{C}^{2}}%
}=\mathbb{C}^{2^{k}}$ over $\mathbb{C}^{2}$, we will now have to work rather
with
$\underset{k\text{ times}}{\underbrace{\mathbb{C}^{N}\otimes
\dots\otimes\mathbb{C}^{N}}}
=\mathbb{C}^{N^{k}}$.

The use of the 
algorithmic
relations 
%(\ref{eqIntApr23.star2}) 
in engineering
and operator%
\index{operator!algebra@--- algebra} algebra theory predates their
more recent use in wavelet theory and
wavepacket analysis%
\index{analysis!wavepacket}.

\subsubsection{\label{S2.2.1}Pyramid algorithms}

For $N>2$, the algorithm of the previous section takes the following form.
\bigskip

\label{Exe2Oct5-21}
\textbf{The pyramid algorithm%
\index{algorithm!pyramid} and the Cuntz%
\index{Cuntz relations} relations revisited.}
%This exercise is a generalization of Exercise
%\textup{\ref{Exe1Sep7-4}} to the general setting
%of\kern0.2pt{} Theorem \textup{\ref{ThmMul.1}}, and to
%begin with we also assume \textup{(\ref{eqMul.12})}. Again, 
By
Fourier%
\index{Fourier1@Fourier!equivalence@--- equivalence} equivalence of
$L^{2}\left(  \mathbb{T}\right)  $ and $\ell
^{2}\left(  \mathbb{Z}\right)  $ via the Fourier%
\index{Fourier1@Fourier!series@--- series} series,
it follows that the system $\left\{  S_{i}\right\}
_{i=0}^{N-1}$
%defined by \textup{(\ref{eqMul.13})}
is equivalent to a system $\left\{  
\smash{\hat{S}_{i}}
\right\}  _{i=0}^{N-1}$
acting on $\ell^{2}\left(  \mathbb{Z}\right)  $. Specifically, $\left(
S_{i}f\right)  \sphat=\hat{S}_{i}\hat{f}$, $i=0,\dots,N-1$, where $\hat{f}\left(
n\right)  :=\int_{\mathbb{T}}z^{-n}f\left(  z\right) \,d\mu\left(  z\right)
$.
%, satisfy the identities \textup{(\ref{eqMul.14}%
%)}. 
For $c:=\left(  c_{n}\right)
_{n\in\mathbb{Z}}$ in $\ell^{2}\left(  \mathbb{Z}\right)  $, and functions $f$
on $\mathbb{R}$, set
\[
f_{-1}\left(  x\right)  :=N^{-\frac{1}{2}}f\left(
\frac{x}{N}\right)  ,
\]
and
\[
\left(  c\ast f\right)  \left(  x\right)
:=\sum_{n\in\mathbb{Z}}c_{n}f\left(  x-n\right)  .
\]
Let $\left\{
m_{i}\right\}  _{i=0}^{N-1}$ be
low-pass and high-pass wavelet filter%
\index{filter!wavelet}%
\index{filter!high-pass}%
\index{filter!low-pass}s,
%as in \textup{(\ref{eqMul.7})--(\ref{eqMul.8})}, 
and let
%\linebreak
$\varphi$, $\psi_{1},\dots,\psi_{N-1}$ be the corresponding
scaling%
\index{scaling!function@--- function} function, resp., wavelet%
\index{wavelet!function@--- function} functions.
%given by \textup{(\ref{eqMul.9})--(\ref{eqMul.10})}.
Now introduce the corresponding
operators $S_{i}$, 
%of
%\textup{(\ref{eqMul.13})} and
%\textup{(\ref{eqMul.14})}, 
and their cousins $\hat
{S}_{i}$. The adjoint%
\index{operator!adjoint}s $\hat{S}_{i}^{\ast}$ are
also
called \emph{filter%
\index{filter|textbf}s}.

%\begin{enumerate}
%\item \label{Exe2Oct5-21(1)}
%Show that%
Then
\begin{equation}
c\ast\varphi=\underset{\text{coarser%
\index{resolution!coarser} resolution%
\index{resolution}}}{\underbrace{\left(  \left(
\hat{S}_{0}^{\ast}c\right)  \ast\varphi\right)  _{-1}}}+\underset
{\text{detail%
\index{detail}}}{\underbrace{
\sum_{i=1}^{N-1}
\left(  \left(  \hat{S}_{i}^{\ast}c\right)
\ast\psi_{i}\right)  _{-1}}}\text{\qquad for all }c\in\ell^{2}\left(
\mathbb{Z}\right)  .\label{eqExe2Oct5-21(1).1}%
\end{equation}

%\item \label{Exe2Oct5-21(2)}
The scaling%
\index{scaling!function@--- function} function $\varphi$ defines a resolution%
\index{resolution!subspace@--- subspace} subspace
$\mathcal{V}_{0}\subset L^{2}\left(  \mathbb{R}\right)  $.
%by \textup{(\ref{eqMul.15})}. Show that
%\textup{(\ref{eqExe2Oct5-21(1).1})} represents the orthogonal%
%\index{orthogonal!decomposition@--- decomposition} decomposition of
%functions in $\mathcal{V}_{0}$
%into an orthogonal sum of a function with
%coarser%
%\index{resolution!coarser} resolution and a function in the intermediate detail%
%\index{detail!subspace@--- subspace} subspace.
%For an illustration, see
%again the
%table
%in Exercise \textup{\ref{Exe1Sep7-4}(\ref{Exe1Sep7-4(2)})},
%modified 
For the case $N>2$:

\begin{center}
$\renewcommand{\arraystretch}{0.25}\addtolength{\tabcolsep}{-3pt}%
%\mkern-27mu
\begin{tabular}
[c]{ccccrrrcl}%
\multicolumn{9}{c}{\parbox{0.9\textwidth}
{\centering\small Discrete%
\index{wavelet!discrete} vs.\ continuous%
\index{wavelet!continuous} wavelets, i.e., $\ell\sp 2$ vs.\ $L\sp 
2\left(  \mathbb{R}\right)  $:\rule[-12pt]{0pt}{12pt}}}\\
\mathstrut &  &  &  &  &  &  &  & \\
$\left\{  0\right\}  $ & $\longleftarrow$ & $\!\cdots\!$ & $\longleftarrow$ &
$\mathcal{V}_{2}\raisebox{-12pt}{$\searrow\vphantom{\raisebox{-2pt}{$\searrow
$}}$\hskip-6pt}$ & $\mathcal{V}_{1}\raisebox{-12pt}{$\searrow\vphantom
{\raisebox{-2pt}{$\searrow$}}$\hskip-6pt}$ & $\mathcal{V}_{0}\raisebox
{-12pt}{$\searrow\vphantom{\raisebox{-2pt}{$\searrow$}}$\hskip-6pt}$ &  &
finer%
\index{scaling!finer} scales\\\hline
&  &  &  & \multicolumn{1}{c}{} & \multicolumn{1}{c}{} & \multicolumn{1}{c}{}
&  & \\\cline{1-7}
&  &  &  & \multicolumn{1}{c}{} & \multicolumn{1}{c}{} & \multicolumn{1}{c}{}
& \multicolumn{1}{|c}{} & \\\cline{1-6}
&  &  &  & \multicolumn{1}{c}{} & \multicolumn{1}{c}{} & \multicolumn{1}{|c}{}
& \multicolumn{1}{|c}{} & \\\cline{1-5}
&  &  &  & \multicolumn{1}{c}{} & \multicolumn{1}{|c}{} &
\multicolumn{1}{|c}{} & \multicolumn{1}{|c}{} & \\\cline{1-4}
&  & $\!\cdots\!$ &  & \multicolumn{1}{|c}{$\mathcal{W}_{3}$} &
\multicolumn{1}{|c}{$\mathcal{W}_{2}$} & \multicolumn{1}{|c}{$\mathcal{W}_{1}%
$} & \multicolumn{1}{|c}{$\cdots\!\!$} & rest of $L^{2}\left(  \mathbb{R}\right)
$\\\cline{1-4}
&  &  &  & \multicolumn{1}{c}{} & \multicolumn{1}{|c}{} &
\multicolumn{1}{|c}{} & \multicolumn{1}{|c}{} & \\\cline{1-5}
&  &  &  & \multicolumn{1}{c}{} & \multicolumn{1}{c}{} & \multicolumn{1}{|c}{}
& \multicolumn{1}{|c}{} & \\\cline{1-6}
&  &  &  & \multicolumn{1}{c}{} & \multicolumn{1}{c}{} & \multicolumn{1}{c}{}
& \multicolumn{1}{|c}{} & \\\cline{1-7}
&  &  &  & \multicolumn{1}{c}{} & \multicolumn{1}{c}{} & \multicolumn{1}{c}{}
&  & \\\hline
\rule[6pt]{0pt}{6pt} &  & $\!\cdots\!$ & \multicolumn{1}{r}{$\rlap{$%
\underset{\textstyle U}{\longleftarrow}$}\hskip3pt$} & $\rlap{$\underset
{\textstyle U}{\longleftarrow}$}\hskip3pt$ & $\rlap{$\underset{\textstyle
U}{\longleftarrow}$}\hskip3pt$ & \multicolumn{1}{c}{} &  & \\
&  & $\llap{$W$}\makebox[6pt]{\raisebox{4pt}{\makebox[0pt]{\hss$\uparrow$\hss
}}\raisebox{-4pt}{\makebox[0pt]{\hss$|$\hss}}}$ &  &
\multicolumn{1}{c}{$\makebox[6pt]{\raisebox{4pt}{\makebox[0pt]{\hss$%
\uparrow$\hss}}\raisebox{-4pt}{\makebox[0pt]{\hss$|$\hss}}}$} &
\multicolumn{1}{c}{$\makebox[6pt]{\raisebox{4pt}{\makebox[0pt]{\hss$%
\uparrow$\hss}}\raisebox{-4pt}{\makebox[0pt]{\hss$|$\hss}}}$} &
\multicolumn{1}{c}{$\makebox[6pt]{\raisebox{4pt}{\makebox[0pt]{\hss$%
\uparrow$\hss}}\raisebox{-4pt}{\makebox[0pt]{\hss$|$\hss}}}\rlap{$W$}$} &  &
\\
\rule[-6pt]{0pt}{6pt}$\left\{  0\right\}  $ & $\longleftarrow$ & $\!\cdots\!$ &
\multicolumn{1}{r}{$\rlap{$\overset{\textstyle S_0}{\longleftarrow}$}%
\hskip3pt$} & $\rlap{$\overset{\textstyle S_0}{\longleftarrow}$}\hskip3pt$ &
$\rlap{$\overset{\textstyle S_0}{\longleftarrow}$}\hskip3pt$ &
\multicolumn{1}{c}{} &  & \\\cline{1-7}
&  &  &  & \multicolumn{1}{c}{} & \multicolumn{1}{c}{} & \multicolumn{1}{c}{}
& \multicolumn{1}{|c}{} & \\\cline{1-6}
&  &  &  & \multicolumn{1}{c}{} & \multicolumn{1}{c}{} & \multicolumn{1}{|c}{}
& \multicolumn{1}{|c}{} & \\\cline{1-5}
&  &  &  & \multicolumn{1}{c}{} & \multicolumn{1}{|c}{} &
\multicolumn{1}{|c}{} & \multicolumn{1}{|c}{} & \\\cline{1-4}
&  & $\!\cdots\!$ &  & \multicolumn{1}{|c}{$S_{0}^{2}\mathcal{L}$} &
\multicolumn{1}{|c}{$S_{0}\mathcal{L}$} & \multicolumn{1}{|c}{$
\,\mathcal{L}
=\smash{\bigvee\limits_{i=1}^{N-1}}
S_{i}\ell^{2}$\,} & \multicolumn{1}{|c}{} & \\\cline{1-4}
&  &  &  & \multicolumn{1}{c}{} & \multicolumn{1}{|c}{} &
\multicolumn{1}{|c}{} & \multicolumn{1}{|c}{} & \\\cline{1-5}
&  &  &  & \multicolumn{1}{c}{} & \multicolumn{1}{c}{} & \multicolumn{1}{|c}{}
& \multicolumn{1}{|c}{} & \\\cline{1-6}
&  &  &  & \multicolumn{1}{c}{} & \multicolumn{1}{c}{} & \multicolumn{1}{c}{}
& \multicolumn{1}{|c}{} & \\\cline{1-7}
&  &  &  & \makebox[36pt]{\hfill}\llap{$S_{0}^{2}\ell^{2}$}\raisebox
{12pt}{$\nearrow\vphantom{\raisebox{2pt}{$\nearrow$}}$\hskip-6pt} &
\makebox[36pt]{\hfill}\llap{$S_{0}\ell^{2}$}\raisebox{12pt}{$\nearrow
\vphantom{\raisebox{2pt}{$\nearrow$}}$\hskip-6pt} & $\ell^{2}\raisebox{12pt}%
{$\nearrow\vphantom{\raisebox{2pt}{$\nearrow$}}$\hskip-6pt}$ &  &
\end{tabular}
\ $
\end{center}

%\emph{Hint:} Use the
%method in Exercise \textup{\ref{Exe1Sep7-4}},
%the proof of\kern0.2pt{} Theorem \textup{\ref{ThmMul.1}},
%and Remark \textup{\ref{RemCoarse}}.

%\item \label{Exe2Oct5-21(3)}
%%\exeonestar
%Discuss to what extent \textup{(\ref{Exe2Oct5-21(1)})} and
%\textup{(\ref{Exe2Oct5-21(2)})} are
%still true if one does not assume \textup{(\ref{eqMul.12})}, i.e., if
%one merely is in the tight frame%
%\index{frame!tight} context; i.e., end
%the discussion started in Remark \textup{\ref{RemCoarse}}.
%\end{enumerate}

More refined pyramid algorithms yield wavelet packets as follows.

The Haar%
\index{Haar1@Haar!wavelet@--- wavelet} wavelet
is supported in $\left\lbrack 0,1\right\rbrack $,
and if $j\in\mathbb{Z}_{+}$ and $k\in\mathbb{Z}$, then the modified function
$x\mapsto\psi\left( 2^{j}x-k\right) $ is supported in the smaller interval
$\frac{k}{2^{j}}\leq x\leq\frac{k+1}{2^{j}}$. When $j$ is fixed, these intervals
are contained in $\left\lbrack 0,1\right\rbrack $ for
$k\in\left\{ 0,1,\dots,2^{j}-1\right\} $.
This
is not the case
for the other wavelet%
\index{wavelet!function@--- function}
functions.
For one thing, the non-Haar wavelets $\psi$ have support intervals
of length more than one,
and this
forces periodicity considerations;
see \cite{CDV93}%
\index{Daubechies0@Ingrid Daubechies}.
For this reason, Coifman%
\index{Coifman0@R.R. Coifman} and
Wickerhauser%
\index{Wickerhauser0@M.V. Wickerhauser} \cite{CoWi93}%
\index{Coifman0@R.R. Coifman}%
\index{Wickerhauser0@M.V. Wickerhauser} invented the
concept of wavelet packet%
\index{wavelet!packet@--- packet}s.
They are built from functions
with prescribed smoothness,
and yet they have localization
properties that rival
those of the (discontinuous)
Haar%
\index{Haar1@Haar!wavelet@--- wavelet} wavelet.

There
are powerful but nontrivial theorem%
\index{theorem!wavelet restriction}s on restriction
algorithm%
\index{algorithm!restriction}s for wavelets
$\psi_{j,k}\left( x\right) =2^{\frac{j}{2}}\psi\left( 2^{j}x-k\right) $
from $L^{2}\left( \mathbb{R}\right) $ to $L^{2}\left( 0,1\right) $. We
refer the reader to \cite{CDV93}%
\index{Daubechies0@Ingrid Daubechies} and \cite{MiXu94} for
the details of this construction.
The
underlying idea of 
%Exercises \ref{Exe1Sep28-8} and \ref{Exe1Sep28-9}
%dates back to 
Alfred Haar%
\index{Haar0@Alfred Haar}
%, but
%it 
has found a recent renaissance
in the work of Wickerhauser%
\index{Wickerhauser0@M.V. Wickerhauser} \cite{Wic93}%
\index{Wickerhauser0@M.V. Wickerhauser}
on \emph{wavelet packet%
\index{wavelet!packet@--- packet}s}. The idea there,
which is also motivated by the
Walsh function algorithm%
\index{algorithm!Walsh function}, is to replace
the refinement%
\index{refinement!equation@--- equation} equation (\ref{eqInt.1b}) by a
related recursive system as follows:
Let $m_{0}\left( z\right) =\sum_{k}a_{k}z^{k}$,
$m_{1}\left( z\right) =\sum_{k}b_{k}z^{k}$,
for example $b_{k}=\left( -1\right) ^{k}\bar{a}_{1-k}$,
$k\in\mathbb{Z}$,
be a given low-pass/high-pass%
\index{filter!high-pass}%
\index{filter!low-pass} system, $N=2$. Then
consider the following \emph{refinement system} on $\mathbb{R}$:
\begin{align}
W_{2n}\left( x\right) 
&=\sqrt{2}\sum_{k\in\mathbb{Z}}a_{k}W_{n}\left( 2x-k\right) ,\label{eqWaveletPacketRefinement(1)}
\\
W_{2n+1}\left( x\right) 
&=\sqrt{2}\sum_{k\in\mathbb{Z}}b_{k}W_{n}\left( 2x-k\right) .
\label{eqWaveletPacketRefinement}
\end{align}
Clearly the function $W_{0}$ can be identified
with the traditional scaling%
\index{scaling!function@--- function} function $\varphi$
of (\ref{eqInt.1}). 
A theorem%
\index{theorem!Coifman Wickerhauser@Coifman--Wickerhauser} of Coifman%
\index{Coifman0@R.R. Coifman}
and Wickerhauser%
\index{Wickerhauser0@M.V. Wickerhauser} (Theorem 8.1, \cite{CoWi93}%
\index{Coifman0@R.R. Coifman}%
\index{Wickerhauser0@M.V. Wickerhauser})
states that if $\mathcal{P}$ is a partition of
$\left\{ 0,1,2,\dots\right\} $ into subsets of the form
\[
I_{k,n}=\left\{ 2^{k}n,2^{k}n+1,\dots ,2^{k}\left( n+1\right) -1\right\} ,
\]
then the function system
\[
\left\{ 2^{\frac{k}{2}}W_{n}\left( 2^{k}x-l\right) 
\Bigm| I_{k,n}\in\mathcal{P},\;l\in\mathbb{Z}\right\} 
\]
is an orthonormal%
\index{orthonormal!basis@--- basis} basis for $L^{2}\left( \mathbb{R}\right) $.
Although it is not spelled out
in \cite{CoWi93}%
\index{Coifman0@R.R. Coifman}%
\index{Wickerhauser0@M.V. Wickerhauser}, this construction
of bases in $L^{2}\left( \mathbb{R}\right) $ divides
itself into the two cases, the
true orthonormal%
\index{orthonormal!basis@--- basis} basis (ONB), and the
weaker property of forming a function
system which is only a tight
frame%
\index{frame!tight}. As in the wavelet case, to get the $\mathcal{P}$-system to
really be an ONB for $L^{2}\left( \mathbb{R}\right) $,
we must assume the transfer
operator%
\index{operator!transfer} $R_{\left| m_{0}\right| ^{2}}$ to have
\emph{%
Per\-ron--Fro\-be\-ni\-us
spectrum%
\index{Perron Frobenius@Perron--Frobenius!spectrum@--- spectrum}} on $C\left( \mathbb{T}\right) $.
This means that the intersection of the
point spectrum%
\index{spectrum!peripheral point} of $R_{\left| m_{0}\right| ^{2}}$ with $\mathbb{T}$ is
the singleton $\lambda =1$, and that
$\dim\ker ( ( 1-R_{\left| m_{0}\right| ^{2}}) 
|_{C( \mathbb{T}) } ) =1$.
%;
%see Table \ref{TabDauSum} in the next
%chapter.

\subsubsection{\label{S2.2.2}Subdivision algorithms}

The algorithms for wavelets and wavelet packets involve the pyramid idea as
well as subdivision. Each subdivision produces a multiplication of subdivision
points. If the scaling is by $N$, then $j$ subdivisions multiply the number of
subdivision points by $N^{j}$. If the scaling is by a $d\times d$ integral
matrix $\mathbf{N}$, then the multiplicative factor is $\left\vert
\det\mathbf{N}\right\vert ^{j}$ in the number of subdivision points placed in
$\mathbb{R}^{d}$.

In the discussion below, we restrict attention to $d=1$, but the conclusions
hold with only minor modification in the general case of $d>1$ and matrix scaling.

%We have already
%noted that if 
If
$W$ is a continuous function on $\mathbb{T}$, the \emph{transfer
operator%
\index{operator!transfer}} or \emph{kneading operator%
\index{operator!kneading}} $R_{W}$ 
%defined by (\ref{eqSpe.1}),%
\begin{equation}
R_{W}\xi\left(  z\right)  =\frac{1}{N}\sum_{w^{N}=z}W\left(  w\right)
\xi\left(  w\right)  =S_{0}^{\ast}W\xi\left(  z\right)  , \label{eqk&c.1bis}%
\end{equation}
with the alias 
%given by (\ref{eqSpe.12}),%
\begin{equation}
\left(  R_{W}f\right)  _{n}=\sum_{k}c_{Nn-k}f_{k} \label{eqk&c.2}%
\end{equation}
in the Fourier%
\index{Fourier1@Fourier!transform@--- transform} transformed space, has an adjoint%
\index{operator!adjoint} which is the \emph{subdivision
operator%
\index{operator!subdivision}} or chopping operator%
\index{operator!chopping} 
%given by (\ref{eqIsoMar5.6}),%
\begin{equation}
\left(  R_{W}^{\ast}\xi\right)  \left(  z\right)  =\overline{W\left(
z\right)  }\,\xi\left(  z^{N}\right)  \label{eqk&c.3}%
\end{equation}
on functions $\xi$ on $\mathbb{T}$, with the alias 
%given by
%(\ref{eqSpbMay10.5}),
\begin{equation}
\left(  R_{W}^{\ast}f\right)  _{n}=\sum_{k}\overline{c_{Nk-n}}\,f_{k}
\label{eqk&c.4}%
\end{equation}
on sequences.

We will analyze the duality between $R_{W}$ and $R_{W}^{\ast}$ and their
spectra. 
%in Sections \ref{k&c} and \ref{Gre}.
%
Specializing to $W=\left|  m_{0}\right|  ^{2}$, we note that
$R_{W}$ is then the transfer operator of
orthogonal type%
\index{orthogonal!type@--- type} wavelet%
\index{orthogonal!wavelet@--- wavelet}s.
%is (\ref{eqIntNew.4}).
%or (\ref{eqSpe.8}),
%while biorthogonal type%
%\index{biorthogonal!type@--- type}.
% (referring to the wider class of
%biorthogonal%
%\index{wavelet!biorthogonal} wavelets) is (\ref{eqIntNew.7}) or (\ref{eqRTO.3}) 
%in Proposition \ref{ProRTO.1} or (\ref{eqRTO.21}) in Theorem \ref{ProRTO.4}.
In the following, 
%we shall refer 
%to these two conditions in
%connection with 
%a given 
$W$ 
%which 
is assumed only to satisfy 
%(\ref{eqRTO.19}),
%i.e., 
$W\in\operatorname*{Lip}\nolimits%
\index{Lipschitz!space@--- space}_{1}\left(  \mathbb{T}\right)  $ and
%(\ref{eqRTO.21}), i.e., 
$W\geq0$. Other conditions 
%will then be added.
%, such as
%(\ref{eqSpe.8}), or (\ref{eqRTO.21}).
are discussed in \cite{BrJo02b}.

In the engineering terminology of 
%Exercise 
\S\ \ref{Exe1Sep7-4}, the operation
(\ref{eqk&c.2}) is composed of a local filter%
\index{filter} with the numbers $c_{k}$ as coefficients, followed
by the down-sampling 
$\vcenter{\hbox{\Large $\bigcirc \llap{$\vcenter{\hbox{\normalsize $\scriptscriptstyle N$}}\mkern-1mu\vcenter{\hbox{\normalsize $\scriptscriptstyle\downarrow$}}\mkern3.5mu$}$}}$,
%$N\downarrow$, 
while
(\ref{eqk&c.4}) is composed of up-sampling 
$\vcenter{\hbox{\Large $\bigcirc \llap{$\vcenter{\hbox{\normalsize $\scriptscriptstyle N$}}\mkern-1mu\vcenter{\hbox{\normalsize $\scriptscriptstyle\uparrow$}}\mkern3.5mu$}$}}$,
%$N\uparrow$,
followed by an application of a dual filter%
\index{filter!dual}. In signal processing%
\index{signal!processing@--- processing}, 
$\vcenter{\hbox{\Large $\bigcirc \llap{$\vcenter{\hbox{\normalsize $\scriptscriptstyle N$}}\mkern-1mu\vcenter{\hbox{\normalsize $\scriptscriptstyle\downarrow$}}\mkern3.5mu$}$}}$
%$N\downarrow$
is referred to as ``decimation'' even if $N$ is not $10$.

The operator $S$ ($=R_{W}^{\ast}$) 
%in (\ref{eqSpbMay10.5}) 
is called the
subdivision operator%
\index{operator!subdivision}, or the \emph{woodcutter operator%
\index{operator!woodcutter}}, because of its use in
computer graphics. Iterations of $S$ will generate a shape which (in the case
of one real dimension) takes the form of the graph of a function $f$ on
$\mathbb{R}$. If $\xi\in\ell^{\infty}\left(  \mathbb{Z}\right)  $ is given,
and if the differences%
\begin{equation}
D_{n}\left(  i\right)  =f\left(  \frac{i}{2^{n}}\right)  -\left(  S^{n}%
\xi\right)  \left(  i\right)  ,\qquad i\in\mathbb{Z},
\label{eqSpbMay14.poundbis}%
\end{equation}
are small, for example if
\begin{equation}
\lim_{n\rightarrow\infty}\sup_{i\in\mathbb{Z}}\left|  D_{n}\left(  i\right)
\right|  =0, \label{eqSpbMay14.small}%
\end{equation}
then we say that $\xi$ represents \emph{control points}, or a control polygon,
and the function $f$ is the limit of the \emph{subdivision%
\index{subdivision} scheme}.

It follows that the subdivision operator%
\index{operator!subdivision} $S$ on the sequence spaces,
especially on $\ell^{\infty}\left(  \mathbb{Z}\right)  $, governs
\emph{pointwise approximation%
\index{approximation!pointwise}} to refinable limit functions. 
%But we will see
%in Theorem \ref{CorCas.1} that 
The dual version of $S$, i.e., $R=S^{\ast}$ (=
the transfer operator%
\index{operator!transfer}) governs the corresponding \emph{mean approximation%
\index{approximation!mean@--- in mean}}
problem, i.e., approximation%
\index{approximation!L two@$L\sp {2}$-} relative to the $L^{2}\left(  \mathbb{R}\right)
$-norm.

In Scholium 4.1.2 of \cite{BrJo02b}, we 
%will 
consider the eigenvalue problem%
\begin{equation}
S\xi=\lambda\xi,\qquad\lambda\in\mathbb{C}, \label{eqSpbMay14.poundpound}%
\end{equation}
and $\xi\neq0$ in some suitably defined space of sequences. The formula
(\ref{eqSpbMay14.poundbis}) for the limit of a given subdivision%
\index{subdivision} scheme $S$
makes it clear that the case (\ref{eqSpbMay14.poundpound}) must be excluded.
For if (\ref{eqSpbMay14.poundpound}) holds, for some $\lambda\in\mathbb{C}$,
and some sequence $\xi$ of control points, then there is not a corresponding
regular function $f$ on $\mathbb{R}$ with its values given on the finer grids
$2^{-n}\mathbb{Z}$, $n=1,2,\dots$, by%
\begin{equation}
f_{\xi}\left(  i2^{-n}\right)  \approx\left(  S^{n}\xi\right)  \left(
i\right)  =\lambda^{n}\xi\left(  i\right)  ,\qquad i\in\mathbb{Z}.
\label{eqSpbMay14.poundpoundpound}%
\end{equation}
We 
%will 
show in Example 4.1.3 of \cite{BrJo02b} that there are no such control
points $\xi$ in $\ell^{2}\left(  \mathbb{Z}\right)  \setminus\left\{
0\right\}  $. Hence the stability of the algorithm!

\subsubsection{\label{S2.2.3}Wavelet packet algorithms}

The main difference between the algorithms of wavelets and those of wavelet
packets is that for the wavelets the path in the pyramid is to one side only:
a given resolution is split into a coarser one and the intermediate detail.
The intermediate detail may further be broken down into frequency bands. With
the operators $S_{j}f\left(  z\right)  =m_{j}\left(  z\right)  f\left(
z^{N}\right)  $ acting on $L^{2}\left(  \mathbb{T}\right)  $, the coarser
subspace after $j$ steps is modelled on $S_{0}^{j}L^{2}\left(  \mathbb{T}%
\right)  $, and the projection onto this subspace is $S_{0}^{j}S_{0}^{\ast
\,j}$ where $S_{0}$ is the isometry of $L^{2}\left(  \mathbb{T}\right)
\cong\mathcal{V}_{0}$ defined by the low-pass filter $m_{0}$. But in the
construction of the wavelet packet, the subspace resulting by running the
algorithm $j$ times is $S_{i_{1}}S_{i_{2}}\cdots S_{i_{j}}L^{2}\left(
\mathbb{T}\right)  $, and the projection onto this subspace is
\[
S_{i_{1}}S_{i_{2}}\cdots S_{i_{j}}S_{i_{j}}^{\ast}\cdots S_{i_{2}}^{\ast
}S_{i_{1}}^{\ast}.
\]
If $n\in\mathbb{Z}_{+}$, the wavelet function $W_{n}$ is computed from the iteration
$i_{1},\dots,i_{j}$ corresponding to the representation%
\[
n=i_{1}+i_{2}N+i_{3}N^{2}+\dots+i_{j}N^{j-1},
\]
where $i_{1},\dots,i_{j}\in\left\{  0,1,\dots,N-1\right\}  $ are unique from
the Euclidean algorithm.

\subsubsection{\label{S2.2.4}Lifting algorithms: Sweldens and more}

%As the next exercise shows,
%while 
The discussion 
%in the present
%chapter 
centers around the matrix%
\index{matrix!function@--- function}
functions
$A\colon\mathbb{T}\rightarrow\mathrm{GL}_{2}\left( \mathbb{C}\right) $.
%things
%are nicer if $\det A\left( z\right) \equiv 1$.
%The upper/lower factorization
%\textup{(\ref{eqUL})} in this exercise does
%\emph{not} apply to, say, \textup{(\ref{eqIsoFeb12.v})}
%in Remark \textup{\ref{CorIsoFeb12.1star}(a)}.

%\begin{exercise}
\label{Exe6DauSwe-11}
\textbf{The case $\det A\equiv 1$.}
Recall that we call a
finite sum
$\sum_{k=-n_{0}}^{n_{1}}A_{k}z^{k}$, $n_{0},n_{1}\geq 0,$
a Fourier%
\index{Fourier1@Fourier!polynomial@--- polynomial} polynomial
%\textup{(}see Definition \textup{\ref{DefRTOMar9.8})}
both if the coefficients $A_{k}$ are numbers, and if they are matrices.
The matrix%
\index{matrix!function@--- function}-valued Fourier%
\index{Fourier1@Fourier!polynomial@--- polynomial}
polynomials
$\mathbb{T}\ni z\mapsto A\left( z\right) \in M_{2}\left( \mathbb{C}\right) $
such that
$\det A\left( z\right) \equiv 1$ form a subgroup of
$C\left( \mathbb{T},\mathrm{GL}_{2}\left( \mathbb{C}\right) \right) $
which we denote $\mathcal{SL}_{2}$.

%\begin{enumerate}
%\item \label{Exe6DauSwe-11-a}
%\exetwostar
%Show that 
For every
$A\left( z\right) $ in $\mathcal{SL}_{2}$ there are
$m\in\mathbb{Z}_{+}$, $K\in\mathbb{C}\setminus\left\{ 0\right\} $,
and scalar-valued
Fourier%
\index{Fourier1@Fourier!polynomial@--- polynomial} polynomials
$u_{1}(z),\dots,u_{m}(z),l_{1}(z),\dots,l_{m}(z)$ such that
\begin{multline}
A\left( z\right) 
=\left( 
\begin{matrix}
K & 0 \\
0 & K^{-1}
\end{matrix}
\right) 
\cdot
\left( 
\begin{matrix}
1 & 0 \\
l_{1}\left( z\right)  & 1
\end{matrix}
\right) 
\cdot
\left( 
\begin{matrix}
1 & u_{1}\left( z\right)  \\
0 & 1
\end{matrix}
\right) 
\cdot \\
\cdot\left( 
\begin{matrix}
1 & 0 \\
l_{2}\left( z\right)  & 1
\end{matrix}
\right) 
\cdot
\left( 
\begin{matrix}
1 & u_{2}\left( z\right)  \\
0 & 1
\end{matrix}
\right) 
\cdots
\left( 
\begin{matrix}
1 & 0 \\
l_{m}\left( z\right)  & 1
\end{matrix}
\right) 
\cdot
\left( 
\begin{matrix}
1 & u_{m}\left( z\right)  \\
0 & 1
\end{matrix}
\right) .
\label{eqUL}
\end{multline}
%
%\item \label{Exe6DauSwe-11-b}
%Spell out an algorithm%
%\index{algorithm} for the factorization.
%
%\item \label{Exe6DauSwe-11-c}
%Ponder the extent of nonuniqueness in the factorization.
%\end{enumerate}
%
%\noindent\emph{Hint:} 
See \textup{\cite{DaSw98}%
\index{Sweldens0@W. Sweldens}%
\index{Daubechies0@Ingrid Daubechies}}.
This is the first step in the
Daubechies%
\index{Daubechies1@Daubechies!Sweldens lifting@--- --Sweldens lifting}--Sweldens
lifting algorithm%
\index{algorithm!Daubechies Sweldens lifting@Daubechies--Sweldens lifting} for the
discrete%
\index{wavelet!discrete} wavelet%
\index{wavelet!transform@--- transform} transform.
Thus the case $\det\left( A\left( z\right) \right) =1$ gives a constructive
lifting algorithm%
\index{algorithm} for wavelets, and such an
algorithm%
\index{algorithm} has not been established in the
$C\left( \mathbb{T},\mathrm{GL}_{2}\left( \mathbb{C}\right) \right) $
case. The decomposition%
\index{decomposition} could
also be compared with \textup{Proposition 3.3 of \cite{BrJo02a}},
which was mentioned in connection with the proof of \textup{(\ref{eqHom.3})}.
%\end{exercise}

%\begin{exercise}
\label{Exe6Dec3-last}
Recall the correspondence between
matrix%
\index{matrix!function@--- function} functions and wavelet
filter%
\index{filter!wavelet}s: If
%\linebreak
$A\colon\mathbb{T}\rightarrow\mathrm{GL}_{2}\left( \mathbb{C}\right) $ is
a matrix%
\index{matrix!function@--- function} function, then the corresponding dyadic
wavelet%
\index{wavelet!dyadic} filter%
\index{filter!wavelet}s are
\[
m_{i}^{\left( A\right) }\left( z\right) 
=\sum_{j=0}^{1}A_{i,j}\left( z^{2}\right) z^{j},\qquad i=0,1.
\]
%
%\begin{enumerate}
%\item
%\label{Exe6Dec3-last-a}
%Show 
It follows
that the two matrix%
\index{matrix!function@--- function} functions $A$ and $B$
satisfy
\[
A=\left( 
\begin{matrix}
1 & 0 \\
l & 1
\end{matrix}
\right) B
\]
for some $l$ in the ring $\mathcal{F}$ of Fourier%
\index{Fourier1@Fourier!polynomial@--- polynomial}
polynomials if and only if
$m_{0}^{\left( A\right) }=m_{0}^{\left( B\right) }$ and
$m_{1}^{\left( A\right) }\left( z\right) 
=m_{1}^{\left( B\right) }\left( z\right) 
+l\left( z^{2}\right) m_{0}^{\left( A\right) }\left( z\right) $.

%\item
%\label{Exe6Dec3-last-b}
%Show 
Similarly note that the two matrix%
\index{matrix!function@--- function} functions $A$ and $B$
satisfy
\[
A=\left( 
\begin{matrix}
1 & u \\
0 & 1
\end{matrix}
\right) B
\]
for some $u\in\mathcal{F}$
if and only if
$m_{1}^{\left( A\right) }=m_{1}^{\left( B\right) }$ and
$m_{0}^{\left( A\right) }\left( z\right) 
=m_{0}^{\left( B\right) }\left( z\right) 
+u\left( z^{2}\right) m_{1}^{\left( A\right) }\left( z\right) $.
%\end{enumerate}
%\end{exercise}

\noindent\textbf{Remark.}
The conclusion 
%from
%Exercises \ref{Exe6Nov19-12} and \ref{Exe6Dec3-last} 
is that
the wavelet algorithm%
\index{algorithm!wavelet}
for a general
wavelet filter%
\index{filter!wavelet} corresponding to a matrix%
\index{matrix!function@--- function}
function, say $A$, may be broken down
in a sequence of zig-zag steps acting alternately
on the high-pass and the low-pass signal
components%
\index{low-pass!signal component@--- signal component}%
\index{high-pass!signal component@--- signal component}.

\subsection{\label{S2.2bis}Factorization theorems for matrix functions}

We mentioned that for matrix functions corresponding to finite impulse
response (FIR) filters which are unitary, we need only the constant matrix
(which is chosen such as to achieve the high-pass and low-pass conditions) and
factors of the form%
\[
U_{P}\left(  z\right)  =zP+P^{\perp}\cong\left(
\addtolength{\arraycolsep}{2pt}
\begin{array}
[c]{c|c}%
z & 0\\\hline
0 & 1
\end{array}
_{\mathstrut}^{\mathstrut}
\right)
\]
where $P$ is a rank-one projection in $\mathbb{C}^{N}$ and $N$ is the scaling
number of the subdivision.

Unfortunately, no such factorization theorem is available for the non-unitary
FIR filters. But the matrix functions take values in the non-singular complex
$N\times N$ matrices. The Sweldens--Daubechies factorization and the lifting
algorithm serve as a substitute. There are still the general non-unimodular
FIR-matrix functions where factorizations are so far a bit of a mystery. The
matrix functions are called \emph{polyphase matrices} in the engineering
literature. The following summary serves as a classification theorem for the
orthogonal wavelets of compact support: the wavelets correspond to FIR
polyphase matrices which are unitary.

In summary, an algorithm%
\index{algorithm} to construct all the wavelet%
\index{wavelet!function@--- function} functions $\psi$ of
scale%
\index{scaling!number@--- number} $2$ with support in $\left[  0,2k+1\right]  $ can be established as follows:

\begin{enumerate}
\renewcommand{\theenumi}{\arabic{enumi}}
\renewcommand{\labelenumi}{$\lbrack\!\lbrack\theenumi\rbrack\!\rbrack$}
\item \label{GenAlg(1)}Pick $k$ one-dimensional orthogonal projection%
\index{operator!projection}s
$Q_{1},\dots,Q_{k}$ in $M_{2}\left(  \mathbb{C}\right)  $ and define the
unitary-valued matrix%
\index{matrix!unitary}%
\index{matrix!function@--- function} function $A\left(  z\right)  $ on $\mathbb{T}$ by%
\begin{equation}
A\left(  z\right)  =V\left(  1-Q_{1}+zQ_{1}\right)  \left(  1-Q_{2}%
+zQ_{2}\right)  \cdots\left(  1-Q_{k}+zQ_{k}\right)  ,\label{eqGen.13}%
\end{equation}
where%
\begin{equation}
V=
%\smash{
\frac{1}{\sqrt{2}}\left( 
\begin{array}
[c]{cc}%
1 & 1\\
1 & -1
\end{array}
\right)  
%}
.\label{eqGen.14}%
\end{equation}
Then each $Q_{j}$ has the form%
\begin{equation}
Q_{j}=\left(
\begin{array}
[c]{cc}%
\lambda_{j} & \sqrt{\lambda_{j}\left(  1-\lambda_{j}\right)  }e^{i\theta_{j}%
}\\
\sqrt{\lambda_{j}\left(  1-\lambda_{j}\right)  }e^{-i\theta_{j}} &
1-\lambda_{j}%
\end{array}
\right)  ,\label{eqGen.15}%
\end{equation}
where $\lambda_{j}\in\left[  0,1\right]  $ and $\theta_{j}\in\left[
0,2\pi\right)  $. (See Proposition 3.3 of \cite{BrJo02a}.)

\item \label{GenAlg(2)}Define the filter%
\index{filter}s $m_{0}\left(  z\right)  $ and
$m_{1}\left(  z\right)  $ by 
\begin{equation}
m_{i}\left(  z\right)  =\sum_{j=0}^{N-1}z^{j}A_{ij}\left(  z^{N}\right)
,\qquad i,j=0,\dots,N-1, \label{eqHom.17}%
\end{equation}
with $N=2$.

\item \label{GenAlg(3)}Define $\hat{\varphi}$ by 
\begin{equation}
\hat{\varphi}\left(  t\right)  =\prod_{k=1}^{\infty}\left(  \frac{m_{0}\left(
tN^{-k}\right)  }{\sqrt{N}}\right) . \label{eqMul.9}%
\end{equation}
If the
condition 
\begin{equation}
\operatorname*{PER}\left(  \left\vert \hat{\varphi}\right\vert ^{2}\right)
\left(  t\right)  :=\sum_{n\in\mathbb{Z}}\left\vert \hat{\varphi}\left(
t+2\pi n\right)  \right\vert ^{2}=1 \label{eqMul.12}%
\end{equation}
fails, then the algorithm%
\index{algorithm} stops.

\item \label{GenAlg(4)}If the condition (\ref{eqMul.12}) holds, one may
alternatively define $\varphi$ by the cascade algorithm%
\index{algorithm!cascade}
%in Theorem
%\ref{CorCas.1}, 
%i.e., 
\begin{align}
\varphi\left(  x\right)  &=\sqrt{N}\sum_{n\in\mathbb{Z}}a_{n}\varphi\left(
Nx-n\right)  , \label{eqInt.1}%
\\
\chi\left(  x\right)  &=%
\left\{ \begin{array} [c]{ll}1, & 0\leq x<1,\\ 0, & x\in\mathbb
{R}\setminus\left[ 0,1\right\rangle, \end{array} \right. 
\label{eqIntApr25.4}%
\\
M_{a}&\colon\psi\longmapsto\sqrt{N}\sum_{n}a_{n}\psi\left(  Nx-n\right)  .
\label{eqIntMar28.4}%
\end{align}

\item \label{GenAlg(5)}The wavelet%
\index{wavelet!function@--- function} function $\psi$ is then defined by
\begin{equation}
\psi_{i}\left(  x\right)  =\sqrt{N}\sum_{n\in\mathbb{Z}}a_{n}^{\left(
i\right)  }\varphi\left(  Nx-n\right)  , \label{eqHom.21}%
\end{equation}
where $a_{n}^{\left(  i\right)  }$ are the Fourier%
\index{Fourier1@Fourier!coefficients@--- coefficient(s)} coefficients of $m_{i}$,%
\begin{equation}
m_{i}\left(  z\right)  =\sum_{n}a_{n}^{\left(  i\right)  }z^{n},
\label{eqHom.22}%
\end{equation}
and $z=e^{-it}$;
this is the most general
wavelet%
\index{wavelet!function@--- function} function with support in $\left[  0,2k+1\right]  $.

\item \label{GenAlg(6)}All other wavelet%
\index{wavelet!compactly supported}%
\index{wavelet!function@--- function} functions with compact support can be
obtained from the ones in
$\lbrack\!\lbrack\ref{GenAlg(5)}\rbrack\!\rbrack$
by integer translation. 
%See
%Exercise \ref{Exe2Sep6-20}.
\end{enumerate}

\subsubsection{\label{S2.2bis.1}The case of polynomial functions [the polyphase matrix, joint work with 
Ola Bratteli]}

\label{Bio}
One problem occurring in the biorthogonal context which does not have an
analogue in the orthogonal setting stems from the fact that the
duality
relations
\begin{equation}
\sum_{w^{N}=z}\overline{m_{i}\left(  w\right)  }\,\tilde{m}_{j}\left(
w\right)  =N\delta_{i,j} \text{\qquad for }i,j=0,\dots,N-1
\label{eqBio.1}%
\end{equation}
do not give any absolute restrictions on the size of $m_{i}$
and $\tilde{m}_{j}$, e.g., a bound on the inner product%
\index{product!inner} of two vectors in
$\mathbb{C}^{N}$ does not give a bound on the size of the vectors if they are
not equal. This is reflected in the bi-Cuntz%
\index{bi Cuntz relations@bi-Cuntz relations} relations defined by $m_{i}$,
$\tilde{m}_{i}$. 
%In analogy with (\ref{eqBio.13}), 
Let us now define%
\begin{equation}
\left(  S_{i}f\right)  \left(  z\right)  =m_{i}\left(  z\right)  f\left(
z^{N}\right)  ,\qquad\left(  \smash{\tilde{S}_{i}}f\right)  \left(  z\right)
=\tilde{m}_{i}\left(  z\right)  f\left(  z^{N}\right)  \label{eqBio.16}%
\end{equation}
for $z\in\mathbb{T}$, $f\in L^{2}\left(  \mathbb{T}\right)  $. Instead of the
usual Cuntz%
\index{Cuntz relations} relations, 
%(\ref{eqMul.14}), 
the $S_{i}$, $\tilde{S}_{i}$ now
satisfy%
\begin{align}
S_{i}^{\ast}\tilde{S}_{j}  &  =\delta_{i,j}1,\label{eqBio.17}\\
\sum_{i}S_{i}\tilde{S}_{i}^{\ast}  &  =1. \label{eqBio.18}%
\end{align}
If $A,\tilde{A}\in C\left(  \mathbb{T},\mathrm{GL}_{N}\left(  \mathbb{C}%
\right)  \right)  $ are the matrix%
\index{matrix!function@--- function}-valued functions associated to $m_{i}$ and
$\tilde{m}_{i}$
by
\begin{equation}
m\left(  z\right)  =A\left(  z^{N}\right)  v\left(  z\right)  ,\qquad\tilde
{m}\left(  z\right)  =\tilde{A}\left(  z^{N}\right)  v\left(  z\right)  ,
\label{eqBio.10}%
\end{equation}
we compute%
\begin{equation}
S_{i}^{\ast}S_{j}=\left(  AA^{\ast}\right)  _{j,k} \label{eqBio.19}%
\end{equation}
in the sense that $S_{i}^{\ast}S_{j}$ is contained in the commutative algebra
of multiplication operator%
\index{operator!multiplication}s on $L^{2}\left(  \mathbb{T}\right)  $ defined by
$C\left(  \mathbb{T}\right)  $, and $\left(  AA^{\ast}\right)  _{j,i}\in
C\left(  \mathbb{T}\right)  $. Correspondingly,%
\begin{equation}
\tilde{S}_{i}^{\ast}\tilde{S}_{j}=\left(  \smash{\tilde{A}\tilde{A}^{\ast }%
}\right)  _{j,i} \label{eqBio.20}%
\end{equation}
so all the operators $S_{i}^{\ast}S_{j}$, $\tilde{S}_{i}^{\ast}\tilde{S}_{j}$
are contained in the abelian algebra $C\left(  \mathbb{T}\right)  $. We may
introduce operators $S$, $\tilde{S}$ from
\begin{equation}
L^{2}\left(  \mathbb{T}\right)  ^{N}=\underset{0}{L^{2}\left(  \mathbb{T}%
\right)  }\oplus\dots\oplus\underset{N-1}{L^{2}\left(  \mathbb{T}\right)  }
\label{eqBio.21}%
\end{equation}
into $L^{2}\left(  \mathbb{T}\right)  $ by%
\begin{equation}
S=\left(  S_{0},S_{1},\dots,S_{N-1}\right)  ,\qquad\tilde{S}=\left(
\smash{\tilde{S}_{0},\dots ,\tilde{S}_{N-1}}\right)  \label{eqBio.22}%
\end{equation}
and then $S^{\ast}$ maps $L^{2}\left(  \mathbb{T}\right)  $ into
(\ref{eqBio.21}), etc., and the relations (\ref{eqBio.17})\textendash
(\ref{eqBio.20}) take the form%
\begin{align}
&  \left\{
\begin{array}
[c]{ll}%
S^{\ast}\tilde{S}=1, & \text{where }1\text{ is the identity in }M_{N}\left(
\mathbb{C}\right)  \otimes C\left(  \mathbb{T}\right)  ,\\
S\tilde{S}^{\ast}=1, & \text{where }1\text{ is the identity in }C\left(
\mathbb{T}\right)  ,
\end{array}
\right. \label{eqBio.23}\\
&  \left\{
\begin{array}
[c]{l}%
S^{\ast}S=AA^{\ast},\\
\tilde{S}^{\ast}\tilde{S}=\tilde{A}\tilde{A}^{\ast}.
\end{array}
\right.  \label{eqBio.24}%
\end{align}
These relations say that all combinations of
product%
\index{product!operator}s
of $S$ and $S^{\ast}$
with $\tilde{S}$ and $\tilde{S}^{\ast}$ lie in the algebra $M_{N}\left(
\mathbb{C}\right)  \otimes C\left(  \mathbb{T}\right)  $. But in addition $A$
and $\tilde{A}$ are matrix%
\index{matrix!function@--- function}-valued functions on $\mathbb{T}$, so
%(\ref{eqBio.11}) implies%
\begin{equation}
AA^{\ast}\tilde{A}\tilde{A}^{\ast}=A\tilde{A}^{\ast}=1=\tilde{A}\tilde
{A}^{\ast}AA^{\ast} \label{eqBio.25}%
\end{equation}
and hence%
\begin{equation}
S^{\ast}S=\left(  \tilde{S}^{\ast}\tilde{S}\right)  ^{-1} \label{eqBio.26}%
\end{equation}
and all the matrix%
\index{matrix!function@--- function}-valued functions commute.

This discussion can be summarized by saying that the bi-Cuntz%
\index{bi Cuntz relations@bi-Cuntz relations} relations are
much less rigid than the original Cuntz%
\index{Cuntz relations} relations, i.e.:

\begin{scholium}
\label{SchBio.1}Given \emph{any} bijective operator $S$ from $L^{2}\left(
\mathbb{T}\right)  ^{N}$ into $L^{2}\left(  \mathbb{T}\right)  $ one may
define $\tilde{S}=\left(  S^{\ast}\right)  ^{-1}$ and the bi-Cuntz%
\index{bi Cuntz relations@bi-Cuntz relations} relations
\textup{(\ref{eqBio.23})} are satisfied. If, more specifically, $S$ is given by
\textup{(\ref{eqBio.22})} and \textup{(\ref{eqBio.16})}, then operators
$\tilde{S}_{0},\dots,\tilde{S}_{N-1}$ exist such that the bi-Cuntz%
\index{bi Cuntz relations@bi-Cuntz relations} relations
\textup{(\ref{eqBio.17})\textendash (\ref{eqBio.18})} are satisfied if and
only if the operator%
\index{operator!invertible} $A\in M_{N}\left(  \mathbb{C}\right)  \otimes C\left(
\mathbb{T}\right)  $ defined by 
 (\ref{eqBio.10})
is invertible, in
which case one must use $\tilde{A}=\left(  A^{\ast}\right)  ^{-1}$,
\textup{(\ref{eqBio.10})}, and \textup{(\ref{eqBio.16})} to define $\tilde
{S}_{0},\dots,\tilde{S}_{N-1}$.
\end{scholium}

Let us now connect the filter%
\index{filter!wavelet}s to the wavelets. We have already defined the
scaling%
\index{scaling!function@--- function} functions $\varphi$, $\tilde{\varphi}$ and wavelet%
\index{wavelet!function@--- function} functions $\psi
_{i}^{{}}$, $\psi_{i}^{{}}$, $i=1,\dots,N$. 
%by (\ref{eqBio.3})\textendash
%(\ref{eqBio.4}). The further discussion follows the lines from (\ref{eqMul.15}%
%) in Chapter \ref{ChHom}. 
The expansions for $\varphi$ and $\tilde{\varphi}$ 
%in
%(\ref{eqBio.3})
converge%
\index{convergence!uniform} uniformly on compacts, thus
$\hat{\varphi}$ and $%
\Hat{\Tilde{\varphi}}%
$ are continuous functions on $\mathbb{R}$. To decide that these functions are
in $L^{2}\left(  \mathbb{R}\right)  $ one again forms%
\begin{equation}
f_{\varphi}\left(  t\right)  =\operatorname*{PER}\left(  \left\vert
\hat{\varphi}\right\vert ^{2}\right)  \left(  t\right)  =\sum_{n\in\mathbb{Z}%
}\left\vert \hat{\varphi}\left(  t-2\pi n\right)  \right\vert ^{2}%
\label{eqBio.27}%
\end{equation}
and $f_{\tilde{\varphi}}$ similarly, and one deduces again from
the nonlinear intertwining relation%
\begin{equation}
R^{k}\left(  p\left(  \psi_{1},\psi_{2}\right)  \right)  =p\left(  M_{m_{0}%
}^{k}\psi_{1},M_{m_{0}}^{k}\psi_{2}\right)  ,\qquad k\in\mathbb{N}
\label{eqCas.7}%
\end{equation}
%(\ref{eqCas.7})
that%
\begin{equation}
R_{m_{0}}\left(  f_{\varphi}\right)  =f_{\varphi},\qquad R_{\tilde{m}_{0}%
}\left(  f_{\tilde{\varphi}}\right)  =f_{\tilde{\varphi}}.\label{eqBio.28}%
\end{equation}

\subsubsection{\label{S2.2bis.2}General results in mathematics on matrix functions}

\label{ChExaTut}
In the standard case of the good old
orthogonal%
\index{orthogonal!wavelet@--- wavelet} wavelets in $L^{2}\left(  \mathbb{R}\right)  $ of $N$ subbands, we
will look for functions $\psi_{1},\dots,\psi_{N-1}$ in $L^{2}\left(
\mathbb{R}\right)  $ such that, if $k$ and $n$ run independently over all the
integers $\mathbb{Z}$, i.e., $-\infty<k,n<\infty$, then the countably infinite
system of functions%
\begin{equation}
\left\{  \smash{N^{k/2}}\psi_{i}\left(  \smash{N^{k}}x-n\right)
\bigm|i=1,\dots,N-1,\;k,n\in\mathbb{Z}\right\}  \label{eqChExaTut.1}%
\end{equation}
is an \emph{orthonormal%
\index{orthonormal!basis@--- basis|textbf} basis} in the Hilbert%
\index{Hilbert!space@--- space} space $L^{2}\left(
\mathbb{R}\right)  $. The second half of the
%\linebreak
word \textquotedblleft
orthonormal\textquotedblright \ refers to the restricting requirement that all
the functions
%\linebreak
$\psi_{1},\dots,\psi_{N-1}$ satisfy%
\begin{equation}
\int_{\mathbb{R}}\left\vert \psi_{i}\left(  x\right)  \right\vert ^{2}\,dx=1,
\label{eqChExaTut.2}%
\end{equation}
or stated more briefly,
\begin{equation}
\left\Vert \psi_{i}\right\Vert _{L^{2}\left(
\mathbb{R}\right)  }=1; 
\label{eqChExaTut.2bis}%
\end{equation}
or yet more briefly,
\begin{equation}
\left\Vert \psi_{i}\right\Vert
=1. 
\label{eqChExaTut.2ter}%
\end{equation}
{}From familiar properties of the Lebesgue measure on $\mathbb{R}$, it then
follows that all the functions%
\begin{equation}
\psi_{i,k,n}\left(  x\right)  :=N^{k/2}\psi_{i}\left(  \smash{N^{k}%
}x-n\right)  ,\qquad1\leq i<N,\;k,n\in\mathbb{Z}, \label{eqChExaTut.3}%
\end{equation}
satisfy the normalization, i.e., that%
\begin{equation}
\left\Vert \psi_{i,k,n}\right\Vert =1\text{\qquad for all }i,k,n.
\label{eqChExaTut.4}%
\end{equation}
The functions (\ref{eqChExaTut.3}) are said to be \emph{orthogonal} if%
\begin{equation}
\int_{\mathbb{R}}\overline{\psi_{i,k,n}\left(  x\right)  }\,\psi_{i^{\prime
},k^{\prime},n^{\prime}}\left(  x\right)  \,dx=0 \label{eqChExaTut.5}%
\end{equation}
whenever $\left(  i,k,n\right)  \neq\left(  i^{\prime},k^{\prime},n^{\prime
}\right)  $. We say that the two triple indices are different if $i\neq
i^{\prime}$ or $k\neq k^{\prime}$ or $n\neq n^{\prime}$. If, for example,
$i=i^{\prime}$ and $k=k^{\prime}$, then when the same function is translated
by different amounts $n$ and $n^{\prime}$, the two resulting functions are
required to be orthogonal%
\index{orthogonal!translates@--- translates}. It is an elementary geometric fact from
the theory of
Hilbert%
\index{Hilbert!space@--- space} space
that if the functions in (\ref{eqChExaTut.3}) form an
orthonormal%
\index{orthonormal!basis@--- basis} basis, then for every function $f\in L^{2}\left(  \mathbb{R}%
\right)  $, i.e., every measurable function $f$ on $\mathbb{R}$ such that
\begin{equation}
\left\Vert f\right\Vert ^{2}=\int_{\mathbb{R}}\left\vert f\left(  x\right)
\right\vert ^{2}\,dx<\infty , 
\label{eqChExaTut.6pre}%
\end{equation}
we have the identity%
\begin{equation}
\left\Vert f\right\Vert ^{2}=\sum_{i,k,n}\left\vert \int_{\mathbb{R}}%
\overline{\psi_{i,k,n}\left(  x\right)  }\,f\left(  x\right)  \,dx\right\vert
^{2}, \label{eqChExaTut.6}%
\end{equation}
where the triple summation in (\ref{eqChExaTut.6}) is over all configurations
$1\leq i<N$, $k,n\in\mathbb{Z}$. It is convenient to rewrite
(\ref{eqChExaTut.6}) in the following more compact form:%
\begin{equation}
\left\Vert f\right\Vert ^{2}=\sum_{i,k,n}\left\vert \left\langle\, \psi _{i,k,n}\mid
f\,\right\rangle \right\vert ^{2}. \label{eqChExaTut.7}%
\end{equation}
Surprisingly, it turns out that (\ref{eqChExaTut.7}) may hold even if the
functions $\psi_{i,k,n}$ of (\ref{eqChExaTut.3}) do not form an orthonormal
basis. It may happen that one of the initial functions
%\linebreak
$\psi_{1}$, $\dots$, or
$\psi_{N-1}$ satisfies $\left\Vert \psi_{i}\right\Vert <1$, and yet that
(\ref{eqChExaTut.7}) holds for all $f\in L^{2}\left(  \mathbb{R}\right)  $.
These more general systems are still called wavelets, but since they are
special, they are referred to as \emph{tight frame%
\index{frame!tight|textbf}s}, as opposed to
orthonormal bases. In either case, we will talk about a \emph{wavelet%
\index{wavelet!expansion@--- expansion|textbf}
expansion} of the form%
\begin{equation}
f\left(  x\right)  =\sum_{i,k,n}\left\langle\, \psi _{i,k,n}\mid f\,\right\rangle \psi_{i,k,n}\left(
x\right)  . \label{eqChExaTut.8}%
\end{equation}
It follows that the sum on the right-hand side in (\ref{eqChExaTut.8})
converge%
\index{convergence!L two@$L\sp {2}$-}s in the norm of $L^{2}\left(  \mathbb{R}\right)  $ for all functions
$f$ in $L^{2}\left(  \mathbb{R}\right)  $ if (\ref{eqChExaTut.7}) holds.

But 
%as mentioned in the introduction,
there is a yet more general form
of wavelets, called \emph{biorthogonal%
\index{wavelet!biorthogonal|textbf}}. The conditions on the functions
$\psi_{1},\dots,\psi_{N-1}$ are then much less restrictive than the
orthogonality axioms. Hence these wavelets are more flexible and adapt better
to a variety of applications, for example, to data compression, or to computer
graphics. But the biorthogonality conditions are also a little more technical
to state. We say that some given functions $\psi_{i}$, $i=1,\dots,N-1$, in
$L^{2}\left(  \mathbb{R}\right)  $ are part of a biorthogonal wavelet%
\index{wavelet!biorthogonal} system
if there is a second system of functions $\tilde{\psi}_{i}$, $i=1,\dots,N-1$,
in $L^{2}\left(  \mathbb{R}\right)  $, such that every $f\in L^{2}\left(
\mathbb{R}\right)  $ admits a representation%
\begin{equation}
f\left(  x\right)  =\sum_{i,k,n}\left\langle\, \psi _{i,k,n}\mid f\,\right\rangle \tilde{\psi}%
_{i,k,n}\left(  x\right)  =\sum_{i,k,n}\left\langle\, \smash{\tilde{\psi}}\vphantom{\psi}
_{i,k,n}\mid f\,\right\rangle \psi_{i,k,n}\left(  x\right)  , \label{eqChExaTut.9}%
\end{equation}
and%
\begin{equation}
\tilde{\psi}_{i,k,n}\left(  x\right)  =N^{k/2}\tilde{\psi}_{i}\left(
\smash{N^{k}}x-n\right)  . \label{eqChExaTut.10}%
\end{equation}
In the standard normalized case where $\left\langle\, \psi_{i}\mid\smash{\tilde{\psi}%
}\vphantom{\psi} _{i}\,\right\rangle =1$, then you will notice that condition
(\ref{eqChExaTut.7}) turns into%
\begin{equation}
\left\Vert f\right\Vert ^{2}=\sum_{i,k,n}\overline{\left\langle\, \psi _{i,k,n}\mid f\,\right\rangle 
}\,\left\langle\, \smash{\tilde{\psi}}\vphantom{\psi} _{i,k,n}\mid f\,\right\rangle \label{eqChExaTut.11}%
\end{equation}
for all
$f\in L^{2}\left(  \mathbb{R}\right)  $.

%As you recall, we demonstrated 
%that 
The orthogonal%
\index{orthogonal!wavelet@--- wavelet} wavelets correspond
to matrix%
\index{matrix!function@--- function} functions $\mathbb{T}\rightarrow\mathrm{U}_{N}\left(  \mathbb{C}%
\right)  $, while the wider class of biorthogonal wavelet%
\index{wavelet!biorthogonal}s corresponds to the
much bigger group of matrix%
\index{matrix!function@--- function} functions $\mathbb{T}\rightarrow\mathrm{GL}%
_{N}\left(  \mathbb{C}\right)  $, via the associated wavelet filter%
\index{filter!wavelet}s. You may
ask, why bother with the more technical-looking biorthogonal systems? It turns
out that they are forced on us by the engineers. They tell us that the real
world is not nearly as orthogonal as the mathematicians would like to make it
out to be.
There is a paucity of symmetric orthogonal%
\index{orthogonal!wavelet@--- wavelet} wavelets,
and symmetry (``linear phase'')
is prized by engineers and workers in image processing,
where the more general wavelet families and their duality play a crucial role.
Now what if we could change the biorthogonal wavelet%
\index{wavelet!biorthogonal}s into the
orthogonal%
\index{orthogonal!wavelet@--- wavelet} ones, and still keep the essential spectral%
\index{spectral!properties of biorthogonal wavelet@--- properties of biorthogonal wavelet} properties intact? Then
everyone will be happy. This last chapter shows that it is possible, and even
in a fairly algorithmic fashion, one that is amenable to computations.

%As we saw in Chapters \ref{Int} and \ref{ChHom}, 
Wavelet
filter%
\index{filter!wavelet}s may be understood as matrix%
\index{matrix!function@--- function}
functions, i.e., functions from the one-torus
$\mathbb{T}\subset\mathbb{C}$ into some group of invertible
matrices%
\index{matrix!function@--- function}%
\index{matrix!invertible}. If the scale%
\index{scaling!number@--- number} number is
$N$, then there are three such
matrix groups which are especially
relevant for wavelet analysis%
\index{analysis!wavelet}:
\[
\framebox[0.24\textwidth]{\parbox[t]{0.22\textwidth}{\raggedright
%\epigraph 
$\mathrm{U}_{N}\left( \mathbb{C}\right) $: all unitary
$N\times N$ complex matrices%
\index{matrix!unitary}}}
\subset
\framebox[0.24\textwidth]{\parbox{0.22\textwidth}{\raggedright
%\epigraph 
$\mathrm{GL}_{N}\left( \mathbb{C}\right) $: all invertible
$N\times N$ complex matrices%
\index{matrix!invertible}}}
\supset
\framebox[0.24\textwidth]{\parbox[t]{0.22\textwidth}{\raggedright
%\epigraph 
$\mathrm{SL}_{N}\left( \mathbb{C}\right) $: all
$N\times N$ complex matrices $A$ with $\det A=1$.}}
\]
It is possible to reduce
some questions in the $\mathrm{GL}_{N}$ case to better
understood results for $\mathrm{U}_{N}\left( \mathbb{C}\right) $;
see Chapter 6 of \cite{BrJo02b}. The
$\mathrm{SL}_{2}$ case
is especially interesting in view of Daubechies%
\index{Daubechies1@Daubechies!Sweldens lifting@--- --Sweldens lifting}--Sweldens
lifting for dyadic wavelet%
\index{wavelet!dyadic}s; see
%Exercises 
\S\ \ref{Exe6DauSwe-11}.
% and \ref{Exe6Nov19-12}.

\subsubsection{\label{S2.2bis.3}Connection between matrix functions and wavelets}

\noindent\textbf{Definitions:}
A function, or a distribution%
\index{distribution}, $\varphi$ satisfying (\ref{eqInt.1}) is said to
be \emph{refinable}, the equation (\ref{eqInt.1}) is called the
\emph{refinement%
\index{refinement!equation@--- equation} equation}, or also, as noted above, the
``scaling%
\index{scaling!identity@--- identity} identity'', and $\varphi$ is called the scaling%
\index{scaling!function@--- function} function.
The coefficients $a_{n}$ of (\ref{eqInt.1}) are called the
\emph{masking%
\index{masking coefficients} coefficients}. 
%(see Definition \ref{DefRTOMar9.8}).

We will
mainly concentrate on the case when the set $\left\{  a_{n}\right\}  $ is
finite. But in general, a function $\varphi\in L^{2}\left(  \mathbb{R}\right)
$ is said to be refinable with scale%
\index{scaling!number@--- number} number $N$ if
$\varphi\left(  x/N\right)  $ is in the $L^{2}$-closed linear span of the
translates $\left\{  \varphi\left(  x-k\right)  \right\}  _{k\in\mathbb{Z}%
}\subset L^{2}\left(  \mathbb{R}\right)  $; see, e.g.,
\cite{HSS96,SSZ99,StZh98,StZh01}%
\index{Strang0@Gilbert Strang}.

Since there are refinement operations which are
more general than scaling \textup{(}see for example \textup{\cite{DLLP01})},
there are variations of \textup{(\ref{eqInt.1})} which are correspondingly
more general, with regard to both the refinement steps that are
used and the dimension of the spaces. The term
``scaling%
\index{scaling!identity@--- identity} identity'' is usually, but not always, reserved
for \textup{(\ref{eqInt.1})}, while more general refinement%
\index{operator!refinement}s lead to
``refinement%
\index{refinement!equation@--- equation} equations''. However, \textup{(\ref{eqInt.1})} often
goes under both names. The vector versions of the
identities get the prefix ``multi-'', for example \emph{multiscaling}
and \emph{multiwavelet%
\index{multiwavelet}}.
%; see Section \textup{\ref{Ope}}.
%\end{definitions}

If $m_{0}$ satisfies a condition for obtaining orthogonal wavelet%
\index{orthogonal!wavelet@--- wavelet}s,
\begin{equation}
\sum_{w^{N}=z}\left\vert m_{0}\left(  w\right)  \right\vert ^{2}=N,
\label{eqIntNew.4}%
\end{equation}
together with the normalization%
\begin{equation}
m_{0}\left(  1\right)  =\sqrt{N}, \label{eqIntNew.5}%
\end{equation}
then (\ref{eqInt.1}) has a solution $\varphi$ in $L^{2}\left(  \mathbb{R}%
\right)  $ which can be obtained by taking the inverse Fourier%
\index{Fourier1@Fourier!inverse transform@inverse --- transform} transform of
the product%
\index{product!infinite representation@infinite --- representation} expansion%
\begin{equation}
\hat{\varphi}\left(  t\right)  =\prod_{k=1}^{\infty}\left(  \frac{m_{0}\left(
tN^{-k}\right)  }{\sqrt{N}}\right)  .\label{eqIntNew.6}%
\end{equation}
(Here and later we use the convention that
if $m\left(  z\right)  $ is a function of $z\in\mathbb{T}$, then
$m\left( t\right)  =m\left(  e^{-it}\right)  $.)
That (\ref{eqIntNew.6}) gives a solution $\varphi$
of (\ref{eqInt.1})
follows from the relation
\begin{equation}
\hat{\varphi}\left(  t\right)  =
\frac{1}{\sqrt{N}}
m_{0}\left( \frac{t}{N}\right)  
\hat{\varphi}\left( \frac{t}{N}\right)  .
\label{eqIntAug15.15}
\end{equation}
%which is equivalent to (\ref{eqInt.1}).

\subsubsection*{\label{S2.2bis.3.1}\textup{\ref{S2.2bis.3}.1}\enspace Multiresolution wavelets}

We mentioned that there is a direct connection between $m_{0}=\sum a_{n}z^{n}$
and the scaling%
\index{scaling!function@--- function} function $\varphi$ on $\mathbb{R}$ given in (\ref{eqIntNew.2}%
),
%\textendash 
(\ref{eqInt.1}), and (\ref{eqIntNew.6}). There is a similar
correspondence between the high-pass filters%
\index{filter!high-pass} $m_{i}$ and the wavelet%
\index{wavelet!generator@--- generator} generators $\psi_{i}\in L^{2}\left(  \mathbb{R}\right)  $. In
the \emph{biorthogonal%
\index{biorthogonal}} case, there is a second system $\tilde{m}_{i}%
\leftrightarrow\tilde{\psi}_{i}$ and the two systems
\begin{multline}
\left\{  N^{\frac{j}{2}}\psi_{i}\left(  N^{j}x-k\right)  \right\}  \text{\quad
and\quad}\left\{  N^{\frac{j^{\prime}}{2}}\tilde{\psi}_{i^{\prime}}\left(
N^{j^{\prime}}x-k^{\prime}\right)  \right\}  ,\label{eqIntFeb22.10}\\
i,i^{\prime}\in\left\{  1,2,\dots,N-1\right\}  ,\;j,j^{\prime},k,k^{\prime}%
\in\mathbb{Z},
\end{multline}
then form a dual wavelet%
\index{wavelet!basis@--- basis} basis%
\index{dual wavelet basis}, or dual wavelet frame%
\index{frame} for $L^{2}\left(  \mathbb{R}\right)  $ in the sense of
\cite{Dau92}%
\index{Daubechies0@Ingrid Daubechies}, Chapter 5.
We considered this biorthogonal%
\index{biorthogonal} case in more detail in 
\S\ \ref{Bio} above.
Much more detail can be found in Chapter 6 of \cite{BrJo02b}.

The idea of constructing maximally smooth%
\index{vanishing moments} wavelets when some side conditions
are specified has been central to much of the activity in wavelet analysis%
\index{analysis!wavelet} and
its applications since the mid-1980's. As a supplement to \cite{Dau92}%
\index{Daubechies0@Ingrid Daubechies}, the survey
article \cite{Stra93}%
\index{Strang0@Gilbert Strang} is enjoyable reading.
The paper \cite{LaHe96} treats the issue in a more specialized setting and is
focussed on the moment method. Some of the early applications to data
compression and image coding are done very nicely in \cite{HSS+95}%
\index{Strang0@Gilbert Strang},
\cite{SHS+99}%
\index{Strang0@Gilbert Strang}, and \cite{HSW95}. An interesting, related but different,
algebraic and geometric approach to the problem is offered in \cite{PeWi99}%
\index{Wickerhauser0@M.V. Wickerhauser}.

We now turn to an
interesting variation of this
setup, which includes higher dimensions, i.e.,
when the Hilbert%
\index{Hilbert!space@--- space} space is
$L^{2}\left( \mathbb{R}^{d}\right) $, $d=2,3,\dots$. Staying for the moment with
$d=1$, and $N$ fixed, we will take
the viewpoint of what is called \emph{resolution%
\index{resolution|textbf}s}, but here
understood in a broad sense of closed subspaces: A
closed linear subspace $\mathcal{V}\subset L^{2}\left( \mathbb{R}\right) $ is
said to be an $N$-resolution%
\index{resolution} if
it is invariant under the unitary operator%
\index{operator!unitary}%
\index{operator!scaling}
\begin{equation}
U=U_{N}\colon f\longmapsto N^{-\frac{1}{2}}f\left( \frac{x}{N}\right) ,
\label{eqIntOct16.U}
\end{equation}
i.e., if $U$ maps $\mathcal{V}$ into a proper
subspace of itself. The subspace
$\mathcal{V}$ is said to be \emph{translation
invariant} if
\begin{equation}
f\in\mathcal{V}\iff f\left( \,\cdot \,-k\right) \in\mathcal{V}
\text{\qquad for all }k\in\mathbb{Z}.
\label{eqIntOct16.7}
\end{equation}
If there is a function $\varphi$
such that $\mathcal{V}=\mathcal{V}_{\varphi}$ is the closed
linear span of
\begin{equation}
\left\{ \varphi\left( \,\cdot \,-k\right) \mid k\in\mathbb{Z}\right\} ,
\label{eqIntOct16.8}
\end{equation}
then clearly $\mathcal{V}$ is translation invariant.
The translation-invariant resolution%
\index{resolution!subspace@--- subspace}
subspaces $\mathcal{V}$ are actively studied and reasonably well
understood. If $\mathcal{V}$ is
of the form $\mathcal{V}_{\varphi}$ in (\ref{eqIntOct16.8}), then
we say that it is \emph{singly
generated}, and that $\varphi$
is a scaling%
\index{scaling!function@--- function} function of scale%
\index{scaling!number@--- number}
$N$. 

\subsubsection*{\label{S2.2bis.3.2}\textup{\ref{S2.2bis.3}.2}\enspace Generalized multiresolutions [joint work with L.~Baggett, K.~Merrill,
and J.~Packer]}

The case when
the resolution subspace
$\mathcal{V}$ is not singly generated
is also interesting, and these
resolution%
\index{resolution!subspace@--- subspace} subspaces are frequently
called \emph{generalized multiresolution%
\index{multiresolution!generalized}%
\index{multiresolution!subspace@--- subspace}
subspaces} (GMRA). There is much current
and very active research on them;
see, for example, \cite{BaLa99}, \cite{LPT01}%
\index{Packer0@Judy Packer},
\cite{BaMe99}, \cite{HLPS99}, \cite{HSS01},
\cite{SSZ99}%
\index{Strang0@Gilbert Strang}, and \cite{Jor01a}. The case
when $\mathcal{V}$ is not singly generated
as a resolution%
\index{resolution!subspace@--- subspace} subspace of scale%
\index{scaling!number@--- number}
$N>2$, i.e., when $\mathcal{V}$ is not
of the form (\ref{eqIntOct16.8}), occurs
in the study of \emph{wavelet%
\index{wavelet!set@--- set} sets}.
A wavelet set in $\mathbb{R}^{d}$ is defined relative to an expansive
$d\times d$ matrix $\mathbf{N}$ over $\mathbb{Z}$. A subset $E\subset
\mathbb{R}^{d}$ is said to be an $\mathbf{N}$-wavelet set if there is a single
wavelet function $\psi\in L^{2}\left(  \mathbb{R}^{d}\right)  $ such that
$\hat{\psi}=\chi_{E}$. Specifically, the condition states that the family 
\begin{equation}
\left\{  \,\left\vert \det\mathbf{N}\right\vert ^{j/2}\psi\left(
\mathbf{N}^{j}x-\mathbf{k}\right)  :j\in\mathbb{Z},\;\mathbf{k}\in
\mathbb{Z}^{d}\,\right\}  \label{eq2.2bis.3.1}%
\end{equation}
is an orthonormal basis for $L^{2}\left(  \mathbb{R}^{d}\right)  $. This can
be checked to be equivalent to the combined set of two tiling properties for
$E$ as a subset of $\mathbb{R}^{d}$:

\begin{enumerate}
\renewcommand{\theenumi}{\alph{enumi}}
\item \label{S2.2bis.3.2(1)}the family of subsets $\left\{  \,\mathbf{N}^{j}E:j\in\mathbb{Z}%
\,\right\}  $ tiles $\mathbb{R}^{d}$;

\item \label{S2.2bis.3.2(2)}the translates $\left\{  \,E+2\pi\mathbf{k}:\mathbf{k}\in\mathbb{Z}%
^{d}\,\right\}  $ tile $\mathbb{R}^{d}$.
\end{enumerate}

We define tiling by the requirement that the sets in the family have overlap
at most of measure zero relative to Lebesgue measure on $\mathbb{R}^{d}$.
Similarly, the union%
\begin{equation}
\mathbb{R}^{d}=\bigcup_{j\in\mathbb{Z}}\mathbf{N}^{j}E=\bigcup_{\mathbf{k}%
\in\mathbb{Z}^{d}}E+2\pi\mathbf{k}\label{eq2.2bis.3.2}%
\end{equation}
is understood to be only up to measure zero.

It is easy to see that compactly supported wavelets in $L^{2}\left(
\mathbb{R}^{d}\right)  $ are MRA wavelets, while most wavelets $\psi=\left(
\chi_{E}\right)  \spcheck$ from wavelet sets $E$ are not. Thess wavelets are
typically (but not always) frequency localized.

The main difference between the GMRA (stands for generalized multiresolution
analysis) wavelets and the more traditional MRA ones may be understood in
terms of multiplicity. Both come from a fixed resolution subspace
$\mathcal{V}_{0}\subset L^{2}\left(  \mathbb{R}^{d}\right)  $ which is
invariant under the translations $\left\{  \,T_{n}:n\in\mathbb{Z}%
^{d}\,\right\}  $ where
\begin{equation}
\left(  T_{n}f\right)  \left(  x\right)  :=f\left(  x-n\right)  \text{\qquad
for }x\in\mathbb{R}^{d}\text{ and }n\in\mathbb{Z}^{d}.\label{eq2.2bis.3.3}%
\end{equation}
Hence $\left\{  T_{n}|_{\mathcal{V}_{0}}\right\}  _{n\in\mathbb{Z}^{d}}$ is a
unitary representation of $\mathbb{Z}^{d}$ on the Hilbert space $\mathcal{V}%
_{0}$. As a result of Stone's theorem, we find that there are subsets%
\[
E_{1}\supset E_{2}\supset\dots\supset E_{j}\supset\cdots
\]
of $\mathbb{T}^{d}$ such that the spectral measure of the (restricted)
representation has multiplicity $\geq j$ on the subset $E_{j}$, $j=1,2,\dots$.
It can be checked that the projection-valued spectral measure is absolutely
continuous. Moreover, there is an intertwining unitary operator%
\begin{equation}
J\colon\mathcal{V}_{0}\longrightarrow
\sideset{}{^{\smash{\oplus}}}{\sum}\limits_{j\geq1}L^{2}\left(  E_{j}\right)
\label{eq2.2bis.3.4}%
\end{equation}
such that%
\begin{equation}
P_{L^{2}\left(  E_{j}\right)  }JT_{n}f\left(  z\right)  =z^{n}\left(
Jf\right)  \left(  z\right)  \label{eq2.2bis.3.5}%
\end{equation}
holds for all $f\in\mathcal{V}_{0}$ and $z\in E_{j}$. We may then consider the
functions $\varphi_{j}\in\mathcal{V}_{0}$ ($\subset L^{2}\left(
\mathbb{R}^{d}\right)  $) defined by
\begin{equation}
\varphi_{j}:=J^{-1}(0,\dots,0,\underbrace{\,\chi_{E_{j}}\,}%
_{\!\!\!\!j\text{'th place}\!\!\!\!},0,0,\dots).\label{eq2.2bis.3.6}%
\end{equation}
It was proved by Baggett and Merrill \cite{BaMe99} that $\left\{
\,\varphi_{j}:j\geq1\,\right\}  $ generates a normalized tight frame for
$\mathcal{V}_{0}$: specifically, that%
\begin{equation}
\sum_{j\geq1\vphantom{n\in \mathbb{Z}^{d}}}\sum_{n\in\mathbb{Z}^{d}}\left\vert
\left\langle \,T_{n}\varphi_{j}\mid f\,\right\rangle \right\vert
^{2}=\left\Vert f\right\Vert _{L^{2}\left(  \mathbb{R}^{d}\right)  }%
^{2}\label{eq2.2bis.3.7}%
\end{equation}
holds for all $f\in\mathcal{V}_{0}$.

Treating $\left(  \varphi_{1},\varphi_{2},\dots\right)  $ as a vector-valued
function, denoted simply by $\varphi$, we see that there is a matrix function%
\[
H\colon\mathbb{T}^{d}\longrightarrow\left(  \text{complex square
matrices}\right)
\]
such that%
\begin{equation}
\hat{\varphi}\left(  \mathbf{N}^{\operatorname*{tr}%
\rlap{\!\raisebox{5pt}[0pt][0pt]{$\scriptstyle\nearrow$\raisebox{5pt}[0pt][0pt]{\small for transpose}}}}%
t\right)  =H\left(  e^{it}\right)  \hat{\varphi}\left(  t\right)
,\label{eq2.2bis.3.7bis}%
\end{equation}
where $t=\left(  t_{1},\dots,t_{d}\right)  \in\mathbb{R}^{d}$, and
$e^{it}:=\left(  e^{it_{1}},e^{it_{2}},\dots,e^{it_{d}}\right)  $.

But this method takes the Hilbert space $L^{2}\left(  \mathbb{R}^{d}\right)  $
as its starting point, and then proceeds to the construction of wavelet
filters in the form (\ref{eq2.2bis.3.7bis}). Our current joint work with
Baggett, Merrill, and Packer reverses this. It begins with a matrix function
$H$ defined on $\mathbb{T}^{d}$, and then offers subband conditions on the
matrix function which allow the construction of a GMRA for $L^{2}\left(
\mathbb{R}^{d}\right)  $ with generator $\varphi=\left(  \varphi_{1}%
,\varphi_{2},\dots\right)  $ given by (\ref{eq2.2bis.3.7bis}). So the Hilbert
space $L^{2}\left(  \mathbb{R}^{d}\right)  $ shows up only at the end of the
construction, in the conclusions of the theorems.

\subsubsection{\label{SNew5Aug2.3.4} Matrix completion}

In using the polyphase matrices, one may only have the first few rows, and be 
faced with the problem of completing to get the entire function $A$ from a torus 
into the matrices of the desired size. The case when only the first row is 
given, say corresponding to a specified low-pass filter, is treated in
\cite{BrJo02b} and \cite{BrJo02a}, and we refer the reader to the references given there, 
especially \cite{JiSh94}, \cite{RiSh91}, \cite{RiSh92}, and \cite{Vai93}.

The wavelet transfer operator is used in a variety of wavelet applications not
covered here, or only touched upon tangentially: stability of refinable
functions, regularity, approximation order, unitary matrix extension
principles, to mention only a few. The reader is referred to the following
papers for more details on these subjects:
\cite{DHRS03}, \cite{RoSh03}, \cite{RST01}, \cite{JJS01}, \cite{RoSh00}, \cite{RiSh00}, 
\cite{JiSh99}, \cite{She98}, \cite{RoSh98}, \cite{LLS98}, \cite{RoSh97}, \cite{BJMP03b},
and \cite{BJMP04}.

  The unitary extension principle (UEP) of \cite{DHRS03} involves the interplay 
between a finite set of filters (functions on $\mathbb{R}/\mathbb{Z}$),
and a corresponding tight 
frame (alias Parseval frame) in $L^2\left(\mathbb{R}^d\right)$.

  For the sake of illustration, let us take $d = 1$, and scaling 
number $N = 2$, i.e., the case of dyadic framelets. Naturally, the notion of 
tight frame is weaker than that of an orthonormal basis (ONB), and it is shown 
in \cite{DHRS03} that when a system of wavelet filters $m_i$, $i = 0, 1, \dots, r$ is 
given ($m_0$ must be low-pass), then the orthogonality condition on the $m_i$'s 
always gets us a framelet in $L^2\left(\mathbb{R}\right)$, i.e., the functions $\psi_i$ corresponding to 
the high-pass filters $m_i$, $i = 1,\dots,r$, generate a tight frame for $L^2\left(\mathbb{R}\right)$, 
also called a framelet. The correspondence $m_i$ to $\psi_i$ is called the UEP in 
\cite{DHRS03}.
 
  The orthogonality condition for $m_i$, $i = 0, 1, \dots, r$, referred to in the UEP 
is simply this: Form an $\left(r + 1\right)$-by-$2$ matrix-valued function $F\left(x\right)$ by using $m_i\left(x\right)$, 
$i = 0, 1, \dots, r$ in the first column, and the translates of the $m_i$'s by a half 
period, i.e., $m_i\left(x + 1/2\right)$, $i = 0, 1, \dots, r$ in the second. The condition on 
this matrix function $F\left(x\right)$ is that the two columns are orthogonal and have unit
norm in $\ell^2$ for all $x$. Note that we still get the unitary matrix functions 
acting on these systems, in the way we outlined above. But there is redundancy 
as the unitary matrices are $\left(r + 1\right)$-by-$\left(r + 1\right)$. The reader is refered to \cite{DHRS03} 
for further details.

  We emphasize that several of these, and other related topics, invite the kind 
of probabilistic tools that we have stressed here. But a more systematic 
discussion is outside the scope of this brief set of notes. 
We only hope to offer a modest introduction to a variety of more specialized 
topics.

\begin{rem}
\label{SNew5Aug2.3.4(a)}
The orthogonality condition for $m_i$, $i = 0, 1, \dots, r$, may be stated in 
terms of the operators $S_i$ from equation (\ref{eq2.9}), $N = 2$. For each $i = 0,1,\dots, r$, 
define an operator on $L^2\left(\mathbb{R}/\mathbb{Z}\right)$ as in (\ref{eq2.9}). Then the arguments from Section \ref{S2} 
show that the orthogonality condition for $m_i$, $i = 0, 1, \dots, r$, i.e., the UEP 
condition, is equivalent to the operator identity (\ref{eq2.8}) where the summation now 
runs from $0$ to $r$. Operator systems $S_i$ satisfying (\ref{eq2.8}) are called row-%
isometries.
\end{rem}

\begin{rem}
\label{SNew5Aug2.3.4(b)} There are two properties of the low-pass filter $m_0$ which we have glossed 
over. First, $m_0$ must be such that the corresponding scaling function $\varphi$ is in 
$L^2\left(\mathbb{R}\right)$. Without an added condition on $m_0$, $\varphi$ might only be a distribution. 
Secondly, when the dyadic scaling in $L^2\left(\mathbb{R}\right)$ is restricted to the resolution 
subspace $V_0\left(\varphi\right)$, the corresponding unitary part must be zero. These two 
issues are addressed in \cite{BJMP03b}, \cite{BJMP04}, and \cite{DHRS03}.
\end{rem}

\subsubsection{\label{S2.2bis.4} Connections between matrix functions and signal processing}

\label{Mat}
Since our joint work with Baggett, Merrill, and Packer on the GMRA wavelets is
still in progress, we restrict the discussion of matrix functions here to the
MRA case.

The two groups of matrix%
\index{matrix!function@--- function} functions $C\left(  \mathbb{T},\mathrm{U}_{N}\left(
\mathbb{C}\right)  \right)  $ and $C\left(  \mathbb{T},\mathrm{GL}_{N}\left(
\mathbb{C}\right)  \right)  $, i.e., the continuous functions from the torus
into the respective groups, enter wavelet analysis%
\index{analysis!wavelet} via the associated wavelet filter%
\index{filter!wavelet}s $\left(  m_{i}\right)  _{i=0}^{N-1}$.

In 
%Sections \ref{Hom} and \ref{Iso}, 
\cite{BrJo02b} (see also \S\ \ref{Hom} above),
we give the details of the multiple
correspondence between:
\begin{enumerate}
\renewcommand{\theenumi}{\roman{enumi}}

\item \label{IntMultCorr(1)}matrix%
\index{matrix!function@--- function} functions, $A\colon\mathbb{T}\rightarrow\mathrm{GL}_{N}\left(
\mathbb{C}\right)  $,

\item \label{IntMultCorr(2)}high- and low-pass wavelet filter%
\index{filter!wavelet}s%
\index{filter!high-pass}%
\index{filter!low-pass} $m_{i}$, $\tilde{m}_{i^{\prime}}$, $i,i^{\prime
}=0,1,\dots,N-1$, and

\item \label{IntMultCorr(3)}wavelet%
\index{wavelet!generator@--- generator} generators $\psi_{i}$, $\tilde{\psi}_{i^{\prime}}$,
$i,i^{\prime}=1,\dots,N-1$, together with scaling%
\index{scaling!function@--- function} functions
$\smash{
\varphi
}
$, $
\smash{
\tilde{\varphi}
}
$.
\end{enumerate}
\noindent In particular,%
\begin{align}
\qquad A_{i,j}\left(  z\right)   &  
=
\smash[t]{
\frac{1}{N}\sum_{w^{N}=z}m_{i}\left(
w\right)  w^{-j},
}
 &  &  z\in\mathbb{T},\qquad\label{eqIntFeb22.11}\\
\qquad\left(  A^{-1}\right)  _{i,j}  &  =\frac{1}{N}\sum_{w^{N}=z}%
\overline{\tilde{m}_{j}\left(  w\right)  }\,w^{i}, &  &  z\in\mathbb{T}%
.\qquad\label{eqIntFeb22.12}%
\end{align}
The dependence of the $L^{2}\left(  \mathbb{R}\right)  $-functions in
(\ref{IntMultCorr(3)}) on the group elements $A$ from (\ref{IntMultCorr(1)})
gives rise to homotopy%
\index{wavelet!homotopy of@homotopy of ---s} properties.
%, and the results in Sections \ref{PFe},
%\ref{Cyc}, \ref{SPB}, and \ref{Iso} are building up to that, while the final
%results are stated in Section \ref{Hom}. 
The standard orthogonal wavelet%
\index{orthogonal!wavelet@--- wavelet}s represent the special case when $m_{i}=\tilde
{m}_{i}$, or equivalently, $A\left(  z\right)  =\left(  \left(  A\left(
z\right)  \right)  ^{\ast}\right)  ^{-1}$, $z\in\mathbb{T}$. Hence, the matrix%
\index{matrix!function@--- function}
functions are unitary in this case.

The scaling%
\index{scaling!function@--- function}/wavelet function%
\index{wavelet!function@--- function}s $\varphi,\psi_{1},\dots,\psi_{N-1}$ with support on a
fixed compact interval, say $\left[  0,kN+1\right]  $, $k=0,1,\dots$, can be
parameterized with a finite number of parameters since the unitary%
\index{matrix!function@--- function}%
\index{matrix!unitary}-valued
function $z\rightarrow A\left(  z\right)  $ in (\ref{eqIntFeb22.11}) then is a
polynomial in $z$ of degree%
\index{degree!of a Fourier polynomial@--- of a Fourier polynomial} at
most $k\left(  N-1\right)  $. It is well-known
folklore from computer-generated pictures that the shape of the
scaling%
\index{scaling!function@--- function}/wavelet%
\index{wavelet!shape of} function%
\index{wavelet!function@--- function}s depends continuously on these parameters; see Figures
%\ref{FigChart1}--\ref{FigChartd} 
1.1--1.7 in \cite{BrJo02b} and \cite{Tre01b}.

The scaling%
\index{scaling!function@--- function} function $\varphi\in L^{2}\left(  \mathbb{R}\right)  $ of
(\ref{eqInt.1}) is illustrated
% in Figures \ref{FigChart1}--\ref{FigChartd}
there, in
the case $N=2$, and for orthogonal%
\index{orthogonal!translates@--- translates} $\mathbb{Z}$-translates, i.e., the case
(\ref{eqIntNew.4}). These pictures illustrate the dependence of $\varphi$ on
the masking%
\index{masking coefficients} coefficients $\left(  a_{n}\right)  $ in the case of
\cite{Tre01b}:%
\begin{equation}
\begin{aligned} a_{0} & =
%\frac{1}{4}%
(\eta_{0}-\eta_{1}-\eta_{2}+\eta_{3}+\eta_{4})
/4
,
\\
a_{1} & =
%\frac{1}{4}%
(\eta_{0}+\eta_{1}-\eta_{2}+\eta_{3}-\eta_{4})
/4
,\\
a_{2} & =
%\frac{1}{2}
(\eta{_0}-\eta_{3}-\eta_{4})
/2,
\\
a_{3} & =
%\frac{1}{2}
(\eta_{0}-\eta_{3}+\eta_{4})
/2,\\
a_{4} & =
%\frac{1}{4}%
(\eta_{0}+\eta_{1}+\eta_{2}+\eta_{3}+\eta_{4})
/4
,
\\
a_{5} & =
%\frac{1}{4}%
(\eta_{0}-\eta_{1}+\eta_{2}+\eta_{3}-\eta_{4})
/4
,
\end{aligned}
\label{eqA.11}%
\end{equation}
where
\begin{equation}
\begin{aligned}
\eta_{0} & =
%\frac{1}{\sqrt{2}}
1/\sqrt{2}, & & \\
\eta_{1} & =
%\frac{1}{\sqrt{2}}
(\cos2\theta+\cos2\rho)
/\sqrt{2}, \quad 
& \eta_{2} & =
%\frac{1}{\sqrt{2}}
(\sin2\theta+\sin2\rho)
/\sqrt{2}, \\
\eta_{3} & =
%\frac{1}{\sqrt{2}}
\cos(2\theta-2\rho)
/\sqrt{2},
& \eta_{4} & = 
%\frac{1}{\sqrt{2}}
\sin(2\theta-2\rho)
/\sqrt{2}.
\end{aligned}
\label{eqA.11b}%
\end{equation}
These formulas arise from an independent pair of rotations by angles $\theta$
and $\rho$ of two ``spin vectors%
\index{spin vector}'', i.e., by taking the matrix%
\index{matrix!function@--- function} function $A$ in (\ref{eqIntFeb22.11}) unitary%
\index{matrix!unitary}, $\mathbb{T}\ni
z\rightarrow A_{\theta,\rho}\left(  z\right)  \in\mathrm{U}_{2}\left(
\mathbb{C}\right)  $, and setting
\begin{equation}
A(z)=V(Q_{\theta}^{\perp}+zQ_{\theta}^{{}})(Q_{\rho}^{\perp}+zQ_{\rho}^{{}})
=VU_{\theta}(z)U_{\rho}(z)
\label{eqA.8}%
\end{equation}
with%
\begin{align}
V&=
%\smash[t]{
\frac{1}{\sqrt{2}}\left(
\begin{array}
[c]{cc}%
1 & 1\\
1 & -1
\end{array}
\right)  
%}
, \label{eqA.8bis}%
\\[3\jot]
Q_{\theta}&=\left(
\begin{array}
[c]{cc}%
\cos^{2}\theta & \cos\theta\sin\theta\\
\cos\theta\sin\theta & \sin^{2}\theta
\end{array}
\right)  \nonumber \\
&=\frac{1}{2}\left(  \left(
\begin{array}
[c]{cc}%
1 & 0\\
0 & 1
\end{array}
\right)  +\left(
\begin{array}
[c]{cc}%
\cos2\theta & \sin2\theta\\
\sin2\theta & -\cos2\theta
\end{array}
\right)  \right)  , \label{eqA.9}%
\end{align}
and the orthogonal complement to the one-dimensional projection%
\index{operator!projection} $Q_{\theta}$,%
\begin{equation}
Q_{\theta}^{\perp}=\rlap{$\displaystyle 
Q_{\theta+\left(  \pi/2\right)  }^{{}}.$}
\phantom{\frac{1}{2}\left(  \left(
\begin{array}
[c]{cc}%
1 & 0\\
0 & 1
\end{array}
\right)  +\left(
\begin{array}
[c]{cc}%
\cos2\theta & \sin2\theta\\
\sin2\theta & -\cos2\theta
\end{array}
\right)  \right)  ,}
\label{eqA.10}%
\end{equation}

With the coefficients $a_{0}$, $a_{1}$, $a_{2}$, $a_{3}$, $a_{4}$, $a_{5}$
given by (\ref{eqA.11}), the algorithmic approach
to graphing the solution $\varphi$ to the
scaling%
\index{scaling!identity@--- identity} identity (\ref{eqInt.1}) is as
follows (see \cite{Jor01b}, \cite{Tre01b} for details):
the relation (\ref{eqInt.1}) for $N=2$
is interpreted as giving the values of the left-hand $\varphi$
by an operation performed on those of the $\varphi$ on the right,
and a binary digit inversion transforms this into the form
\begin{equation}
\mathbf{f}_{k+1}^{\prime}\left( x+\frac{1}{2^{k+1}}\right) 
=\mathbf{A}\mathbf{f}_{k}\left( x\right) ,
\label{eqCascadeAlgorithm}
\end{equation}
where $\mathbf{A}$
is the $2\times 3$ matrix $\mathbf{A}_{i,j}
=\sqrt{2}a_{4+i-2j}$ constructed from the
coefficients in (\ref{eqInt.1}),
and
$\mathbf{f}_{j}$ and $\mathbf{f}_{j}^{\prime}$
are the vector functions
\begin{equation}
\mathbf{f}_{j}\left( x\right) 
=\left( 
\begin{matrix}
\varphi\left( x-\frac{2}{2^{j}}\right) \\
\varphi\left( x-\frac{1}{2^{j}}\right) \\
\varphi\left( x\right) 
\end{matrix}
\right) ,
\qquad
\mathbf{f}_{j}^{\prime}\left( x\right) 
=\left( 
\begin{matrix}
\varphi\left( x-\frac{1}{2^{j}}\right) \\
\varphi\left( x\right) 
\end{matrix}
\right) .
\label{eqCascadeVector}
\end{equation}

The signal processing aspect can be understood from the description of subband
filters in the analusis and synthesis of time signals, or more general signals
for images. In either case, we have two subband systems $m=(m_{0},m_{1}%
,\dots)$ and $\tilde{m}=(\tilde{m}_{0},\tilde{m}_{1},\dots)$ where the
functions%
\[
m_{j}\left(  z\right)  =\sum_{n}a_{n}^{\left(  j\right)  }z^{n}\text{\quad
and\quad}\tilde{m}_{j}\left(  z\right)  =\sum_{n}\tilde{a}_{n}^{\left(
j\right)  }z^{n}%
\]
are the generating functions defined from the filter coefficients $\smash{\left(
a_{n}^{\left(  j\right)  }\right)  }$ and $\left(  \tilde{a}_{n}^{\left(
j\right)  }\right)  $, $n\in\mathbb{Z}^{d}$.

\section*{\label{App}Appendix A: Topics for further research} 

Originally we had anticipated adding two more chapters to these 
tutorials, but time and space prevented this. Instead we include the table of 
contents for this additional material. The details for the remaining chapters 
will be published elsewhere. But as the items in the list of contents suggest, 
there are still many exciting open problems in the subject that the reader may 
wish to pursue on his/her own.
We feel that the following list of topics offers at least an outline of several 
directions that the reader, could take in his/her own study and research on 
wavelet-related mathematics.

\section{\label{S3}Connection between the discrete signals and the wavelets}

\subsection{\label{S3.1}Wavelet geometry in $L^2(R^n)$}

\subsection{\label{S3.2}Intertwining operators between sequence spaces $l^2$ and $L^2(R^n)$}

\subsection{\label{S3.3}Infinite products of matrix functions}

\subsubsection{\label{S3.3.1}Implications for $L^2(R^n)$}

\subsubsection{\label{S3.3.2}Wavelets in other Hilbert spaces of fractal measures}

\subsection{\label{S3.4}Dependence of the wavelet functions on the matrix functions
which define the wavelet filters}

\subsubsection{\label{S3.4.1}Cycles}

\subsubsection{\label{S3.4.2}The Ruelle-Lawton wavelet transfer operator}

%\newpage

\section{\label{S4}Other topics in wavelets theory}

\subsection{\label{S4.1}Invariants}

\subsubsection{\label{S4.1.1}Invariants for wavelets: Global theory}

\subsubsection{\label{S4.1.2}Invariants for wavelet filters: Local theory}

\subsection{\label{S4.2}Function classes}

\subsubsection{\label{S4.2.1}Function classes for wavelets}

\subsubsection{\label{S4.2.2}Function classes for filters}

\subsection{\label{S4.3}Wavelet sets}

\subsection{\label{S4.4}Spectral pairs}
\medskip
%\newpage

\section*{\label{Dual}Appendix B: Duality principles in analysis}

Several versions of spectral duality are presented. On the two sides we
present (1) a basis condition, with the basis functions indexed by a
frequency variable, and giving an orthonormal basis; and (2) a geometric
notion which takes the form of a tiling, or a Iterated Function System
(IFS). Our initial motivation derives from the Fuglede conjecture, see \cite{B2,B6,B7}: For a
subset $D$ of $\mathbb{R}^n$ of finite positive measure, the Hilbert space $L^2(D)$ admits
an orthonormal basis of complex exponentials, i.e., $D$ admits a Fourier basis
with some frequencies $L$ from $\mathbb{R}^n$, if and only if $D$ tiles $\mathbb{R}^n$ (in the
measurable category) where the tiling uses only a set $T$ of vectors in $\mathbb{R}^n$.
If some $D$ has a Fourier basis indexed by a set $L$, we say that $(D,L)$ is a spectral pair.
We recall from \cite{B8} that if $D$ is an $n$-cube, then the sets $L$ in (1) are precisely the
sets $T$ in (2). This begins with work of Jorgensen and Steen Pedersen \cite{B8}
where the admissible sets $L = T$ are characterized.
Later it was shown, \cite{B5} and \cite{B10} that the identity $T = L$ holds
for all $n$. The proofs are based on general Fourier duality, but they do not
reveal the nature of this common set $L = T$. A complete list is known only
for $n=1$, $2$, and $3$, see \cite{B8}.

We then turn to the scaling IFS's built from the $n$-cube with a given
expansive integral matrix $A$. Each $A$ gives rise to a fractal in the small,
and a dual discrete iteration in the large. In a different paper \cite{JoPe98},
Jorgensen and Pedersen characterize those IFS fractal limits which admit
Fourier duality. The surprise is that there is a rich class of fractals that
do have Fourier duality, but the middle third Cantor set does not.
We say that an affine IFS, built on affine maps in $\mathbb{R}^n$ defined by a given
expansive integral matrix $A$ and a finite set of translation vectors, admits
Fourier duality if the set of points $L$, arising from the iteration of the
$A$-affine maps in the large, forms an orthonormal Fourier basis (ONB) for the
corresponding fractal $\mu$ in the small, i.e., for the iteration limit built
using the inverse contractive maps, i.e.,  iterations of the dual affine
system on the inverse matrix $A^{-1}$. By ``fractal in the small'', we mean the
Hutchinson measure $\mu$ and its compact support, see \cite{Hut81}. (The best known
example of this is the middle-third Cantor set, and the measure $\mu$ whose
distribution function is corresponding Devil's staircase.)

In other words, the condition is that the complex exponentials indexed by $L$
form an ONB for $L^2(\mu)$. Such duality systems are indexed by complex
Hadamard matrices $H$, see \cite{B8} and \cite{JoPe98}; and the duality issue is connected to
the spectral theory of an associated Ruelle transfer operator, see \cite{BrJo02b}.
These matrices $H$ are the same Hadamard matrices which index a
certain family of quasiperiodic  spectral pairs $(D,L)$ studied in \cite{B6} and
\cite{B7}. They also are used in a recent construction of Terence Tao \cite{Tao} of a
Euclidean spectral pair $(D,L)$ in $\mathbb{R}^5$ for which $D$ does not a tile $\mathbb{R}^5$ with
any set of translation vectors $T$ in $\mathbb{R}^5$; see also \cite{IKT03}.

We finally report on joint research with Dorin Dutkay \cite{DuJo03}, \cite{DuJo04a}, \cite{DuJo04b}, \cite{DuJo04c}
where we show that all the affine IFS's, and more general limit systems 
from dynamics and probability theory,  admit wavelet constructions, i.e., admit 
orthonormal bases of wavelet functions in Hilbert spaces which are constructed 
directly from the geometric data. A substantial part of the picture involves 
the construction of limit sets and limit measures, a part of geometric measure 
theory.
%\bigskip
%\newpage

\noindent\textbf{Acknowledgements:}\enspace   
\label{Ack}We are happy to thank the organizing committee at the
National University of Singapore for all their dedicated work in planning
and organizing a successful conference, of which this tutorial is a part. We
especially thank Professors Wai Shing Tang, Judith Packer, Zuowei Shen, and
the head of the Department of Mathematics of the NUS for all their work in
making my visit to Singapore possible. We thank the Institute for Mathematical Sciences at the
NUS, and the US National Science Foundation, for partial financial support
in the preparation of these lecture notes. We discussed various parts of the
mathematics with our colleagues, Professors Larry Baggett, David Larson, Ola
Bratteli, Kathy Merrill, Judy Packer, and we thank them for their
encouragements and suggestions. The typesetting and graphics were expertly
done at the University of Iowa by Brian Treadway. We also thank Brian
Treadway for a number of corrections, and for very helpful suggestions,
leading to improvements of the presentation.

\ifx\undefined\bysame
\newcommand{\bysame}{\leavevmode\hbox to3em{\hrulefill}\,}
\fi

\printindex

\end{document}